\documentclass[10pt]{article}

\usepackage{amssymb}
\usepackage{fullpage}

\usepackage{graphicx}        
\usepackage{amscd}

\usepackage{amssymb}
\usepackage{amsmath}
\usepackage{mathrsfs}
\usepackage{verbatim}
\usepackage{proof}

\graphicspath{{./Images/}}

\newtheorem{theorem}{Theorem}
\newtheorem{lemma}[theorem]{Lemma}
\newtheorem{definition}[theorem]{Definition}
\newtheorem{corollary}[theorem]{Corollary}

\newtheorem{remark}[theorem]{Remark}

\DeclareMathAlphabet\gothic{U}{euf}{m}{n}

\usepackage{color}

\usepackage{caption}
\usepackage{subcaption}

\begin{document}




\title{Left-invariant evolutions of wavelet transforms on the Similitude Group}


\author{Upanshu Sharma and Remco Duits 
}
\maketitle

\begin{abstract}
Enhancement of multiple-scale elongated structures in noisy image data is relevant for many biomedical applications but commonly used PDE-based enhancement techniques often fail at crossings in an image. To get an overview of how an image is composed of local multiple-scale elongated structures we construct a multiple scale orientation score, which is a continuous wavelet transform on the similitude group, $SIM(2)$. Our unitary transform maps the space of images onto a reproducing kernel space defined on $SIM(2)$, allowing us to robustly relate Euclidean (and scaling) invariant operators on images to left-invariant operators on the corresponding continuous wavelet transform. 
Rather than often used wavelet (soft-)thresholding techniques, we employ the group structure in the wavelet domain to arrive at left-invariant evolutions and flows (diffusion), for contextual crossing preserving enhancement of multiple scale elongated structures in noisy images. We present experiments that display benefits of our work compared to recent PDE techniques acting directly on the images and to our previous work on left-invariant diffusions on orientation scores defined on Euclidean motion group.
\end{abstract}







\section{Introduction \label{ch:1}}
Elongated structures in the human body such as fibres and blood vessels often require analysis for diagnostic purposes. A wide variety of medical imaging techniques such as magnetic
resonance imaging (MRI), microscopy, X-ray 
fluoroscopy, fundus imaging etc. exist to achieve this.
Many (bio)medical questions related to such images require detection and tracking of the
elongated structures present therein. Due to the desire to reduce acquisition time and radiation dosage the acquired  medical images are often noisy, of low contrast and suffer from occlusions and incomplete data.  Furthermore multiple-scale elongated structures exhibit crossings and bifurcations which is a notorious problem in (medical) imaging. Hence crossing-preserving enhancement of these structures 
is an important preprocessing step for subsequent detection.
 
In recent years PDE based techniques have gained popularity in the field of image processing. Due to well posed mathematical results these techniques lend themselves to stable algorithms and also allow mathematical and geometrical interpretation of classical methods such as Gaussian and morphological filtering, dilation or erosion etc. on $\mathbb{R}^d$. 
These techniques typically regard the original image, $f\in\mathbb{R}^2\rightarrow\mathbb{R}$, as an initial state of a parabolic (diffusion like) evolution process yielding  filtered versions, $u_f:\mathbb{R}^2\times\mathbb{R}^+\rightarrow\mathbb{R}$. Here $u_f$ is called the scale space representation of image $f$. The domain of $u_f$ is scale space $\mathbb{R}^2\times\mathbb{R}^+$.
A typical scale space evolution is of the form
\begin{align}
\begin{cases}
\partial_{s}u_f (\mathbf{x},s)&=\nabla_{\mathbf{x}}\cdot(C(
u_{f}(\cdot,s))(\mathbf{x})\nabla_{\mathbf{x}}u_f)(\mathbf{x},s)\\
u_{f}(\mathbf{x},0)&=f(\mathbf{x}),
\end{cases}\label{Gen_ScaleSpace}
\end{align}
where $C(u_{f}(\cdot,s))(\mathbf{x})$ models the diffusivity depending on the differential structure at $(\mathbf{x},s,u_{f}(\mathbf{x},s))\in\mathbb{R}^d\times\mathbb{R}^+\times\mathbb{R}$.
For $C=1$, Eq.\eqref{Gen_ScaleSpace} is the usual linear heat equation. The corresponding evolution is  known in image processing as a Gaussian Scale Space \cite{Lindeberg1993,Alvarez1993,Romeny1997,Lindeberg2013}. 
In their seminal paper \cite{Perona1990}, Perona and Malik proposed nonlinear filters to bridge scale space and restoration ideas. Based on the observation that diffusion should not occur when the (local) gradient value is large 
(to avoid blurring the edges), they pointed out that nonlinear adaptive isotropic diffusion is achieved by replacing $C=1$ by $C(u_{f}(\cdot,s))(\mathbf{x})=c(\|\nabla_{\mathbf{x}}u_{f}
(\mathbf{x},s)\|)$, where $c:\mathbb{R}^+\rightarrow\mathbb{R}^+$ is some smooth strictly decreasing positive function 
vanishing at infinity. An improvement of the Perona-Malik scheme is the ``coherence-enhancing diffusion" (CED) introduced by Weickert \cite{Weickert1999} which additionally uses the direction of the gradient $\nabla_{\mathbf{x}}u_{f}$ leading to diffusion constant $c$ being replaced by a nonlinear matrix. 

However these methods often fail in image analysis applications with crossing or bifurcating curves as the direction of gradient at these structures is ill-defined, see \cite{Franken2009} for more details. Scharr et al. in \cite{Scharr2012} present techniques which effectively deal with the particular case of X-junctions by relying on the $2$-nd order jet of Gaussian derivatives in the image domain.  Passing through higher order jets of Gaussian derivatives and induced Euclidean invariant differential operators does not allow one to generically  deal with complex crossings and/or bifurcating structures. 
Instead we need gauge frames in higher dimensional affine Lie groups to deal with this issue. A gauge frame is a local coordinate system aligned/gauged with locally present (elongated) structures in an image. Differentiating w.r.t. such coordinates provides intrinsically natural derivatives as opposed to differentiating w.r.t. (artificially imposed) global coordinates. At interesting locations in the image, where multiple scale elongated structures cross, one needs multiple gauge frames. Therefore instead of gauge frames per position, $x\in\mathbb{R}^2$ in a (grey-scale) image $f:\mathbb{R}^2\rightarrow\mathbb{R}$, we attach gauge frames to each Lie group element,
\begin{align*}
g=(x,t)\in G=\mathbb{R}^2\rtimes T.
\end{align*} 
in a Coherent state (CS) transform $\mathcal{W}_\psi f:G\rightarrow \mathbb{C}$ of an image $f:\mathbb{R}^2\rightarrow\mathbb{R}$.
In this article we mainly consider $(G=SE(2),T=SO(2))$ and $(G=SIM(2),T=\mathbb{R}^+\times SO(2))$, where (multiple scale) elongated structures are manifestly disentangled within the score, allowing for a crossing preserving flow (steered by gauge frames in the score). 
In medical image processing $\mathcal{W}_\psi f$ for $G=SE(2)$ is also referred to as an orientation score as it provides a score of how an image is decomposed out of local (possibly crossing) orientations.  
\subsection{Why extend to the $SIM(2)$ group?}
In this paper we wish to extend the aforementioned orientation score framework to the case of the similitude group 
(group of planar translations, rotations and scaling),  for the following reasons:
\begin{itemize}
\item{Gauge frames based on Gaussian derivatives are usually aligned with a \emph{single} Gaussian gradient direction, which is unstable in the vicinity of complex structures such as crossings and bifurcations. At these structures one needs multiple spatial frames per position. Therefore following the general idea of Scale spaces on affine Lie groups \cite{Duits2007b,Citti2006,Duits2012,Duits2005},  a natural next step would be to extend the domain of images to the affine Lie group $SIM(2)=\mathbb{R}^2\rtimes (\mathbb{R}^+\times SO(2))$ in order to have well posed gauge frames at crossing structures (adapted to multiple orientations and scales that are locally present).}
\item{In the primary visual cortex both multiple scales and orientations are encoded per position. It is generally believed that receptive field profiles in neurophysiological experiments can be modelled by Gaussian derivatives \cite{DeAngelis1995,Young1987,Landy1991,
Bosking1997} and  this provides a biological motivation to incorporate scales.}
\item{Elongated (possibly crossing) structures often exhibit multiple scales. Earlier work by one of the authors \cite{Duits2007a,Franken2009} proposes generic crossing preserving flows via invertible orientation scores (without explicit multiple scale decomposition). However, this approach treats all scales in the same way. As a result these flows do not always  adequately deal with images containing elongated structures with strongly varying widths (scales). Therefore we must encode and process multiple scales in the scores.}
\end{itemize}
Our framework involves the construction of a continuous wavelet transform of an image using the similitude group. This approach is related to that of directional wavelets \cite{Antoine1996,Antoine1999,
Antoine2004,Jacques2011,Ali2014} and  curvelets \cite{Donoho1999, Donoho2005,Donoho2005a,Candes2006}. In particular the recent popular approach of shearlets \cite{Labate2005,
Guo2007,Dahlke2011,Bodmann2013} also falls in the category of continuous wavelet transform but uses the shearlet group which includes the shearing, translation and scaling group. The main aim of this article is to construct rotation and translation invariant diffusion type flows on the wavelet transformed image. Although highly useful, the shearlet group (which is related to the similitiude group by nilpotent approximations, see \ref{App:MethContraction}) is not very suited for such exact rotation and translation covariant processing because of its group structure.

\subsection{Our main results}
\begin{figure}
  \begin{minipage}[c]{0.5\textwidth}
    \includegraphics[width=\textwidth]{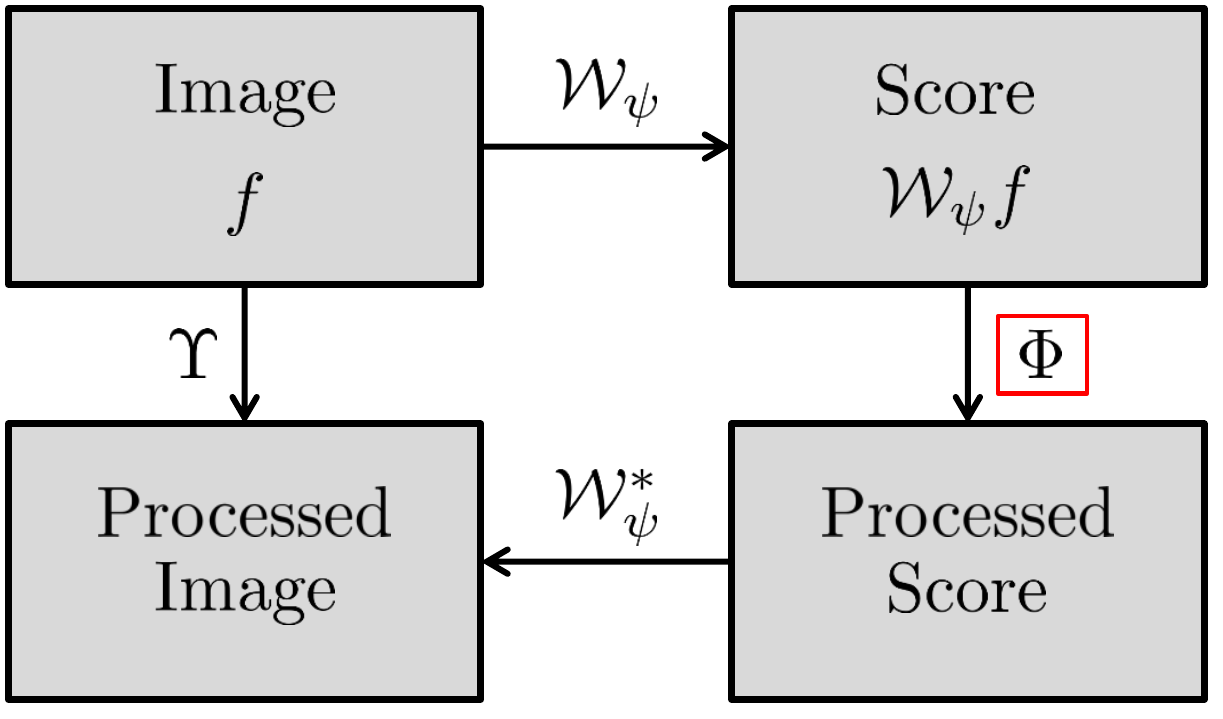}
  \end{minipage} \hspace{4mm}
  \begin{minipage}[c]{0.40\textwidth}
    \caption{
    A schematic view on processing images $f:\mathbb{R}^2\rightarrow\mathbb{R}$ via CS transform $\mathcal{W}_\psi f:G\rightarrow \mathbb{C}$ defined on Lie group $G=\mathbb{R}^2\rtimes T$. Design of well posed transform $\mathcal{W}_\psi$ and of appropriate operators $\Phi$ is the main objective of this article. Note that $\Upsilon=\mathcal{W}_\psi^*\circ\Phi\circ\mathcal{W}_\psi$, i.e. spatial processes are realized via invertible transforms akin to cortical columns in the visual brain \cite{Hubel1993}. 
    } \label{fig1}
  \end{minipage}
\end{figure}

\begin{figure}[t]
\centering
\includegraphics[scale=0.65]{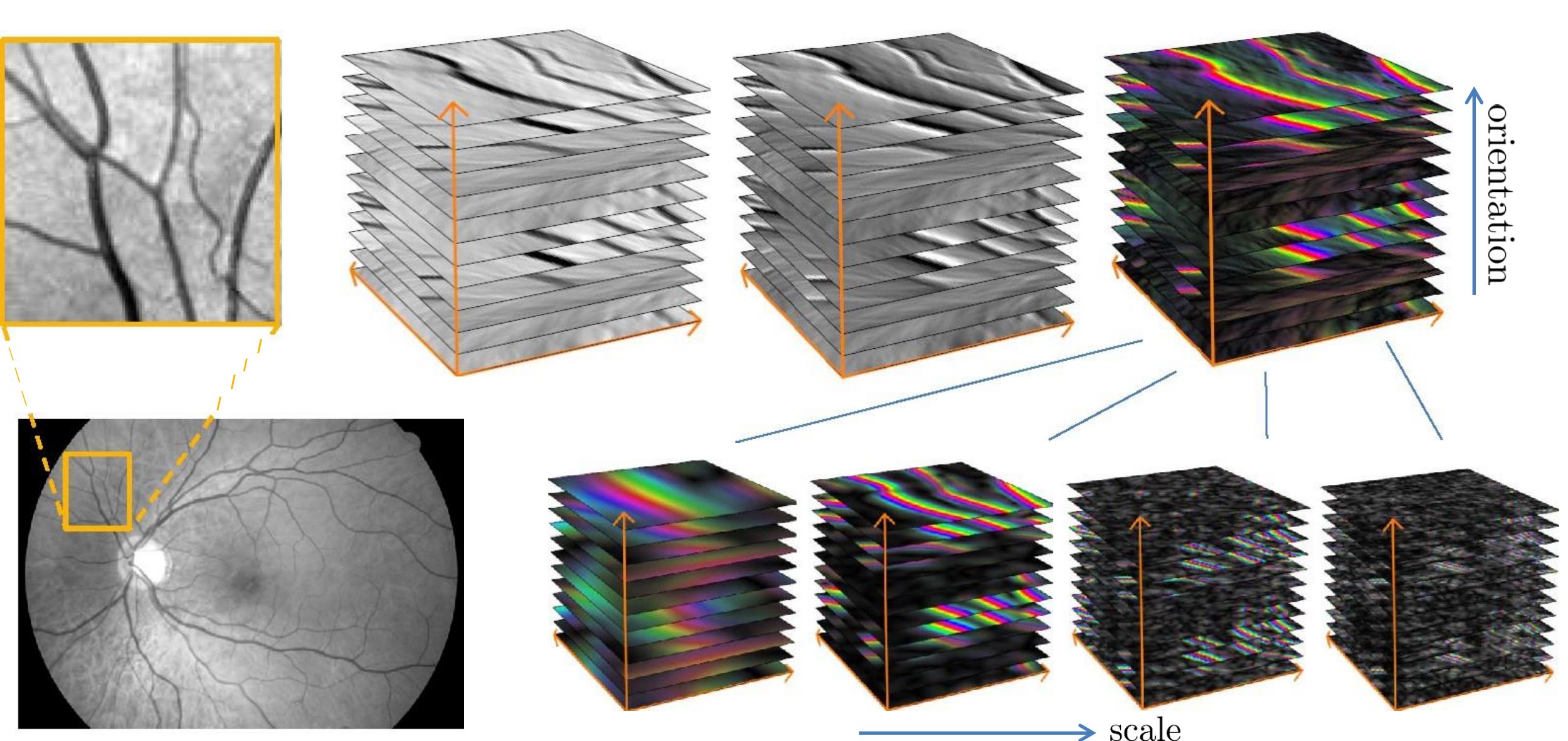}
\caption[ScaleOS]{Top row: original image, the real-part of the CS transform on $SE(2)$ group reflects the centerlines, the imaginary
part the edges of the bloodvessels, the combination (color represents phase direction
and intensity represents the absolute value). Bottom row: visualizations of CW Transform on $SIM(2)$ group that allow us to include scale adaptation in our enhancement and detection.}
\label{fig:ScaleOSIntroIm}
\end{figure}

There are two main motivating questions (see Figure~\ref{fig1}) for the work presented in this article. 
\begin{enumerate}
\item{Can we design  a well-posed invertible score,  which is a complex-valued function on a Lie group $SIM(2)$,  combining the strengths of directional wavelets \cite{Antoine1996,Antoine1999, Antoine2004}, curvelets \cite{Donoho1999, Donoho2005,Donoho2005a,Candes2006} in such a way that allows for accurate and efficient implementation of subsequent contextual-enhancement operators?}
\item{Can we construct contextual flows in the wavelet domain, in order to ensure that only the wavelet coefficients that are coherent (from both probabilistic and group theoretical perspective, \cite{Ali1998}) with the surrounding coefficients become dominant?}
\end{enumerate}
Our answer to both these questions is indeed affirmative. 
In the first part we follow the general approach by Grossmann, Morlet and Paul in \cite{Grossmann1985}, Ali in \cite{Ali1998} and F{\"u}hr in \cite{Fuhr2005} we provide a short review of general results in the context of our case of interest.
In practice there are upper/lower bounds on scaling and so in Theorem~\ref{Thm:StabCond} and we provide stability analysis and condition number of a modified continuous wavelet transform on $SIM(2)$ which can be used for practical applications. The crucial point in our design is that we rely on explicit B-spline decomposition along log-polar type of coordinates in the Fourier domain. These are the canonical coordinates of the second type for $SIM(2)$, which are required for accurate left-invariant processing on the CS transformed images (scores). Figure~\ref{fig:ScaleOSIntroIm} depicts the use of this wavelet transform in practice.\\

The latter and main part of this article is dedicated to answering the second question. Theorem~\ref{thm:EuclideanInv} shows that only left-invariant operators on the transformed image correspond to Euclidean (and scaling) invariant operators on images.  Therefore in this article we restrict ourselves to left-invariant PDE evolutions on $SIM(2)$. Theorem~\ref{thm:StochasticConn} provides a stochastic connection to our left-invariant flows in the wavelet domain. These flows are forward Kolmogorov equations corresponding to stochastic processes for multiple-scale contour enhancement on $SIM(2)$. Using the general theory of coercive operators on Lie groups, in Eq.\eqref{eq:Gaussian_Est}  Gaussian estimates for the Green's function of linear diffusion on $SIM(2)$ are derived. In Eq.\eqref{NonLinDiffOp} we present nonlinear left-invariant adaptive diffusions on the image of a wavelet transform. Finally we present experiments to validate its practical advantages which become clear after comparisions with current state of the art enhancements. 
\subsection{Invertible Orientation Scores}
Based on the early work by Ali, Antoine and Gazeau \cite{Ali2014} and Kalitzin \cite{Kalitzin1999}, Duits et al. \cite{Duits2007a,Duits2005} introduced the framework of invertible orientation score in medical imaging to effectively handle the problem of generic crossing curves in the context of bio-medical applications. In this subsection we briefly explain the ideas developed in \cite{Duits2007a,Duits2005}.

The $2D$-Euclidean motion group (i.e.\ the group of planar rotations and translations) $SE(2)$ is defined as $SE(2)=\mathbb{R}^2\rtimes SO(2)$ where $SO(2)$ is the group of planar rotations. A CS transform, $U_{f}:SE(2)\rightarrow\mathbb{C}$ of an image $f:\mathbb{R}^2\rightarrow\mathbb{R}$ is obtained by means of an anisotropic convolution kernel $\check{\psi}:\mathbb{R}^2\rightarrow\mathbb{C}$ via
\begin{align*}
U_{f}(g)=\int\limits_{\mathbb{R}^2} \overline{\check{\psi}(\mathbf{R}_{\theta}^{-1} (\mathbf{y}-\mathbf{x}))}f(\mathbf{y})
d\mathbf{y}, \ g=(\mathbf{x},\theta)\in SE(2),
\end{align*}
where $\psi(-\mathbf{x})=\check{\psi}(\mathbf{x})$ and $\mathbf{R}_{\theta}\in SO(2)$ is the $2$D counter-clockwise rotation by angle $\theta\in[0,2\pi)$. Assume $\psi\in\mathbb{L}_2(\mathbb{R}^2)$, then the transform $\mathcal{W}_{\psi}$ which maps images $f\in\mathbb{L}_{2}(\mathbb{R}^2)$ can be rewritten as 
\begin{align*}
U_{f}(g)=(\mathcal{W}_{\psi}f)(g)=(\mathcal{U}_{g}\psi,f)_{\mathbb{L}_2
(\mathbb{R}^2)},
\end{align*}
where $g\mapsto\mathcal{U}_{g}$ is a unitary (group-)representation of the Euclidean motion group $SE(2)$ into $\mathbb{L}_2
(\mathbb{R}^2)$ given by $\mathcal{U}_{g}f(\mathbf{y})=f(\mathbf{R}
_{\theta}^{-1}(\mathbf{y}-\mathbf{x}))$ for all $g=(\mathbf{x},\mathbf{R}_\theta)\in SE(2)$ and for all $f\in\mathbb{L}_{2}(\mathbb{R}^2)$.
It is constructed by means of an admissible vector $\psi\in\mathbb{L}(\mathbb{R}^2)$ such that $\mathcal{W}_{\psi}$ is unitary onto the unique reproducing kernel Hilbert space $\mathbb{C}_{K}^{SE(2)}$ of functions on $SE(2)$ with reproducing kernel $K(g,h)=(\mathcal{U}_{g}\psi,\mathcal{U}_{h}\psi)$, which is a closed vector subspace of $\mathbb{L}_{2}(SE(2))$. This leads to the essential Plancherel formula (see \cite{Ali1998,Fuhr2005})
\begin{align*}
\|\mathcal{W}_{\psi}f\| ^{2}_{\mathbb{C}_{K}^{SE(2)}}=\int\limits
_{\mathbb{R}^2}\int\limits_{0}^{2\pi}|(\mathcal{F}&\mathcal{W}_{\psi}f)
(\boldsymbol{\omega},\theta)|^2 \frac{1}{M_{\psi}(\boldsymbol{\omega})}d\boldsymbol{\omega}
d\theta 
=\int\limits_{\mathbb{R}^2}
|(\mathcal{F}f)(\boldsymbol{\omega})|^2 d\boldsymbol{\omega}=\|f\|^2_{\mathbb{L}_2(\mathbb{R}^2)},
\end{align*}
where $M_{\psi}\in C(\mathbb{R}^2,
\mathbb{R})$ is given by $M_{\psi}
(\boldsymbol{\omega})=\int\limits_{0}^{2\pi}
|\mathcal{F}\psi(\mathbf{R}_{\theta}^{T}
\boldsymbol{\omega})|^2 d\theta$. If $\psi$ is chosen such that $M_{\psi}=1$ then we gain $\mathbb{L}_2$ norm preservation. But this is not possible as $\psi\in\mathbb{L}_2(\mathbb{R}^2)\cap
\mathbb{L}_1(\mathbb{R}^2)$  implies that $M_{\psi}$ is a continuous function vanishing at infinity. In practice, however, because of finite grid sampling,  $\mathcal{U}$ is restricted to the space of disc limited images\footnote{This requirement can be avoided in general by using distributional coherent transforms, for e.g. see \cite[App.B]{Bekkers2012}.}.
\begin{figure}[t]
\centering
\includegraphics[scale=0.7]{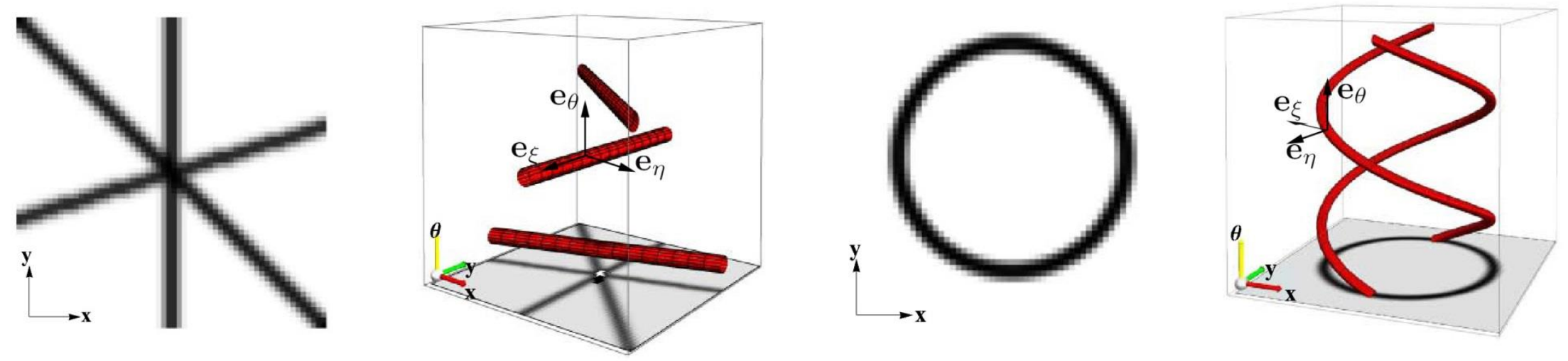}
\caption[Intuitive illustration of a orientation score]{Illustration of a CS transform on $SE(2)$ for a crossing and a circle. This framework unwraps a crossing in the CS domain. In the case of the circle, the resulting response in a $3$D visualization is a spiral since the orientation changes linearly as one traverses the circle. Note that the orientation dimension (displayed vertically) is $2\pi$ periodic.
}
\label{fig:IntuitiveOSCircle}
\end{figure}
Since the CS transform $\mathcal{W}_{\psi}$ maps the space of images $\mathbb{L}_{2}(\mathbb{R}^2)$ unitarily onto the space  $\mathbb{C}_{K}^{SE(2)}$ (provided $M_{\psi}>0$) the original image $f:\mathbb{R}^2\rightarrow\mathbb{R}$ can be reconstructed from  $U_f$ via the adjoint CS transform, $f=\mathcal{W}_{\psi}^*\mathcal{W}_{\psi}[f]$.
\\

For examples of wavelets $\psi$ for which $M_{\psi}=1\large{|}_{B_{0,\varrho}}$ and details on fast approximative reconstruction by integration over angles only,
see \cite{Duits2005}. For details on image processing (particularly enhancement and completion of crossing elongated structures) via CS transform on $SE(2)$ (orientation scores), see \cite{Duits2008,Franken2009,Duits2010,
Duits2010a}. An intuitive illustration of the relation between an elongated structure in a $2$D image and the corresponding elongated structure in CS domain is given in Figure~\ref{fig:IntuitiveOSCircle}. In the top row of Figure~\ref{fig:ScaleOSIntroIm} we also depict a practical example of this transform. 

\subsection{Structure of the article} 
This article is structured as follows.
\begin{itemize}
\item{{\bf (Section 2) Unitary CS transform on images:} Stability of the CS transform is discussed followed by an explicit construction of so called proper wavelets that allow for both a stable (re)construction of the transformed image and accurate implementation of subsequent left-invariant flows.} 
 \item{ {\bf (Section 3) Operators on scores:} Employing the group structure in the wavelet domain a general framework for operators on the scores, involving left-invariant evolutions is discussed. These operators are interpreted in a differential geometric (and probabilistic) setting to provide a strong intuitive rationale for their choice.} 
\item {{\bf (Section 4) Left-invariant diffusions:} Gaussian estimates for the Green's function of linear diffusion on $SIM(2)$ are derived followed by a discussion of non-linear adaptive diffusion.}
\item {{\bf (Section 5) Practical results:} In this section we present experiments that show the advantages of adaptive non-linear diffusion on multi-scale orientation scores in comparison to PDE techniques and previous work on orientation scores. Furthermore we compare our technique to state of the art denoising algorithms. Finally, we also apply these evolutions on well-established vesselness techniques and show benefits.
}
\end{itemize}
\section{Unitary operators between images and scores \label{ch:Hr}}
In this section we present a quick overview of abstract coherent state transforms from a Hilbert space to a functional Hilbert space and thereby arrive at a continuous wavelet transform in our case of interest, the $SIM(2)$ group. This is followed by quantifying the stability of this transform in a sense which will be made precise. We end the section by explicitly constructing so called proper wavelets which would allow us to create appropriate PDE flows in the subsequent sections.

The continuous wavelet transform constructed by unitary irreducible representations of locally compact groups was first formulated by Grossman et al.~\cite{Grossmann1985}. Given a Hilbert space $H$ and a unitary irreducible representation $g\mapsto \mathcal{U}_g$ of any locally compact group $G$ in $H$, a nonzero vector $\psi\in H$ is called admissible if
\begin{align}\label{CpsiDef}
C_{\psi}:= \int\limits_{G}\frac{|(\mathcal{U}_{g}\psi,\psi)|^2}{(\psi,\psi)_{H}}d\mu_{G}(g)<\infty,
\end{align}  
where $\mu_G$ denotes the left invariant Haar measure on $G$. Given an admissible vector $\psi$ and a unitary representation of a locally compact group $G$ in a Hilbert space $H$, the CS transform $\widetilde{\mathcal{W}}_\psi:H\rightarrow \mathbb{L}_{2}(G)$ is defined by 
\begin{align*}
(\widetilde{\mathcal{W}}_\psi[f])(g)=(\mathcal{U}_{g}\psi,f)_H.
\end{align*}
It is well known in mathematical physics~\cite{Ali1998}, that $\widetilde{\mathcal{W}}_\psi$ is a unitary transform onto a closed reproducing kernel space $\mathbb{C}_{K_{\psi}}^{G}$ with $K_{\psi}(g,g')=\frac{1}{C_\psi}(\mathcal{U}_g\psi,\mathcal{U}_{g'}\psi)_{H}$ as an $\mathbb{L}_2$-subspace.
Note that we distinguish between the unitary wavelet transform $\widetilde{\mathcal{W}}:\mathbb{L}_2(\mathbb{R}^2)\rightarrow\mathbb{L}_2(G)$ and the isometric wavelet transform $\mathcal{W}:\mathbb{L}_2(\mathbb{R}^2)\rightarrow \mathbb{C}_K^G$ since the orthogonal projection onto $\mathbb{C}_K^G$ is given by $\mathbb{P}_\psi=\widetilde{\mathcal{W}}_\psi\widetilde{\mathcal{W}}^*_\psi$ whereas $\mathcal{W}_\psi\mathcal{W}_\psi^*=I$.

\subsection{Construction of a Unitary Map from $H$ to $\mathbb{C}^\mathbb{G}_K$}
When $G$ is a  locally compact group the norm on $\mathbb{C}^\mathbb{G}_K$ has a simpler explicit form compared to  usual description, see \cite[Ch:7.2]{Duits2005} and \cite[Lemma 1.7]{Martens1988}. One can give an explicit characterization of $\mathbb{C}^\mathbb{G}_K$, in the case $G=\mathbb{R}^d\rtimes_\varphi T$, with $T$  a linear algebraic group (for definition see \cite{Borel1991}), $\mathcal{U}$ the left-regular action of $G$ onto $H=\mathbb{L}_2(\mathbb{R}^2)$ and thereby formulate a reconstruction theorem for affine groups.
We call $\psi\in\mathbb{L}_2(\mathbb{R}^d)$ an admissible vector if
\begin{align}
0<M_\psi(\omega):=(2\pi)\int\limits_{T}\left|\frac{\mathcal{F}[\mathcal{R}_t\psi(\omega)]}{\sqrt{\text{det} (\varphi(t))}}\right|^2d\mu_T(t)<\infty \text{ for almost every } \omega\in\Omega,\label{Evo21}
\end{align}
where $\mathcal{R}_t$ for each $t\in T$ for any $f\in\mathbb{L}_2(\mathbb{R}^d)$ is defined as 
$$\mathcal{R}_t f(x)=\frac{1}{\sqrt{\text{det} (\varphi(t))}}f((\varphi(t)^{-1})x).
$$ 
For an admissible vector $\psi\in H$, the span of $V_{\psi}=\{\mathcal{U}_{g}\psi|g\in G\}$, is dense in $H$. Further the corresponding CS transform $\mathcal{W}_{\psi}:H\rightarrow\mathbb{C}_{K}^{G}$ is unitary. This follows from general results in \cite[Ch.5]{Fuhr2005}.

\begin{theorem}\label{MPsiRecon}
Let $G=\mathbb{R}^d\rtimes T$ and $\psi$ be an admissible  vector. Then $\mathcal{T}_{M_\psi}\Phi\in \mathbb{L}_{2}(G,d\mu_{G}(g))$ for all $\Phi\in \mathbb{C}_{K}^{G}$, where we define
\begin{align*}
[\mathcal{T}_{M_{\psi}}[\Phi]](\mathbf{b},t):=\mathcal{F}^{-1}\bigg{[}\boldsymbol{\omega}\mapsto(2\pi)^{-d/4}M_{\psi}^{-1/2}(\boldsymbol{\omega})\mathcal{F}[\Phi(\cdot,t)](\boldsymbol{\omega}) \bigg{]}(\mathbf{b}).
\end{align*} 
Therefore $(\cdot,\cdot)_{M_{\psi}}:\mathbb{C}_{K}^{G}\times \mathbb{C}_{K}^{G}\rightarrow\mathbb{C}$ defined by
\begin{equation}\label{MPsiInnProd}
(\Phi,\Psi)_{M_{\psi}}=(\mathcal{T}_{M_{\psi}}[\Phi],\mathcal{T}_{M_{\psi}}[\Psi])_{\mathbb{L}_{2}(G)},
\end{equation}
is an explicit characterization of the inner product on $\mathbb{C}_{K}^{G}$, which is the unique functional Hilbert space with reproducing kernel $K:G\times G\rightarrow \mathbb{C}$ given by
\begin{equation}
K(g,h)=(\mathcal{U}_{g}\psi,\mathcal{U}_{h}\psi)_{\mathbb{L}_{2}(\mathbb{R}^d)}=(\mathcal{U}_{h^{-1}g}\psi,\psi)_{\mathbb{L}_{2}(\mathbb{R}^d)}, \ g,h\in G.
\end{equation}
The wavelet transformation $\mathcal{W}_{\psi}:H\rightarrow\mathbb{C}_K^G$ given by 
\begin{equation}
\mathcal{W}_{\psi}[f](\mathbf{b},t)=(\mathcal{U}_{g}\psi,f)_{\mathbb{L}_{2}(\mathbb{R}^d)}, \ f\in\mathbb{L}_{2}(\mathbb{R}^d), \ g=(\mathbf{b},t)\in\mathbb{R}^d\rtimes_\varphi T,
\end{equation}
is a unitary mapping from $H$ to $\mathbb{C}_{K}^{G}$.
The space $\mathbb{C}_{K}^{G}$ is a closed subspace of the Hilbert space $\mathbb{H}_{\psi}\otimes\mathbb{L}_2(T;\frac{d\mu_{T}(t)}{det(\varphi(t))})$, where
\begin{align*}
\mathbb{H}_{\psi}=\{f\in H|\ M_{\psi}^{-\frac{1}{2}}\mathcal{F}[f]\in\mathbb{L}_{2}(\mathbb{R}^d)\}
\end{align*}  is equipped with the inner product 
\begin{equation*}
(f_1,f_2)=(M_{\psi}^{-\frac{1}{2}}\mathcal{F}[f_1],M_{\psi}^{-\frac{1}{2}}\mathcal{F}[f_2])_{\mathbb{L}_{2}(\mathbb{R}^d;(2\pi)^{-d/2}d\mathbf{x})}, \text{ for all } f_1, \ f_2\in H.
\end{equation*}
The orthogonal projection $\mathbb{P}_{\psi}$ of $\mathbb{H}_{\psi}\otimes\mathbb{L}_2(T;\frac{d\mu_{T}(t)}{det(\varphi(t))})$ onto $\mathbb{C}_{K}^{G}$ is given by $(\mathbb{P}_{\psi}[\Phi])(g)=(K(\cdot,g),\Phi))_{M_{\psi}}$.
\end{theorem}
\textbf{Proof.} See \ref{App:MpsiRecon} for proof. $\hfill\Box$
\begin{remark}
Since $\mathcal{W}_{\psi}:H\rightarrow\mathbb{C}_{K}^{G}$ is unitary, the inverse equals the adjoint and thus the image $f$ can be reconstructed from its orientation score $\mathcal{W}_{\psi}[f]$ by
\begin{align}
f=\mathcal{W}_{\psi}^*[\mathcal{W}_{\psi}[f]]\nonumber 
=\mathcal{F}^{-1}\bigg{[}\boldsymbol{\omega}\mapsto\frac{1}{(2\pi)^{d/2}}\int\limits_{T}\mathcal{F}[\mathcal{W}_{\psi}[f](\cdot,t)](\omega)\mathcal{F}[\mathcal{R}_{t}\psi](\boldsymbol{\omega}) \frac{d\mu_{T}}{|det(\tau(t))|}M_{\psi}^{-1}(\boldsymbol{\omega}) \bigg{]}.
\end{align}
\end{remark}
\subsection{Coherent State Transform on the $SE(2)$ group (Orientation Score)}
In this case $G=SE(2)$, $(\mathcal{U}_{g=(\mathbf{b},\theta)}\psi)(\mathbf{x})=\psi\left(\mathbf{R}_\theta^{-1}(\mathbf{x}-\mathbf{b})\right)$, $\Omega=\mathbb{R}^2$ and $H=\mathbb{L}_2(\mathbb{R}^2)\cap\mathbb{L}_1(\mathbb{R}^2)$ leading to \cite[Thm 1]{Duits2007a}. For more details on construction as well operators in the coherent state domain see \cite{Duits2010,Duits2010a}.
For $\psi\in\mathbb{L}_2(\mathbb{R}^2)\cap\mathbb{L}_1(\mathbb{R}^2)$, $M_\psi$ is a continuous function vanishing at infinity. Note that $C_\psi\neq M_\psi$ in this case as $M_\psi=1$ and therefore unitary operator $\mathcal{W}_\psi:\mathbb{L}_2(\mathbb{R}^2)\rightarrow\mathbb{L}_2(SE(2))$ cannot be obtained. In order to deal with $\mathbb{L}_1-$isometries one needs to either rely on distributions $\psi\in\mathbb{H}^{-s}(\mathbb{R}^2)$ with $s>1$ and appropriate distributional wavelet transforms \cite[App.B]{Bekkers2012} or one needs to restrict to disc-limited images, which is appropriate for imaging applications.   

\subsection{Continuous Wavelet Transform on $SIM(2)$ group (Multiple Scale Orientation Scores)}
Consider the case $T:=SO(2)\times\mathbb{R}^{+}$, $G:=SIM(2)=\mathbb{R}^2\rtimes_\varphi (SO(2)\times\mathbb{R}^{+})$
equipped with the group product
\begin{equation}
(\mathbf{x},a,\theta)(\mathbf{x} ',a',\theta ')=(\mathbf{x}+\varphi(a,\theta)\mathbf{x} ',aa',\theta +\theta ' ), \ \ \text{for all } (\mathbf{x},a,\theta),(\mathbf{x} ',a',\theta ')\in SIM(2),
\end{equation}
where $\varphi(a,\theta)=a\mathbf{R}_{\theta}$.
Consider the unitary representation of $SIM(2)$ in $\mathbb{L}_{2}(\mathbb{R}^2)$ given by,
\begin{equation}
\mathcal{U}_{g=(\mathbf{b},a,\theta)} \psi(\mathbf{x})=\frac{1}{a}\psi\bigg{(}\frac{\mathbf{R}
_{\theta}^{-1}(\mathbf{x}-\mathbf{b})}{a} \bigg{)}, \ a>0, \ \theta\in[0,2\pi], \ \mathbf{b}\in\mathbb{R}^2,
\end{equation}
We denote $\mathcal{U}:(\mathbf{x},t)=(\mathbf{x},a,\theta)\mapsto \mathcal{U}_{(\mathbf{x},t)}$ as 
\begin{equation}
\mathcal{U}_{\mathbf{x},t}f=\mathcal{T}
_{\mathbf{x}}\mathcal{R}_{t}f,   \ \  t=(a,\theta)\in\mathbb{R}^+\times SO(2)
\label{34}
\end{equation}
where $(\mathcal{T}_{\mathbf{x}}f)(\mathbf{x}')=f(\mathbf{x}'-\mathbf{x}) \text{ and } (\mathcal{R}_{t}f)(\mathbf{x}')=\frac{1}{a} f(\frac{1}{a}\mathbf{R}_{\theta}\mathbf{x}')$ for all $\mathbf{x},\mathbf{x}'\in\mathbb{R}^2$.
In the continuous setting of $SIM(2)$ Theorem~\ref{MPsiRecon} yields
\begin{align}
M_{\psi}(\boldsymbol{\omega})=2\pi
\int\limits_{0}^{2\pi}\int\limits_
{\mathbb{R}}
|\hat{\psi}(e^{\tau}\mathbf{R}_{\theta}^{-1}
\boldsymbol{\omega})|^2 d\tau d\theta=C_\psi, \label{25}
\end{align} 
and the reconstruction formula coincides with the known results on unitary CW transforms \cite{Grossmann1985,Ali2014}. Here we have used the notation $\hat{\psi}=\mathcal{F}\psi$.
\subsection{Discrete Implementation}
\label{sec:DiscreteAnalogue}
In this subsection we deal with the practical aspects of the implementation of the continuous wavelet transform discussed earlier.
Recall that $SIM(2)=\mathbb{R}^{2}\rtimes(SO(2)\times \mathbb{R}^{+})$ where $T=(SO(2)\times \mathbb{R}^{+})$ is a locally compact group. We can replace $SO(2)$ to be finite rotation group, denoted by $\mathbb{T}_{N}$ (equipped with discrete topology) which is locally compact  i.e.\ 
\begin{equation}
\mathbb{T}_{N}=\{e^{iks_{\varphi}}|
k\in\{0,1,\ldots ,N-1\},s_{\varphi}=\frac{2\pi}{N}\},  \text{ for } N\in \mathbb{N},
\label{29}
\end{equation}
and obtain $SE(2,N)$, see \cite{Boscain2013} for details.
On the other hand the scaling group $\mathbb{R}^{+}$ cannot be written in terms of a finite scaling group as every finite subgroup of a multiplicative group of a field is a cyclic subgroup.

Further from a practical point of view in the discrete case we need to have a lower and an upper bound on the choice of the scales. We assume that $a\in[a^-,a^+]$ where $0<a^- <a^+ $ and consider the following discretization for scales,
\begin{equation}
\mathbb{D}_{M}=\bigg{\{} 
e^{(\tau^- + ls_\rho)}
\bigg{|}
l\in\{0,1,\ldots\,M-1\},
s_\rho=\frac{\tau^+-\tau^-}
{M}
\bigg{\}}, \text{ for } M\in \mathbb{N},
\label{28}
\end{equation}
where $\tau^- = \log (a^-)$ and $\tau^+=\log (a^+)$.
Using the notation $\theta_k=ks_{\varphi}$ and $a_l=a^-e^{ls_\rho}$ we write the discrete version of \eqref{34}
\begin{align}
U_{f}^{N,M}(\mathbf{b},a_l,\theta_{k})=
(\mathcal{T}_{\mathbf{b}}
\mathcal{R}_{(a_l,\theta_k)}
\psi ,f)_{\mathbb{L}_{2}(\mathbb{R}^2)},
\label{35}
\end{align}
which is the discrete version of wavelet transform of an image $f\in\mathbb{L}_2(\mathbb{R}^2)$. 
 The discrete version of
$M_{\psi}$ is,
\begin{align}
M_{\psi}^D(\boldsymbol{\omega})=\frac{1}{N}\frac{1}{M}\sum\limits_{k=0}^{N-1}\sum\limits_{l=0}
^{M-1}\frac{1}{a}\hat{\psi}\left(
(a_lR_{\theta_k})^{-1}\boldsymbol{\omega}\right).
\end{align}
Note that in the discrete setting we are no longer in the unitary irreducible setting of \cite{Grossmann1985} and $M_\psi^D$ is not a constant (so $M_\psi^D\neq C_\psi$). Furthermore the space of wavelet transforms is embedded in $\mathbb{H}_\psi\otimes \mathbb{L}_2(SO(2)\times\mathbb{R}^+)$ instead of $\mathbb{L}_2(\mathbb{R}^2\times SO(2)\times\mathbb{R}^+)$ 
\subsection{Stable reconstruction of an image from OS}
The unitarity result (Theorem~\ref{MPsiRecon}) with $\mathbb{C}_K^{SIM(2)}\subset \mathbb{H}_{\psi}\otimes\mathbb{L}_2
(SO(2)\times\mathbb{R}^{+};\frac{d\mu_{T}(t)}{det(\tau(t))})$ depends on the wavelet $\psi$. For stability  estimates one requires $\mathbb{L}_2$-norms on both the domain and the range. This means we must impose uniform lower and upper bounds in \eqref{Evo21}, which is possible only when we restrict the space of images to functions in $\mathbb{L}_{2}(\mathbb{R}^2)\cap\mathbb{L}_{1}(\mathbb{R}^2)$ whose Fourier transform is contained in an annulus. The space of these images is a Hilbert space given by,
\begin{equation}
\mathbb{L}_{2}^{\varrho^{-},\varrho^{+}}(\mathbb{R}^2)=\{f\in\mathbb{L}_{2}(\mathbb{R}^2)| \ \text{supp}(\mathcal{F}[f])\subset B_{0,\varrho^{+}}\backslash B_{0,\varrho^{-}}\}, \ \  \varrho^{+}>\varrho^{-}>0,
\label{18}
\end{equation}
where $B_{0,\varrho^{\pm}}$ denotes a ball of radius $\varrho^{\pm}$ around the origin in the Fourier domain.
A practical motivation for the assumption of an upper bound ($\varrho^{+}$) on the support of the Fourier transform of the images is the Nyquist theorem, which states that \textit{every band-limited function is determined by its values on a discrete grid}. \\

The value for $\varrho^{-}$  directly relates to the coarsest scale we wish to detect in the spatial domain. Therefore the removal of  extremely low frequencies  from the image essentially corresponds to background removal in the image which is often an essential pre-processing step in medical image processing, see \cite{Bankman2008,Zhang1997}. 

We wish to construct a wavelet transform 
\begin{align}
\mathcal{W}_{\psi}^{\varrho^-,\varrho^+}:\mathbb{L}_{2}^{\varrho^{-},\varrho^{+}}(\mathbb{R}^2)\rightarrow \mathbb{L}_{2}(\mathbb{R}^2\times SO(2)\times\mathbb{R}^+)
\end{align} 
which requires that (recall Eq.\ref{18}), 
\begin{align*}
\mathcal{U}_{\mathbf{x},a,
\theta}\psi\in\mathbb{L}_{2}^{\varrho^{-},
\varrho^{+}}(\mathbb{R}^2), \text{ where }  
a\in[a^-,a^+] \text{ and } 
\theta\in[0,2\pi], \text{ with }a^+>1>a^- \text{ such that } \frac{\rho^-}{a^-}<\frac{\rho^+}{a^+}.
\end{align*}
We make the choice  
\begin{align}
\psi\in\mathbb{L}_{1}(\mathbb{R}^2)\cap\mathbb{L}_{2}(\mathbb{R}^2) \text{ with } \text{supp}(\mathcal{F}[\psi])\subset B_{0,\varrho^+/a^+}\backslash B_{0,\varrho^-/a^-}
\label{eq:Ass-BL}
\end{align}
and thus we have $\text{supp}(\mathcal{F}[\mathcal{U}_{\mathbf{x},a,\theta}\psi])\subset B_{0,\varrho^+}\backslash B_{0,\varrho^-}$, where  $a\in[a^-,a^+] \text{ and } \theta\in[0,2\pi]$.

Note that in our current context, $\psi\in\mathbb{L}_{2}^{\varrho^-,\varrho^+}(\mathbb{R}^2)$ is called an \textit{admissible wavelet} if
\begin{equation}
0<\widetilde{M}_{\psi}=(2\pi)\int\limits_{0}^{2\pi}
\int\limits_{a^-}^{a^+}
\bigg{|}\frac{\mathcal{F}[\mathcal{R}_{a,\theta}\psi]}{\sqrt{det\tau(t)}}\bigg{|}^2 \frac{da}{a}d\theta<\infty   \text{ on } B_{0,\varrho^+}\backslash B_{0,\varrho^-} ,
\label{26}
\end{equation}
where note that $M_\psi=C_\psi\neq \widetilde{M}_\psi$. 
By compactness of the set $[-\pi,\pi]\times [a^-,a^+]$ and assumption in Eq.\eqref{eq:Ass-BL} it follows that $M_\psi$ is a continuous function vanishing at $\infty$. We define, $SIM^-_+:=\mathbb{R}^2\times [-\pi,\pi]\times [a^-,a^+]$.
\begin{definition}
Let $\psi$ be an admissible wavelet  in the sense of \eqref{26}. Then the wavelet transform $\mathcal{W}_{\psi}^{\varrho^-,\varrho^+}:\mathbb{L}_{2}^{\varrho^-,\varrho^+}(\mathbb{R}^2)\rightarrow\mathbb{L}^{2}(SIM^-_+)$ is given by
\begin{equation*}
(\mathcal{W}_{\psi}^{\varrho^-,\varrho^+}[f])(g)=
\frac{1}{a}\int\limits_{\mathbb{R}^2}\overline{\psi \bigg{(}R_{\theta}^{-1}\bigg{(}\frac{\mathbf{y}-\mathbf{x}}{a}\bigg{)}\bigg{)}}f(\mathbf{y})d\mathbf{y}, \ f\in \mathbb{L}_{2}^{\varrho^-,\varrho^+}(\mathbb{R}^2),
\end{equation*}
for almost every $g=(\mathbf{x},a,\theta)\in SIM^-_+$.
\end{definition}

\textbf{Quantification of Stability}
We quantify stability of an invertible linear transformation $A:V\rightarrow W$ from a Banach space $(V,\|\cdot\| _{V})$ to a Banach space $(W,\|\cdot\| _{W})$ via the condition number 
\begin{equation}
\text{cond}(A)=\| A^{-1}\|\| A\|= \left( \sup\limits_{x\in V}\frac{\|\mathbf{x}\| _{V}}{\| A\mathbf{x} \| _{W}} \right)
\left( \sup\limits_{x\in V}\frac{\| A\mathbf{x} \| _{W}}{\|\mathbf{x}\| _{V}} \right)\geq 1.
\end{equation}
The closer it approximates 1, the more stable the operator and its inverse is. The condition number depends on the norms imposed on $V$ and $W$. 
From a practical point of view, it is appropriate to impose the $\mathbb{L}_{2}(SIM^-_+)$-norm instead of the reproducing kernel norm on the score since it does not depend on the choice of the wavelet $\psi$ 
and we also use a $\mathbb{L}_{2}$-norm on the space of images.  
 
\begin{theorem}\label{Thm:StabCond}
Let $\psi$ be an admissible wavelet, with $\widetilde{M}_{\psi}(\boldsymbol{\omega})>0$ for all $\boldsymbol{\omega}\in B_{0,\varrho^+}\backslash B_{0,\varrho^-}$. Then the condition number $\text{cond}(\mathcal{W}_{\psi}^{\varrho^-,\varrho^+})$ of $\mathcal{W}_{\psi}^{\varrho^-,\varrho^+}:\mathbb{L}_{2}^{\varrho^-,\varrho^+}(\mathbb{R}^2)\rightarrow\mathbb{L}_{2}(G), \ (G=SIM^-_+)$ defined by
\begin{equation*}
\text{cond}(\mathcal{W}_{\psi}^{\varrho^-,\varrho^+})=\|(\mathcal{W}_{\psi}^{\varrho^-,\varrho^+})^{-1}\|\|(\mathcal{W}_{\psi}^{\varrho^-,\varrho^+})\|
=\bigg{(} \sup\limits_{f\in\mathbb{L}_{2}^{\varrho^-,\varrho^+}(\mathbb{R}^2)}\frac{\|f\|_{\mathbb{L}_{2}(\mathbb{R}^2)}}{\|
U_{f}\|_{\mathbb{L}_{2}(G)}}\bigg{)}\bigg{(} \sup\limits_{f\in\mathbb{L}_{2}^{\varrho^-,\varrho^+}(\mathbb{R}^2)}\frac{\|
U_{f}\|_{\mathbb{L}_{2}(G)}}{\|f\|_{\mathbb{L}_{2}(\mathbb{R}^2)}}\bigg{)}
\end{equation*}
satisfies
\begin{equation*}
1\leq (\text{cond}(\mathcal{W}_{\psi}^{\varrho^-,\varrho^+}))^2\leq\bigg{(}\sup\limits_{
\varrho^- \leq\|\boldsymbol{\omega}\|\leq\varrho^+}\widetilde{M}_{\psi}^{-1}(\boldsymbol{\omega})\bigg{)}\bigg{(}\sup
\limits_{
\varrho^- \leq\|\boldsymbol{\omega}\|\leq\varrho^+}\widetilde{M}_{\psi}(\boldsymbol{\omega})\bigg{)}.
\end{equation*}
\end{theorem}
\textbf{Proof.}
Since $\widetilde{M}_{\psi}>0$ and is continuous on the compact set $ B_{0,\varrho^+}\backslash B_{0,\varrho^-}=\{\boldsymbol{\omega}
\in\mathbb{R}^2 \lvert \varrho^- \leq\|
\boldsymbol{\omega}\|\leq\varrho^+\}$, 
$\sup\limits_{\varrho^- <\|
\boldsymbol{\omega}\|<\varrho^+}\widetilde{M}_{\psi}
(\boldsymbol{\omega})=\max\limits_{\varrho^- 
<\|\boldsymbol{\omega}\|<\varrho^+}\widetilde{M}_{\psi}
(\boldsymbol{\omega})$ do exist. The same 
holds for $\widetilde{M}_{\psi}^{-1}$. Furthermore 
for all $f\in\mathbb{L}_{2}^{\varrho^-,
\varrho^+}(\mathbb{R}^2)$, the restriction 
of the corresponding wavelet transform to 
fixed orientations and scales also belong to the same space, i.e. $U_{f}
(\cdot,a,e^{i\theta})\in
\mathbb{L}_{2}^{\varrho^-,
\varrho^+}(\mathbb{R}^2)$, where $\theta\in 
[0,2\pi]$, $a\in[a^-,a^+]$. Using $\check{\psi}(x)=\psi(-x)$,
\begin{align}
(\mathcal{W}_{\psi}f)(\mathbf{x},t)= \int\limits_{\mathbb{R}^2}\overline{\mathcal{R}_{t}\check{\psi}(\mathbf{x}-\mathbf{x}')}f(\mathbf{x}')d\mathbf{x}'
= (\overline{\mathcal{R}_{t}\check{\psi}}*_{\mathbb{R}^{2}} f)(\mathbf{x})=\left(\mathcal{F}^{-1}(\overline{\mathcal{F}\mathcal{R}_{t}\psi} \mathcal{F}f)\right)(\mathbf{x})\label{15}, 
\end{align}
which gives, $\mathcal{F}[U_{f}](\boldsymbol{\omega})=\mathcal{F}[\overline{\mathcal{R}_{t}\check{\psi}}](\boldsymbol{\omega})\mathcal{F}[f](\boldsymbol{\omega})$. By Theorem~\ref{MPsiRecon} we have $\|f\|^{2}_{\mathbb{L}_{2}(\mathbb{R}^2)}=\|U_{f}\|^{2}_{M_{\psi}}:=(U_f,U_f)_{M_\psi}$, 
\begin{align*}
(\text{cond}(\mathcal{W}_{\psi}^{\varrho^-,\varrho^+}))^2&=
\bigg{(} \sup\limits_{f\in\mathbb{L}_{2}^{\varrho^-,\varrho^+}(\mathbb{R}^2)}\frac{\| U_f \|_{\widetilde{M}_{\psi}}}{\|
U_{f}\|_{\mathbb{L}_{2}(G)}}\bigg{)}\bigg{(} \sup\limits_{f\in\mathbb{L}_{2}^{\varrho^-,\varrho^+}(\mathbb{R}^2)}\frac{\|
U_{f}\|_{\mathbb{L}_{2}(G)}}{\| U_f \|_{\widetilde{M}_{\psi}}}\bigg{)}\\
&\leq \bigg{(}\sup\limits_{
\varrho^- \leq\|\boldsymbol{\omega}\|\leq\varrho^+}\widetilde{M}_{\psi}^{-1}(\boldsymbol{\omega})\bigg{)}\bigg{(}\sup
\limits_{
\varrho^- \leq\|\boldsymbol{\omega}\|\leq\varrho^+}\widetilde{M}_{\psi}(\boldsymbol{\omega})\bigg{)}.
\end{align*}
Finally we note that  $1=\|(\mathcal{W}_{\psi}^{\varrho^-,\varrho^+})^{-1}(\mathcal{W}_{\psi}^{\varrho^-,\varrho^+})\|\leq\|(\mathcal{W}_{\psi}^{\varrho^-,\varrho^+})^{-1}\|\|(\mathcal{W}_{\psi}^{\varrho^-,\varrho^+})\|$.
$\hfill\Box$

\begin{corollary}
The stability of the (inverse) wavelet transformation 
$\mathcal{W}_{\psi}^{\varrho^-,\varrho^+}:\mathbb{L}_{2}^{\varrho^-,\varrho^+}(\mathbb{R}^2)\rightarrow\mathbb{L}_{2}(SIM^-_+)$ is optimal if $\widetilde{M}_{\psi}(\boldsymbol{\omega})=constant$ for all $\boldsymbol{\omega}\in\mathbb{R}$, with $\varrho^-\leq\|\boldsymbol{\omega}\|\leq\varrho^+$. 
\end{corollary}
%
In case 
$\widetilde{M}_{\psi}\approx 1_{B_{0,\varrho^+}\backslash  
B_{0,\varrho^- }}$, the reconstruction formula can be simplified to,
\begin{equation}
f\approx\mathcal{F}^{-1}\bigg{[}
\boldsymbol{\omega}\mapsto\frac{1}
{(2\pi)}\int\limits_{a^{-}}^{a^+}
\int\limits_{0}^{2\pi}\mathcal{F}[U_f(\cdot,a,e^{i\theta})]
(\boldsymbol{\omega})\mathcal{F}[\mathcal{R}_{a,e^{i\theta}}\psi](\boldsymbol{\omega})d\theta \frac{da}{a}  \bigg{]}.
\end{equation}

\subsection{Design of Proper Wavelets}
In this sequel a wavelet $\psi\in\mathbb{L}_{2}(\mathbb{R}^2)\cap\mathbb{L}_{1}(\mathbb{R}^2)$ with $\widetilde{M}_{\psi}$ smoothly approximating $1_{B_{0,\varrho^+}\backslash  
B_{0,\varrho^- }}$, is called a \textbf{proper wavelet}. The entire class of proper wavelets allows for a lot of freedom in the choice of $\psi$. In practice it is mostly sufficient to consider wavelets that are similar to the long elongated patch one would like to detect and orthogonal to structures of local patches which should not be detected, in other words employing the basic principle of template matching. We restrict the possible choices by listing below certain practical requirements to be fulfilled by our transform.
\begin{enumerate}
\item{The wavelet transform should yield a finite number of orientations ($N$) and scales ($M$).}
\item{The wavelet should be strongly directional, in order to obtain sharp responses on oriented structures.}
\item{The transformation should handle lines, contours and oriented patters. Thus the wavelet should pick up edge, ridge and periodic profiles.}
\item{In order to pick up local structures, the wavelet should be localized in spatial domain.}
\end{enumerate}
To ensure that the wavelet is strongly directional and minimizes uncertainty in $SIM(2)$, we require that the support of the wavelet be contained in a convex cone in the Fourier domain, \cite{Antoine1996}. 
The following lemma gives a simple but practical approach to obtain proper wavelets $\psi$, with $M_{\psi}
=1_{B_{0,\varrho^+}\backslash B_{0,\varrho^- }}$. 
\begin{lemma}\label{20}
Let $\tau^-,\tau^+$ be chosen such that $\tau^-=\log(a^-)$ and $\tau^+=\log(a^+)$, where $0<a\in[a^-,a^+]$ is the finite interval of scaling. 
Let $A:SO(2)\rightarrow\mathbb{R}^+$ and $ B:[\tau^{-},\tau^+]\rightarrow\mathbb{R}^+$  such that
\begin{align}
2\pi\int\limits_{0}^{2\pi}A(\varphi)d\varphi=1, \ 
\int\limits_{\tau^-}^{\tau^+}B(\rho)d\rho=1, \label{22}
\end{align}
then the wavelet $\psi=\mathcal{F}^{-1}[\boldsymbol{\omega}\rightarrow\sqrt{A
(\varphi)B(\rho)}]$ has $\widetilde{M}_{\psi}(\boldsymbol{\omega})= 1$ for all $\boldsymbol{\omega}=(\rho\cos\varphi,\rho\sin\varphi)\in B_{0,\varrho^+}\backslash  
B_{0,\varrho^- }$.
\end{lemma}
\textbf{Proof.}
From \eqref{25} and \eqref{26}, for all $\boldsymbol{\omega}\in B_{0,\varrho^+}\backslash  
B_{0,\varrho^- }$ we have,
\begin{align*}
M_{\psi}(\boldsymbol{\omega})=2\pi
\int\limits_{0}^{2\pi}\int\limits_
{\varrho^-}^{\varrho^+}
|\hat{\psi}(e^{\tau}\mathbf{R}_{\theta}^{-1}
\boldsymbol{\omega})|^2 d\tau d\theta =
2\pi
\int\limits_{0}^{2\pi}\int\limits_
{\varrho^-}^{\varrho^+}
|\sqrt{A(\varphi-\theta)B(e^{\tau} \rho)}|^2 d\tau d\theta =1. \quad \qquad \qquad \qquad \qquad \qquad \Box
\end{align*}

Lemma \ref{20} can be translated into the discrete framework, recall \eqref{29} and \eqref{28}, making condition \eqref{22},
\begin{equation}
\frac{1}{N}\sum\limits_{k=0}^{N-1}A(\varphi-\theta_k)=1 \ \text{and} \  \frac{1}{M}\sum\limits_{l=0}^{M-1}B(e^{\tau_l}\rho)=1.     
\end{equation}
where we have made use of discrete 
notations introduced in \eqref{35}. 
\begin{figure}[t]
\centering
\includegraphics[scale=1]{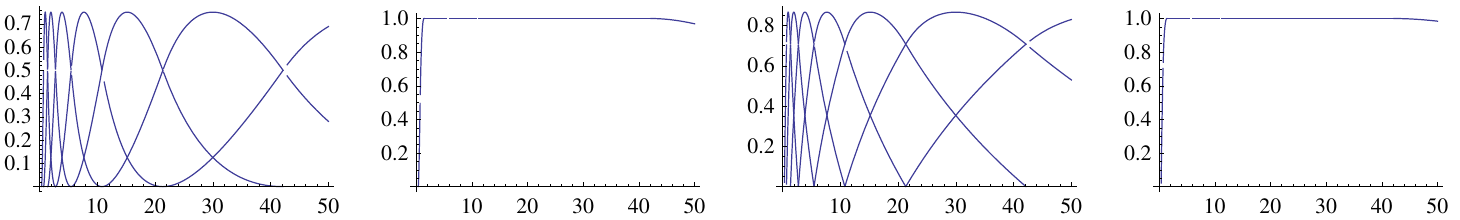}
\caption[B-splines in the log-scale]{
Plot for the B-splines described in \eqref{30}. Values chosen: $a^- =\varrho^- =10^{-8}, \ a^+ =\varrho^+ =50, \ N_s=7$. See Eq.~\eqref{30},~\eqref{eq:BSpline-Properties}  for details.
From Left to Right: (4.1) Plot of $3$rd order B-splines i.e. $k=3$. The B-splines are skewed because of non linear sampling in the scale dimension; (4.2) Summing over the B-splines; (4.3) Plot of square root of the B-splines; (4.4) Summing over the square root of B-Spines.
}
\label{fig:LogSqrtSpline}
\end{figure}
If $2\pi\int\limits_{0}^{2\pi}\sqrt{A(\varphi)}d\varphi\approx 1$ and $\int\limits_{\tau^-}^{\tau^+}\sqrt{B(\rho)}d\rho\approx 1$, we have a fast and simple approximative reconstruction:
\begin{align}
\tilde{f}(\mathbf{x})&=2\pi\int\limits_{0}^{2\pi}
\int\limits_{\tau^-}^{\tau^+}(\mathcal{W}
_{\psi}f)(\mathbf{x},\tau,\theta)d\tau d\theta 
 \approx \mathcal{F}^{-1}[\boldsymbol{\omega}\mapsto (\sqrt{\widetilde{M}_{\psi}}*\mathcal{F}[f](\boldsymbol{\omega}))](\mathbf{x}),\text{ for a.e. } \boldsymbol{\omega}\in B_{0,\varrho^+}\backslash  
B_{0,\varrho^- }. \nonumber\\
&\sqrt{\widetilde{M}_{\psi}}\approx 1_{B_{0,\varrho^+}\backslash  
B_{0,\varrho^- }}\Rightarrow\tilde{f}\approx
f\in \mathbb{L}_{2}^{\varrho^-,\varrho^+}(\mathbb{R}^2).
\label{31}
\end{align} 
We need to fulfil the requirement $\widetilde{M}_{\psi}(\boldsymbol{\omega})\approx 1$ with an appropriate choice of kernel satisfying the additional condition $\sqrt{\widetilde{M}_{\psi}}\approx 1$.
The idea is to ``\textit{fill the cake by pieces of cake}" in the Fourier domain. 
In order to avoid high frequencies in the spatial domain, these pieces must be smooth and they must overlap. A choice of B-spline based functions in the angular and the log-radial direction is an appropriate choice for such a wavelet kernel. This design of wavelets in the Fourier domain is similar to the framework of curvelets \cite{Donoho2005,Donoho2005a,Candes2006}.
However, our decomposition of unity in the Fourier domain is more suited for the subsequent design of left-invariant diffusions in the wavelet domain. 
The reason for us to choose log-polar  B-spline decomposition is that the canonical coordinates of the second kind in $SIM(2)$ are given by
\begin{align*}
\xi=a(x\cos\theta+y\sin\theta), \  \eta=a(-x\sin\theta+y\cos\theta), \  \theta,  \ \tau=\log_e a,
\end{align*}  
and it will turn out that our subsequent evolutions are best expressed in these coordinates.  
We would also like to point out that a similar construction which involves mixing different levels of scale selectivity was presented in \cite{Jacques2007}.

The $k^{\text{th}}$ order B-spline denoted by $B^{k}$ is defined as
\begin{equation}
B^{k}(x)=(B^{k-1}*B^{0})(x),  \ B^{0}(x)=\left\{\begin{array}{c l}
  1  &  if -1/2<x<+1/2 \\
  0  & otherwise
\end{array}
\right.
\end{equation} 
with the property that B-splines add up to 1. For more details see \cite{Unser1999}.
Based on the requirements and considerations above we propose the following kernel
\begin{equation}\label{23}
\psi(\mathbf{x})=
\mathcal{F}^{-1}_{\mathbb{R}^2}[\boldsymbol{\omega}\rightarrow
\sqrt{A(\varphi)B(\rho)}](\mathbf{x})G_{{\sigma}_{s}}(\mathbf{x}),
\end{equation}
where $G_{{\sigma}_{s}}$ is a Gaussian window that enforces spatial locality cf. requirement 5.
$A:\mathbb{T}\rightarrow\mathbb{R}^{+}$ and $B:[\varrho^-,\varrho^+]\rightarrow\mathbb{R}^{+}$ are defined as,
\begin{equation}
A(\varphi)=B^{k}\bigg{(}\frac{(\varphi \ \text{mod} \  2\pi)-\pi/2}{s_{\varphi}}\bigg{)}, \ \ 
B(\varrho)=
B^{k}\left(\frac{\log[\rho]}{s_{\rho}}
\right), 
\label{30}
\end{equation}
where $s_{\varphi}=\frac{2\pi}{N}$ and  $s_{\rho}=(\log[a^+]-\log[a^-])/M$
 where $M$ is the number of chosen scales and $a^-, \ a^+$ are predefined scales, based on   
$\varrho^-, \ \varrho^+$ respectively. See 
Figure \ref{fig:LogSqrtSpline} for a plot of these log B-splines and note that
\begin{align}\label{eq:BSpline-Properties}
\sum\limits_{l=0}^{M-1}B^{k}\left(\frac{\log[\rho]}{s_{\rho}}
+l+\frac{\tau^-}{s_\rho}\right)=1 \text{ and }
\sum\limits_{l=0}^{M-1}\sqrt{B^{k}\left
(\frac{\log[\rho]}{s_{\rho}}
+l+\frac{\tau^-}{s_\rho}\right)}\approx 1.
\end{align}
\begin{figure}[tbp]
\centering
\includegraphics[scale=0.65]{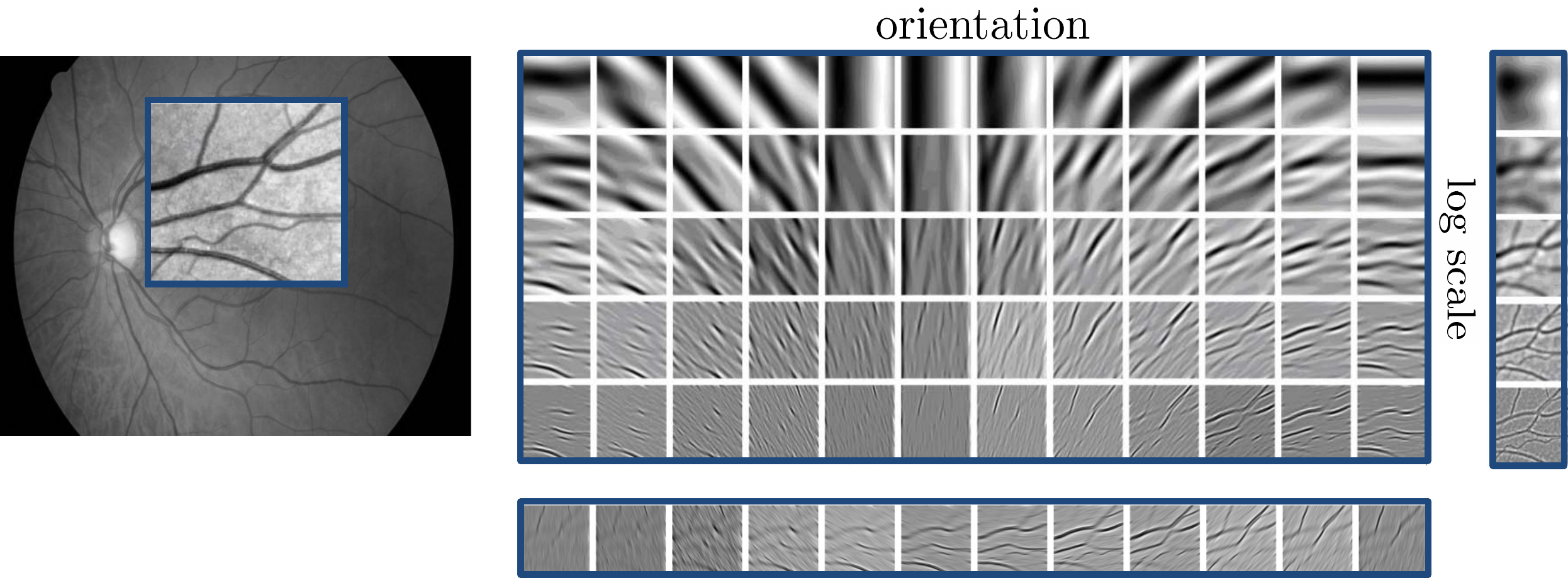}
\caption[ScaleOS for Retinal Image]{Scale-OS of a retinal image. As shown, Scale-OS can be used to create the Orientation Score and Gaussian-Scale Space of the image.
}
\label{fig:ScaleOS_RetinalIm}
\end{figure}
\section{Operators on Scores}
There exists exists a 1-to-1 correspondence between bounded operators $\Phi\in \mathcal{B}(\mathbb{C}_{K}^{G})$ on the range of CW transform and bounded operators $\Upsilon\in \mathcal{B}(\mathbb{L}_{2}(\mathbb{R}^{d}))$ given by
\begin{align}
\Upsilon[f]=(\mathcal{W}
_{\psi}
^{*}\circ\Phi\circ\mathcal{W}_{\psi})[f], \ f\in \mathbb{L}_{2}(\mathbb{R}^{d}),
\end{align}
which allows us to relate operations on transformed images to operations on images in a robust manner. To get a schematic view of the operations see Figure \ref{fig1}. 
By Theorem~\ref{MPsiRecon} that the range of the unitary wavelet transform $\mathcal{W}_\psi:\mathbb{L}_2(\mathbb{R}^2)\rightarrow\mathbb{C}_{K}^{G}$ is a subspace of $\mathbb{L}_{2}(G)$. For proper wavelets we have (approximative) $\mathbb{L}_{2}$-norm preservation and therefore $\mathbb{L}_{2}(G)\cong \mathbb{H}_{\psi}$ (with $\mathbb{L}_2(G)=\mathbb{H}_\psi$ if $M_\psi=1$).
In general if, $\Phi:\mathbb{L}_{2}(G)\rightarrow\mathbb{C}_{K}^{G}$ is a bounded operator  then the range need not be contained in $\mathbb{C}_{K}^{G}$.
Therefore we also consider $\widetilde{\mathcal{W}}_{\psi}:\mathbb{L}_{2}(\mathbb{R}^d) \rightarrow\mathbb{L}_{2}(G)$ given by $\widetilde{\mathcal{W}}_\psi f=\mathcal{W}_\psi f$. Its adjoint is given by, 
\begin{align*}
(\widetilde{\mathcal{W}}_{\psi})^{*}(V)=\int\limits_{G}\mathcal{U}_{g}\psi \  V(g)d\mu_{G}(g), \ V\in \mathbb{L}_{2}(G).
\end{align*}
The operator 
$\mathbb{P}_{\psi}=\widetilde{\mathcal{W}}_{\psi}(\widetilde{\mathcal{W}}_{\psi})^*$ is the orthogonal projection on the space $\mathbb{C}_{K}^{G}$ whereas $\mathcal{W}_\psi\mathcal{W}^*_\psi=I$. This projection can be used to decompose operator $\Phi$:
\begin{align*}
\Phi(U_{f})=\mathbb{P}_{\psi}(\Phi(U_f))+(I-\mathbb{P}_{\psi})(\Phi(U_f)).
\end{align*}
Notice that the orthogonal complement $(\mathbb{C}_{K}^{G})^{\perp}$, which equals $\mathcal{R}(I-\mathbb{P}_{\psi})$, is exactly the null-space of $(\widetilde{\mathcal{W}}_{\psi})^*$ as $\mathcal{N}((\widetilde{\mathcal{W}}_{\psi})^*)=\mathcal{N}((\mathcal{W}_{\psi})^*)=(\mathcal{R}(\mathcal{W}_{\psi}))^{\perp}=(\mathbb{C}^
{K}_{G})^{\perp}$ and so 
\begin{align}
[(\widetilde{\mathcal{W}}
_{\psi})^{*}\circ\Phi\circ\widetilde{\mathcal{W}}_{\psi} ][f]=[(\widetilde{\mathcal{W}}
_{\psi})^{*}\circ\mathbb{P}_{\psi}\circ
\Phi\circ
\widetilde{\mathcal{W}}_{\psi} ][f],
\label{ConvEq3}
\end{align}
for all $f\in \mathbb{L}_{2}(\mathbb{R}^2)$ and all $\Phi\in\mathcal{B}(\mathbb{L}_{2}(G))$, so we see that the net operator associated to $\Phi:\mathbb{L}_2(G)\rightarrow\mathbb{L}_2(G)$ is given by $\mathbb{P}_{\psi}\circ\Phi:\mathbb{L}_2(G)\rightarrow\mathbb{C}_{K}^G$.
In the reminder of this section we present design principles for $\Upsilon$.
\subsection{Design Principles}
We now formulate a few desirable properties of $\Upsilon$, and sufficient conditions for $\Phi$ that guarantee that $\Upsilon$ meets these requirements.
\begin{enumerate}
\item{\textbf{Covariance with respect to rotation, translation and scaling}:
\begin{align}
\Upsilon\circ\mathcal{U}^{SIM(2)}_g=\mathcal{U}^{SIM(2)}_g\circ\Upsilon, \ \ \forall g=(x,y,\tau,\theta)\in SIM(2).\label{ScalingCoV}
\end{align}
This is an important requirement because the net operations on images should not be affected by rotation and translation of the original image. Typically, this is achieved by restricting one self to left-invariant operators $\Phi$. Often we will omit scaling covariance in Eq.\eqref{ScalingCoV} as in many imaging applications this is not natural and therefore we require \eqref{ScalingCoV} only to hold for the $SE(2)$ subgroup.}
\item{\textbf{Left invariant vector fields}: In order to achieve the Euclidean invariance mentioned above, we need to employ left invariant vector fields on $SIM(2)$ as a moving frame of reference.}
\item{\textbf{Nonlinearity}: The requirement that $\Upsilon$ commute with $\mathcal{U}$ immediately rules out linear operators $\Phi$. Recall that $\mathcal{U}$ is irreducible, and by Schur's lemma \cite{Dieudonne1977}, any linear intertwining operator is a scalar multiple of the identity operator.}
\item{\textbf{Left-invariant parabolic evolutions on the Similitude group}:
We consider the following two types of evolutions which include the wavelet transform as a initial condition. \begin{itemize}{
\item{Combine linear diffusions with monotone operations on the co-domain }
\item{Non linear adaptive diffusion}
}
\end{itemize}}
\item{\textbf{Probabilistic models for contextual multi-scale feature propagation in the wavelet domain}: 
Instead of uncorrelated soft-thresholding of wavelet coefficients we aim for PDE flows that amplify the wavelet coefficients which are probabilistically coherent w.r.t. neighbouring coefficients. This coherence w.r.t. neighbouring coefficients is based on underlying stochastic processes (random walks) for multiple-scale contour enhancement.}
\end{enumerate}
In Subsections~\ref{Left-Invariance}-\ref{Stoc_Conn} we will elaborate on these design principles.  
\subsection{Covariance with respect to Rotations and Translations}\label{Left-Invariance}
Let $G=\mathbb{R}^d\rtimes T$ denote an arbitrary affine Lie-group. 
\begin{definition}
\label{ConvDef1}
An operator $\Phi:\mathbb{L}_{2}(G)\rightarrow\mathbb{L}_{2}(G)$ is \textit{left invariant} iff
\begin{equation}
\Phi[\mathcal{L}_{h}f]=\mathcal{L}_{h}[\Phi f], \ \text{for all } h\in G, \  f\in \mathbb{L}_{2}(\mathbb{R}^2),
\end{equation} 
where the left regular action $\mathcal{L}_g$ of $g\in G$ onto $\mathbb{L}_{2}(G)$ is given by
\begin{equation}
\mathcal{L}_{g}\psi(h)=\psi( g^{-1}h).
\end{equation} 
\end{definition} 
\begin{theorem}
\label{thm:EuclideanInv}
Let $\Phi:\mathbb{C}_{K}^{G}\rightarrow\mathbb{L}_{2}(G)$ be a bounded operator. Then the unique corresponding operator $\Upsilon$ on $H$ given by 
$\Upsilon[f]=(\widetilde{\mathcal{W}}_{\psi})^{*}\circ
\Phi\circ\widetilde{\mathcal{W}}_{\psi} [f]$ is Euclidean (and scaling) invariant, i.e.\ $\mathcal{U}_{g}\Upsilon=\Upsilon\mathcal{U}_{g}$ for all $g\in G$ if and only if $\mathbb{P}_{\psi}\circ\Phi$ is left invariant, i.e.\ $\mathcal{L}_{g}(\mathbb{P}_{\psi}\circ\Phi)
=(\mathbb{P}_{\psi}\circ\Phi)
\mathcal{L}_{g}$, for all $g\in G$.
\end{theorem}
\textbf{Proof.}
See Appendix~\ref{app:Covariance-Proof} for proof. $\hfill\Box$\\

\textit{Practical Consequence}: Now let us return to our case of interest $G=SIM(2)$, where $H=\mathbb{L}_2(\mathbb{R}^2)$. Euclidean invariance of $\Upsilon$  is of great practical importance, since the result of operators on scores should not be essentially different if the original image is rotated or translated. In addition in our construction scaling the image also does not affect the outcome of the operation. The latter constraint may not always be desirable in applications. 
\begin{remark}
Other CW transform based approached such as shearlets also carry a group structure. However this is a different group structure corresponding for e.g. with the shearing, translation and scaling group. Although such a group structure follows by a contraction (see \ref{App:MethContraction} \eqref{eq:App-Contarction-Shear}) it is not suited for exact rotation and translation covariant processing.  
\end{remark}
\subsection{Left Invariant Vector fields (differential operators) 
on $SIM(2)$}\label{Sec:LeftInvVecFiel}
Left invariant differential operators are crucial in the construction of appropriate left-invariant evolutions on $G=SIM(2)$. Similar to Definition \ref{ConvDef1}, the right regular action $\mathcal{R}_g$ of $g\in G$  onto $\mathbb{L}_{2}(G)$ is defined by
\begin{align}
\mathcal{R}_g\psi(h)=\psi(hg), \ \forall g,h\in G, \ \psi\in\mathbb{L}_2(G).
\end{align}
A vector field (now considered as a differential 
operator\footnote{Any tangent vector $X\in T(G)$ can be 
considered as a differential operator acting on a function 
$U:G\rightarrow \mathbb{R}$. So, for instance, if we are using  
$X_{e}\in T_{e}(G)$ in the context of differential operators, 
all occurrences of $\mathbf{e}_{i}$ will be replaced by 
$\partial_{i}$, which is the short-hand notation for 
$\frac{\partial}{\partial x_{i}}$. See \cite{Aubin2001} for the equivalence between these two viewpoints.}) $\mathcal{A}$ on a group $G$ 
is called left-invariant if it satisfies
\begin{equation*}
\mathcal{A}_{g}\phi=\mathcal{A}_{e}(\phi\circ L_{g})=\mathcal{A}_{e}(h\mapsto\phi(gh)),
\end{equation*} 
for all smooth functions $\phi\in C_{c}^{\infty}(\Omega_{g})$ where $\Omega_{g}$ is an open set around $g\in G$ and with the left multiplication $L_{g}:G\rightarrow G$ given by $L_g h=gh$. The linear space of left-invariant vector fields $\mathcal{L}(G)$ equipped with the Lie product $[\mathcal{A},\mathcal{B}]=\mathcal{A}\mathcal{B}-\mathcal{B}\mathcal{A}$ is isomorphic to $T_{e}(G)$ by means of the isomorphism,
\begin{equation*}
T_{e}(G)\ni A\leftrightarrow\mathcal{A}\in\mathcal{L}(G)\Leftrightarrow\mathcal{A}_{g}(\phi)=A(\phi\circ L_{g})=A(h\mapsto\phi(gh))=(L_g)_*A(\phi)
\end{equation*} 
for all smooth $\phi:G\supset\Omega_{g}\rightarrow\mathbb{R}$. 

We define an operator $d\mathcal{R}:T_{e}(G)\rightarrow\mathcal{L}(G)$,
\begin{equation}
(d\mathcal{R}(A)\phi)(g):=\lim\limits_{t\downarrow 0}\frac{(\mathcal{R}_{\exp(tA)}\phi)(g)-\phi(g)}{t}, \ A\in T_{e}(G), \ \phi\in\mathbb{L}_{2}(G), \ g\in G, 
\label{SIM2_12}
\end{equation}
and where $\mathcal{R}$ and $\exp$ are the right regular representation and the exponential map respectively. 
Using $d\mathcal{R}$ we obtain the corresponding basis for left-invariant vector fields on $G$:
\begin{equation}
\{\mathcal{A}_{1},\mathcal{A}_{2},
\mathcal{A}_{3},\mathcal{A}_{4}\}:=\{d\mathcal{R}(A_{1}),d\mathcal{R}(A_{2}),
d\mathcal{R}(A_{3}),,d\mathcal{R}(A_{4})\},
\end{equation}
or explicitly in coordinates
\begin{equation}
\{\mathcal{A}_{1},\mathcal{A}_{2},
\mathcal{A}_{3},\mathcal{A}_{4}\}=\{\partial_{\theta},\partial_{\xi},\partial_{\eta},\partial_{\beta}\}=\{\partial_{\theta},
a(\cos\theta\partial_{x}+
\sin\theta\partial_{y}),a(-\sin\theta
\partial_{x}+
\cos\theta\partial_{y}),a\partial_{a}\},
\end{equation}
where we use the short notation $\partial_a:=\frac{\partial}{\partial a}$ for the partial derivatives and where,
\begin{equation*}
\{\mathcal{A}_{1}\big{|}_e,\mathcal{A}_{2}\big{|}_e,\mathcal{A}_{3}\big{|}_e,\mathcal{A}_{4}\big{|}_e\}=\{A_{1},A_{2},A_{3},A_{4}\}=\{\partial_{\theta},\partial_{x},\partial_{y},\partial_{a}\}.
\end{equation*}

\begin{figure}[t]
\centering
\includegraphics[scale=0.6]{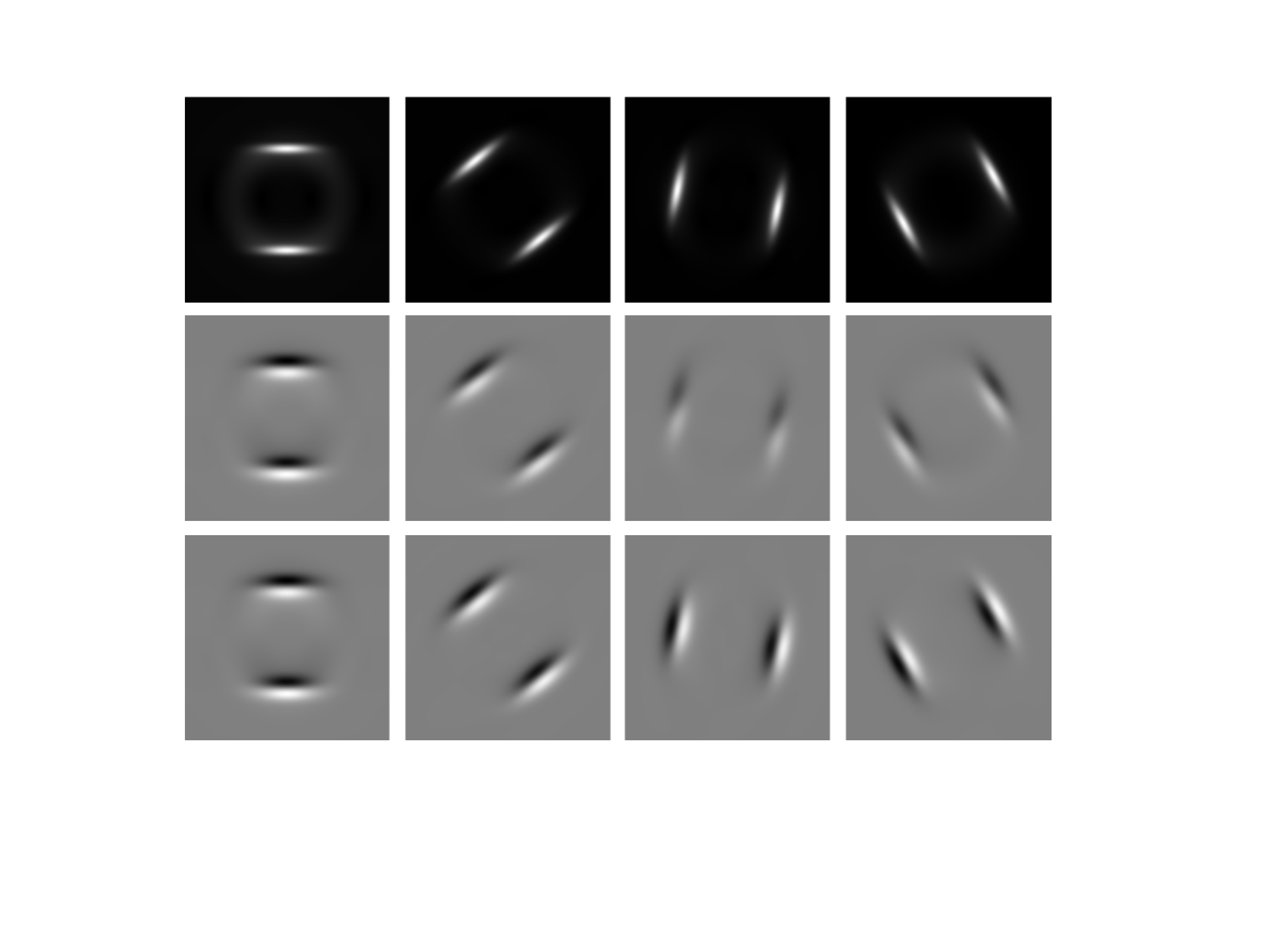}
\caption[Difference between Cartesian and left-invariant derivatives]{The difference between cartesian derivatives and left-invariant derivatives, shown on a scale-OS (at a fixed scale) of an image with a single circle. From left to right, several orientations are shown. 
Row 1: Scale-OS of a circle image at a fixed scale; Row 2: Cartesian derivative $\partial _{y}$; Row 3 : Left-invariant derivative $\partial _{\eta}$. 
Comparing the derivatives $\partial_{y}$ and $\partial_{\eta}$ (Column 3) we observe that $\partial_{\eta}$ is invariant under rotation, i.e.\ the interpretation of $\partial_{\eta}$ stays the same.}
\label{fig:ScematicDerivative}
\end{figure}

To simplify the scale related left invariant differential operator, we introduce a new variable, $\tau=\log a$, which leads to the following change in left invariant derivatives
\begin{equation}
\{\mathcal{A}_{1},\mathcal{A}_{2},
\mathcal{A}_{3},\mathcal{A}_{4}\}=\{\partial_{\theta}
,e^{\tau}(\cos\theta\partial_{x}+
\sin\theta\partial_{y}),e^{\tau}(-\sin\theta\partial_{x}+
\cos\theta\partial_{y}),\partial_{\tau}\}.
\label{2.10}
\end{equation}
The set of differential operators $\{\mathcal{A}_{1},\mathcal{A}_{2},
\mathcal{A}_{3},\mathcal{A}_{4}\}=\{\partial_{\theta},\partial_{\xi},\partial_{\eta},\partial_{\tau}\}$ is the appropriate set of differential operators to be used in orientation scores because all $SIM(2)$-coordinate independent linear and nonlinear combinations of these operators are left invariant. Furthermore at each scale $\partial_{\xi}$ is always the spatial derivative tangent to the orientation $\theta$ and $
\partial_{\eta}$ is always orthogonal to this orientation. Figure \ref{fig:ScematicDerivative} illustrates this for $\partial_{\eta}$ versus $\partial_{y}$.
It is important to note that unlike derivatives $\{\partial_{x},\partial_{y},\partial_{a},\partial_{\theta}\}$ which commute, the left-invariant derivatives $\{\partial_{\xi},\partial_{\eta},\partial_{\tau},\partial_{\theta}\}$ do not commute. However these operators satisfy the same commutator relations as their Lie algebra counterparts as $d\mathcal{R}$ generates a Lie-algebra isomorphism
\begin{align}
\left[A_i,A_j\right]=\sum\limits_{k=1}^4c_{ij}^kA_k\leftrightarrow
[\mathcal{A}_{i},\mathcal{A}_{j}]=\mathcal{A}_{i}\mathcal{A}_{j}-
\mathcal{A}_{j}\mathcal{A}_{i}=\sum\limits_{k=1}^4c^{k}_{ij}\mathcal{A}_k,
\end{align}
where $c^{k}_{ij}$ are the structure constants.

An exponential curve is obtained by using the $\textbf{exp}$ mapping of the Lie algebra elements, i.e\, an exponential curve passing through the identity element $e\in SIM(2)$ at $t=0$ can be written as 
\begin{equation}
\gamma_{c}(t)=\text{exp}\left(t\sum\limits_{i=1}^{4}c^{i}
\mathcal{A}_{i}\bigg{|}_{g=e}\right)=\text{exp}\left( t\sum\limits_{i=1}^{4}c^{i}
A_{i}\right),
\end{equation}
and an exponential curve passing through $g_{0}\in  SIM(2)$ can be obtained by left multiplication with $g_{0}=(x_{0},y_{0},e^{\tau_{0}},\theta_{0})$, i.e.\ $g_{0}\gamma_{c}(t)$. The following theorem applies the method of characteristics (for PDEs) to transport along exponential curves. The explicit formulation  is important because left-invariant convection-diffusion on $SIM(2)$ takes place only along exponential curves, see Theorem~\ref{thm:DiffGeomInt}.
\begin{theorem}\label{thm:Exp_Char}
Let $A\in T_e(SIM(2))$. Then the following holds.
\begin{enumerate}
\item{$U\in\mathcal{D}(d\mathcal{R}(A))\Rightarrow\mathcal{R}_{e^{tA}}U\in\mathcal{D}(d\mathcal{R}(A))$, where $\mathcal{D}(X)$ denotes the domain of operator $X$.}
\item{$e^{t d\mathcal{R}(A)}=\mathcal{R}_{e^{tA}}$, $\forall t>0$ where $d\mathcal{R}$ is defined in \eqref{SIM2_12}. }
\item{$\gamma_{c}(t)=g_{0} \ \text{exp}\left( t\sum\limits_{i=1}^{4}c^{i}
A_{i}\right)$ are  the characteristics for the following PDE,
\begin{align}
\frac{\partial W(g,t)}{\partial t}=-\sum\limits_{i=1}^{4}c^{i}
\mathcal{A}_{i}W(g,t), \ 
W(g,0)=U.
\label{SIM2_14}
\end{align}}
\end{enumerate}
\end{theorem}
\textbf{Proof.}
See \ref{App:Exp_Char} for proof. $\hfill\Box$\\

The exponential map defined on $T_e(SIM(2))$ is bijective, and so we can define the logarithm mapping, $\log=(\exp)^{-1}:SIM(2)\rightarrow T_e(SIM(2))$. For proof see \cite[Chap. 4]{Antoine2004}. For the explicit formulation of the exponential and logarithm curves in our case see \ref{App:ExpCurvesExplicit}.
The explicit form of the $\log$ map will be used to approximate the solution for linear evolutions on $SIM(2)$ in Section~\ref{sec:GausEstSIM}.
\subsection{Quadratic forms on Left Invariant vector fields}\label{sec:OperatorsOnScores}
We apply the general theory of evolutions (convection-diffusion) on Lie groups, \cite{Duits2007b}, to the $SIM(2)$ group and consider the following left-invariant second-order evolution equations,
\begin{equation}
\begin{cases}
\partial _{t}W(g,t)=Q^{\mathbf{D},\mathbf{a}}(\mathcal{A}_{1},\mathcal{A}_{2},\mathcal{A}_{3},\mathcal{A}_{4})W(g,t),\\
W(\cdot,t=0)=\mathcal{W}_{\psi}f(\cdot),
\end{cases}
\label{Evo_0}
\end{equation}    
where $W:SIM(2)\times \mathbb{R}^{+}\rightarrow\mathbb{C}$ and $Q^{\mathbf{D},\mathbf{a}}$ is given by
\begin{equation}
Q^{\mathbf{D},\mathbf{a}}(\mathcal{A}_{1},\mathcal{A}_{2},\mathcal{A}_{3},\mathcal{A}_{4})=\sum\limits_{i=1}^{4}\left(
-a_{i}\mathcal{A}_{i}+ \sum\limits_{j=1}^{4}D_{ij}\mathcal{A}_{i} \mathcal{A}_{j}
\right), \  \ a_{i},D_{ij}\in\mathbb{R}, \ \mathbf{D}:=[D_{ij}]\geq 0, \ \mathbf{D}^{T}=\mathbf{D}.
\label{Evo_5}
\end{equation}
Throughout this article we restrict ourselves to the diagonal case with $a_1=a_2=a_4=0$. This is a natural choice when $a_{i}$ and $D_{ij}$ are constant, as we do not want to impose a-priori curvature and a-priori scaling drifts in our flows. However, when adapting $D$ and $a$ to the initial condition (i.e. wavelet transform data) such restrictions are not necessary. In fact, practical advantages can be obtained when choosing $D$ diagonal w.r.t. optimal gauge frame, \cite{Franken2009,Duits2010a}. Choosing  $D_{ij}=D_{ii}\delta_{ij}, \ i,j\in\{1,2,3,4\}$, the quadratic form becomes,
\begin{equation}
Q^{\mathbf{D},\mathbf{a}}(\mathcal{A}_{1},\mathcal{A}_{2},\mathcal{A}_{3},\mathcal{A}_{4})=[-a_{1}\partial _{\theta}-a_{2}\partial _{\xi}-a_{3}\partial _{\eta}-a_{4}\partial _{\tau}+D_{11}(\partial _{\theta})^{2}+D_{22}(\partial _{\xi})^{2}+D_{33}(\partial _{\eta})^{2}+D_{44}(\partial _{\tau})^{2}].
\label{Evo_1}
\end{equation}
The first order part of \eqref{Evo_1} takes care of transport (convection) along the exponential curves, deduced in Section \ref{Sec:LeftInvVecFiel}. The second order part takes care of diffusion in the $SIM(2)$ group. Note that these evolution equations are left-invariant by construction.

%

H\"{o}rmander in \cite{Hormander1967} gave necessary and sufficient conditions on the convection and diffusion parameters, respectively $a = (a_1, a_2, a_3,a_4)$ and $\mathbf{D}=D_{ij}$, in order to get smooth 
Green's functions of the left-invariant convection-diffusion equation \eqref{Evo_0} with generator \eqref{Evo_1}. By these conditions the non-commutative nature of $SIM(2)$ in certain cases takes
care of missing directions in the diffusion tensor. Applying the H\"{o}rmander's theorem \cite{Hormander1967} produces
%
necessary and sufficient conditions for smooth (resolvent) Green's functions on $SIM(2)\backslash\{e\}$ on the
diffusion and convection parameters $(\mathbf{D},a)$ in the generator \eqref{Evo_1} of \eqref{Evo_0} for diagonal 
$\mathbf{D}$:
\begin{align}
\{1,2,4\}\subset\{i| \ a_i\neq 0 \vee D_{ii}\neq 0\} \vee \{1,3,4\}\subset\{i| \  a_i\neq 0 \vee D_{ii}\neq 0\}.
\end{align}

A covariant derivative of a co-vector 
field $\mathbf{a}$ on the manifold $(SIM(2),
\mathcal{G})$ is a $(0,2)$-tensor field with 
components 
$\nabla_{j}a_i=\mathcal{A}_ja_i-\Gamma^k
_{ij}a_k$, whereas the covariant derivative 
of a vector field $\mathbf{v}$ on $SIM(2)$ 
is a $(1,1)$-tensor field with components 
$\nabla_{j'}v^i=\mathcal{A}_{j'}v^i+
\Gamma^i_{j'k'}v^{k'}$, where we have made use of the notation $\nabla_j:= D_{\mathcal{A}_j}$, when imposing the Cartan connection (tangential to the $SE(2)$-case in \cite{Duits2010a}, $SE(3)$-case in \cite{Duits2011} and the $H(2d+1)$-case in \cite{Duits2012}).
The Christoffel symbols equal minus the (anti-symmetric) structure constants of the Lie algebra $\mathcal{L}(SIM(2))$, i.e.\ $\Gamma_{ij}^k=-c_{ij}^k$. 
The left-invariant equations \eqref{Evo_0} with a diagonal diffusion tensor \eqref{Evo_1} can be rewritten in covariant derivatives as
\begin{align}
\begin{cases}
&\partial_sW(g,s)=\sum\limits_{i,j=1}^4
\mathcal{A}_{i}((D_{ij}(W))(g,s)\mathcal{A}_j W)(g,s)=
\sum\limits_{i,j=1}^4
\nabla_{i}((D_{ij}(W))(g,s)\nabla_j W)(g,s),\\
&W(g,0)=\mathcal{W}_{\psi}f(g), \text{ for all } g\in SIM(2), \ s>0.
\end{cases}
\label{Evo_40}
\end{align}
Both convection and diffusion in the left-invariant evolution equations \eqref{Evo_0}
take place along the exponential curves in $SIM(2)$ which are covariantly constant curves with respect to the Cartan connection.
For proof and various details see \ref{sec:App_DiffGeomInter}.
\subsection{Probabilistic models for contextual feature propagation}\label{Stoc_Conn}  
Section \ref{sec:OperatorsOnScores} described the general form of convection-diffusion operators on the $SIM(2)$ group. For the particular case of contour enhancement i.e.\ diffusion on the $SIM(2)$ group, which corresponds to the choice $D_{ij}=D_{ii}\delta_{ij}, \ i,j\in\{1,2,3,4\}, \ D_{33}=0 \text{ and } \mathbf{a}=\mathbf{0}$, we have the following result.  
\begin{theorem}\label{thm:StochasticConn}
The evolution on $SIM(2)$ given by
\begin{equation}
\begin{cases}
\partial _{t}W(g,t)=[D_{11}(\partial_\theta)^2+D_{22}(\partial_\xi)^2+D_{44}(\partial_\tau)^2]W(g,t),\\
W(\cdot,t=0)=\mathcal{W}_{\psi}f(\cdot),
\end{cases}
\label{Eq:ContourEnhGen}
\end{equation}    
is the forward Kolmogorov (Fokker-Planck) equation of the following stochastic process for multi-scale contour enhancement
\begin{align}
\begin{cases}
\mathbf{X}(s)=\mathbf{X}(0)+\sqrt{2 D_{22}}\epsilon_2\int\limits_0^s\left( \cos(\Theta(t))e_x+\sin(\Theta(t)e_y) e^{\mathfrak{T}(t)}d(\sqrt{t}
\right)\\
\Theta(s)=\Theta(0)+\sqrt{s}\sqrt{2 D_{11}}\epsilon_1\\
\mathfrak{T}(s)=\mathfrak{T}(0)+\sqrt{s}\sqrt{2 D_{44}}\epsilon_4, \label{eq:SDE}
\end{cases}
\end{align} 
where $e_1,  e_2,  e_4\sim \mathcal{N}(0,1)$ are the standard random variables  and  $D_{11},D_{22},D_{44}>0$. 
\end{theorem}
In order to avoid technicalities regarding probability measures on Lie groups, see \cite{Hsu2002} for details, we only provide a short and basic explanation which covers the essential idea of the proof.  

The stochastic differential equation in  \eqref{eq:SDE} can be considered as limiting case of the following discrete stochastic processes on $SIM(2)$:
\begin{align}
\begin{cases}
G_{n+1}:=(\mathbf{X}_{n+1},\Theta_{n+1},\mathfrak{T}_{n+1})=G_n+\sqrt{\Delta s}\sum 
\limits_{i=1,2,4}\sum\limits_{j=1,2,4
}\frac{\epsilon_{i,n+1}}{\sqrt{N}}\sqrt{2D_{ii}} \  e_i\big{|}_{G_N},\\
G_0=(\mathbf{X}_0,\Theta_0,\mathfrak{T}_0),
\end{cases}\label{eq:DiscreteSDE}
\end{align}
where $n=1,\cdots,N-1,N\in\mathbb{N}$ denotes the number of steps with step-size $\Delta s>0$, $\{{\epsilon_{i,n+1}}\}_{i=1,2,4}$ are independent normally distributed $\epsilon_{i,n+1}\sim \mathcal{N}(0,1)$ and $e_j\big{|}_{G_n}\equiv\mathcal{A}_j
\big{|}_{G_n}$, i.e.
\begin{align*}
e_1\big{|}_{G_n}=\begin{pmatrix}
0 \\ 0\\ 0\\ 1
\end{pmatrix}, \ e_2\big{|}_{G_n}=\begin{pmatrix}
e^{\mathfrak{T}}\cos\Theta \\ e^{\mathfrak{T}}\sin\Theta\\ 0\\ 0
\end{pmatrix}, \ e_3\big{|}_{G_n}=\begin{pmatrix}
-e^{\mathfrak{T}}\sin\Theta \\ e^{\mathfrak{T}}\cos\Theta\\ 0\\ 0
\end{pmatrix}, \ e_4\big{|}_{G_n}=\begin{pmatrix}
0 \\ 0\\ 1\\ 0
\end{pmatrix}.
\end{align*} 
Note that the continuous process \eqref{eq:SDE} directly arises from the discrete process \eqref{eq:DiscreteSDE} by recursion and taking the limit $N\rightarrow\infty$. 
\section{Left-invariant Diffusions on $SIM(2)$}
Following our framework of stochastic left-invariant evolutions on $SIM(2)$ we will restrict ourselves to contour enhancement, where the Forward-Kolmogorov equation is essentially a hypo-elliptic
diffusion on the $SIM(2)$ group and therefore we recall \eqref{Eq:ContourEnhGen}
\begin{equation}
\begin{cases}
\partial _{t}W(g,t)=[D_{11}(\partial_\theta)^2+D_{22}(\partial_\xi)^2+D_{44}(\partial_\tau)^2]W(g,t),\\
W(\cdot,t=0)=\mathcal{W}_{\psi}f(\cdot).
\end{cases}
\label{Evo_25}
\end{equation}    
In the remainder of this paper we study linear diffusion (combined with monotone operations on the co-domain) and non-linear diffusion on the $SIM(2)$ group in the context of our imaging application.
\subsection{Approximate Contour Enhancement Kernels for Linear Diffusion on Scale-OS}\label{sec:GausEstSIM}
In \cite{Duits2008,Duits2010}, the authors derive the exact Green's function of \eqref{Evo_0} for the $SE(2)$ case. To our knowledge explicit and exact formulae for heat kernels of linear diffusion on $SIM(2)$ do not exist in the literature. However using the general theory in \cite{Nagel1990,Elst1998}, one can compute Gaussian estimates for Green's function of left-invariant diffusions on Lie groups. As a first step this involves approximating  $SIM(2)$ by a parametrized class of groups $(SIM(2))_q, \ q\in[0,1]$ in between $SIM(2)$ and the nilpotent Heisenberg approximation $SIM(2)_0$. This idea of contraction has been explained in \ref{App:MethContraction}.

According to the general theory in \cite{Elst1998}, the heat-kernels $K_{t}^{q,\mathbf{D}}:(SIM(2))_{q}\rightarrow \mathbb{R}^{+}$ (i.e.\ kernels for contour enhancement whose convolution yields diffusion on $SIM(2))_{q}$) on the parametrized class of groups $(SIM(2))_{q}, \ q\in [0,1]$ in between $SIM(2)$ and $(SIM(2))_{0}$ satisfy the Gaussian estimates
\begin{equation}
|K_{t}^{q,\mathbf{D}}(g)|\leq C t^{-\frac{5}{2}} \exp \left(\frac{-b\| g\| ^{2}_{q}}{4t}\right), \text{  with $C, b>0$ (constant), $g\in (SIM(2))_{q}$ },
\label{Evo_6}
\end{equation}
\begin{figure}[t]
\centering
\includegraphics[scale=0.6]{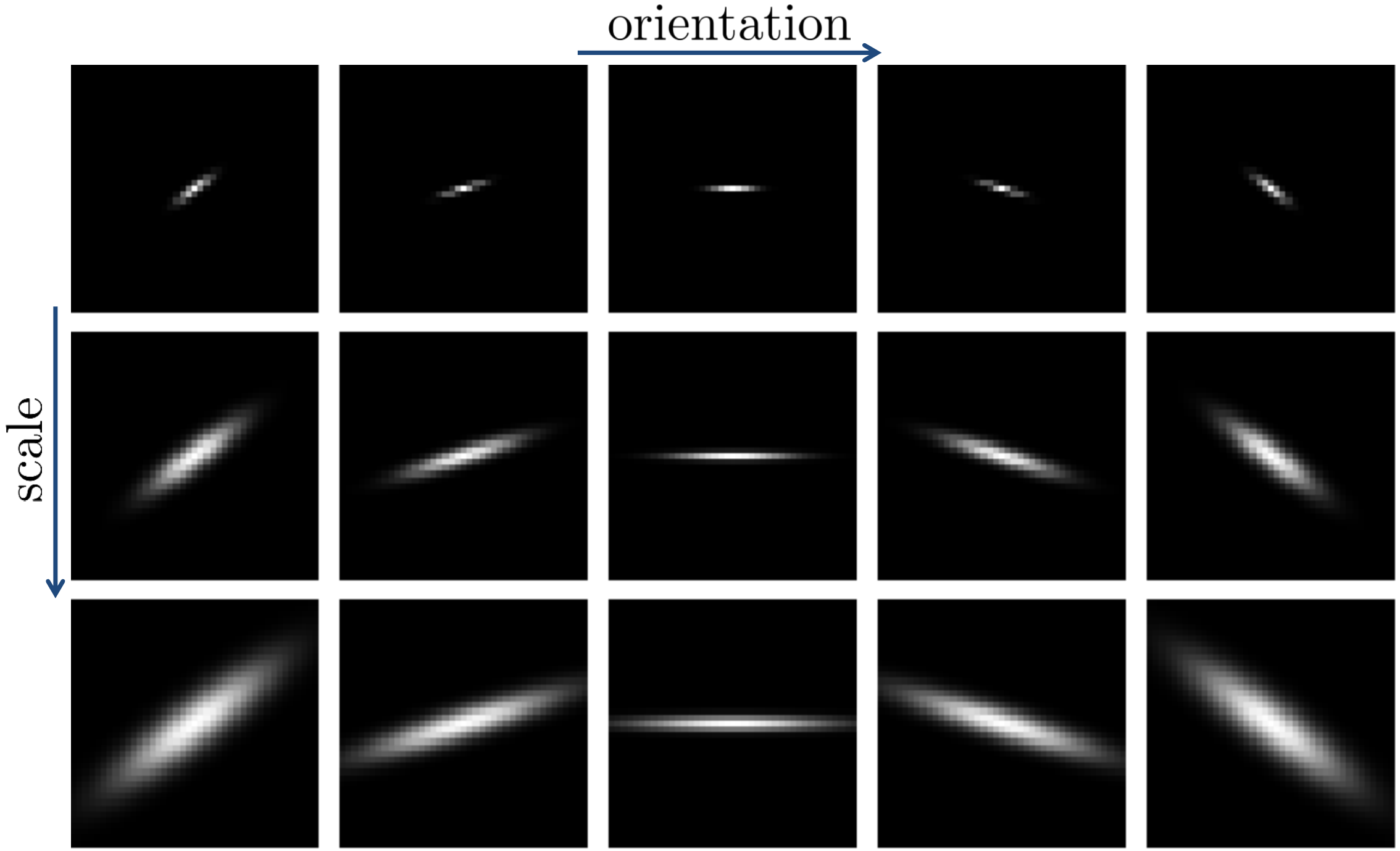}
\caption[Enhancement kernel via Gaussian estimates at different scales and orientations]{Plots of enhancement kernels generated on a $41\times 41$ grid at different scales and orientations. Parameters chosen: 
$D_{11}=0.05, \ D_{22}=1, \ D_{44}=0.01, \
t=0.7$. }
\label{fig:EnhancementKernels_Scales}
\end{figure}

where the norm $\|\cdot\|_{q}:(SIM(2))_{q}\rightarrow\mathbb{R}^+$ is given by $\|g\|_q=|\log_{(SIM(2))_q}(g)|_{q}$. $\log_{(SIM(2))_q}:(SIM(2))_q\rightarrow T_{e}((SIM(2))_q)$, is the logarithmic mapping on $(SIM(2))_q$, which is computed explicitly for the case of $(SIM(2))_{q=1}=SIM(2)$ in \ref{App:ExpCurvesExplicit}, and where the weighted modulus, see \cite[Prop 6.1]{Elst1998}, in our special case of interest is given by,
\begin{align}
\bigg{|}
\sum\limits_{i=1}^{4}c_{q}^{i}A_{i}^{q}
\bigg{|}_q=
\sqrt{|c_{q}^1|^{2/w_1}+|c_{q}^2|^{2/w_2}+|c_{q}^3|^{2/w_3}+|c_{q}^4|^{2/w_4}}=
\sqrt{((c_{q}^1)^{2}+(c_{q}^2)^{2}
+(c_{q}^4)^{2})+|c_{q}^3|},
\end{align}
where recall the weightings from \eqref{Evo_3}. 
\begin{remark}
The constants $b, \ C$ in \eqref{Evo_6} can be taken into account by the transformations  
$t\mapsto\frac{t}{b}, \ f\mapsto Cf$. 
\end{remark}
A sharp estimate for the front factor constant $C$ in \eqref{Evo_6} is given by 
\begin{align}
C=\frac{1}{4\pi D_{11}D_{22}}\frac{1}{\sqrt{D_{44}}}.
\end{align}
This follows from the general theory in \cite{Elst1998} where the choice of constant $C$ is uniform for all groups $(SIM(2))_q, \ q\in[0,1]$. Thus to determine $C$ we need to determine the front factor constant for the green's function $G$ of the following resolvent equation
\begin{align*}
\left((D_{11}(\mathcal{A}_1^0)^2+D_{22}(\mathcal{A}_2^0)^2+D_{44}(\mathcal{A}_4^0)^2)\right)G(x,y,\tau,\theta)=+\delta_{e}
\end{align*}
where $\mathcal{A}^0$ denote the basis for $\mathcal{L}(H)$ (recall $H=(SIM(2))_{q\downarrow 0}$) and $e$ is the identity of $H$. Using $[D_{11}(\mathcal{A}^0_1)^2+D_{22}(\mathcal{A}^0_2)^2,(\mathcal{A}_4^0)^2]=0$ and the results in \cite{Duits2010} 
%
arrive at
\begin{equation}
|K_{t}^{q=1,\mathbf{D}}(g)|\leq \frac{1}{4\pi t^{\frac{5}{2}} D_{11}D_{22}\sqrt{D_{44}}}\exp\left(
\frac{-1}{4t}\left(
\frac{\theta^2}{D_{11}} + \frac{(c^2(g))^2}
{D_{22}}+\frac{\tau ^2}{D_{44}}+\frac{|c^3(g)|}
{\sqrt{D_{11}D_{22}D_{44}}}
\right)
\right),
\end{equation}
where recall the explicit formulae for $c^{i}, \ i\in\{2,3\}$  from \ref{App:ExpCurvesExplicit},
\begin{align*}
c^2(g)=\frac{(y \theta -x \tau) +(-\theta\eta+\tau\xi)}{t \left(1+e^{2 \tau }-2 e^{\tau } cos\theta\right)}, \ 
c^3(g)=\frac{-(x \theta +y \tau) +(\theta\xi+\tau\eta)}{t \left(1+e^{2 \tau }-2 e^{\tau } cos\theta\right)}.
\end{align*}
A problem with these estimates is that they are not differentiable everywhere. This problem can be solved by using the estimate 
\begin{align*}
|a|+|b|\geq\sqrt{a^2+b^2}\geq\frac{1}{\sqrt{2}}(|a|+|b|),
\end{align*}
which holds for all $a,b\in \mathbb{R}$, to the exponents of our Gaussian estimates. Thus we estimate the weighted modulus by the equivalent (for all $q>0$) weighted modulus $|\cdot|_{q}:T_{e}((SIM(2))_q)\rightarrow\mathbb{R}^+$ by $\bigg{|}
\sum\limits_{i=1}^{4}c_{q}^{i}A_{i}^{q}
\bigg{|}_{q}:=\sqrt[4]{((c_{q}^1)^{2}+(c_{q}^2)^{2}
+(c_{q}^4)^{2})^2+|c_{q}^3|^2}$, yielding the Gaussian estimate,
\begin{equation}\label{eq:Gaussian_Est}
|K_{t}^{q=1,\mathbf{D}}(g)|\leq \frac{1}{4\pi t^{\frac{5}{2}} D_{11}D_{22}\sqrt{D_{44}}}\exp\left(
\frac{-1}{4t}\left(
\left[\frac{\theta^2}{D_{11}} + \frac{(c^2(g))^2}
{D_{22}}+\frac{\tau ^2}{D_{44}}\right]^2+\frac{|c^3(g)|^2}
{D_{11}D_{22}D_{44}}
\right)
\right).
\end{equation}

\begin{figure*}[t]
        \centering
        \begin{subfigure}[t]{0.28\textwidth}
          \centering
           \includegraphics[width=\textwidth]{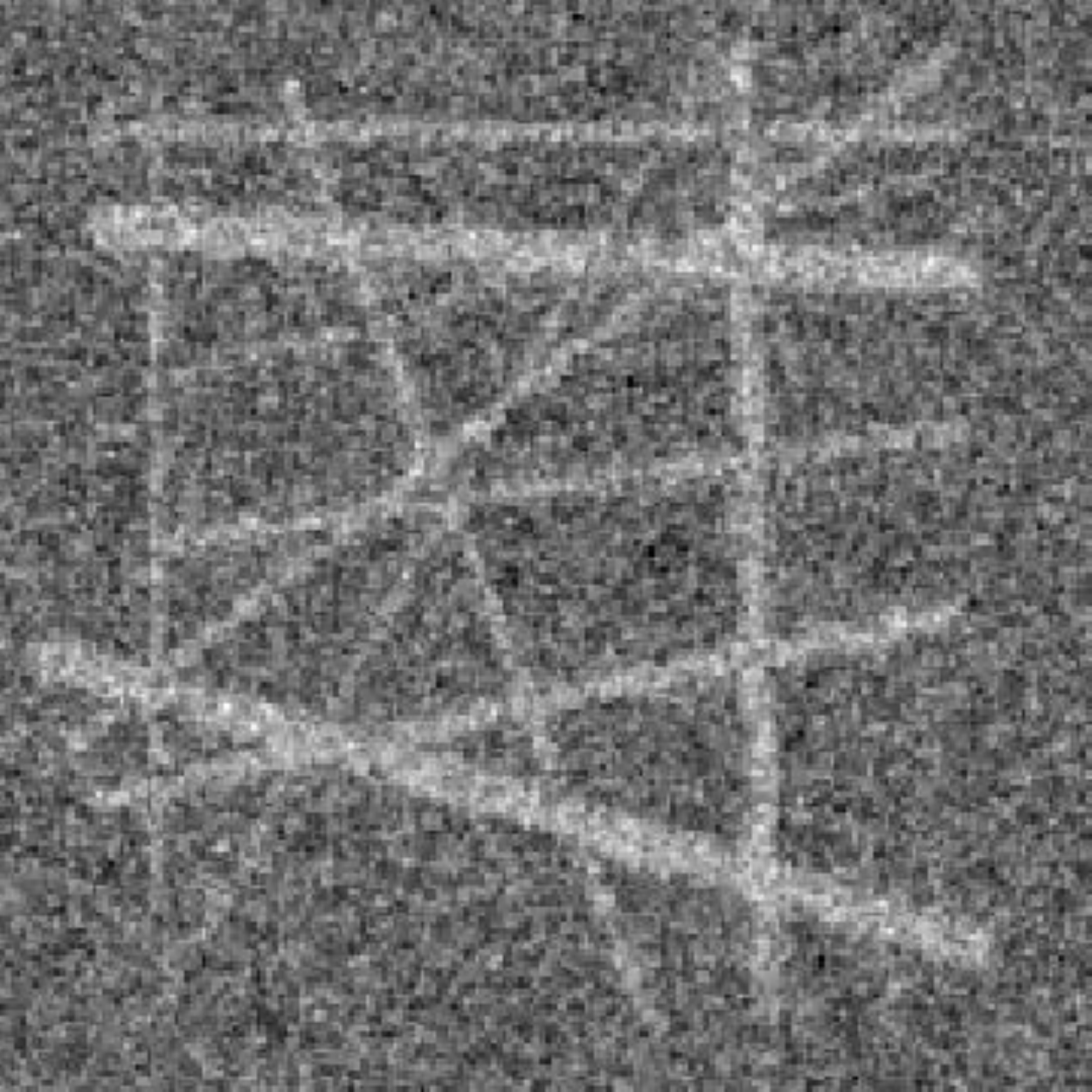}
            \caption{Noisy input image}
        \end{subfigure}
        \begin{subfigure}[t]{0.28\textwidth}
                \centering
                \includegraphics[width=\textwidth]{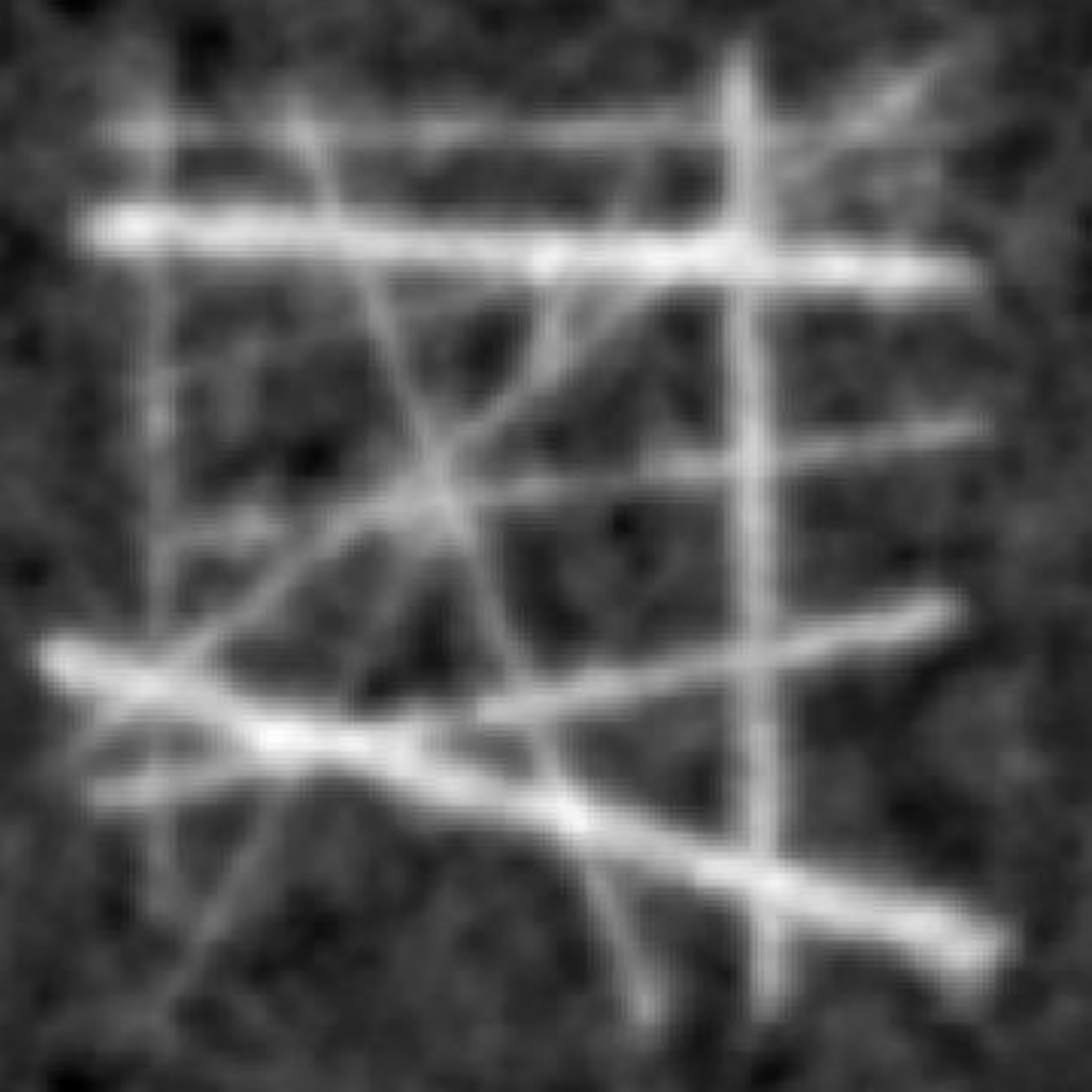}
                \caption{$SE(2)$ Linear Diffusion}
        \end{subfigure}
\begin{subfigure}[t]{0.28\textwidth}
                \centering
                \includegraphics[width=\textwidth]{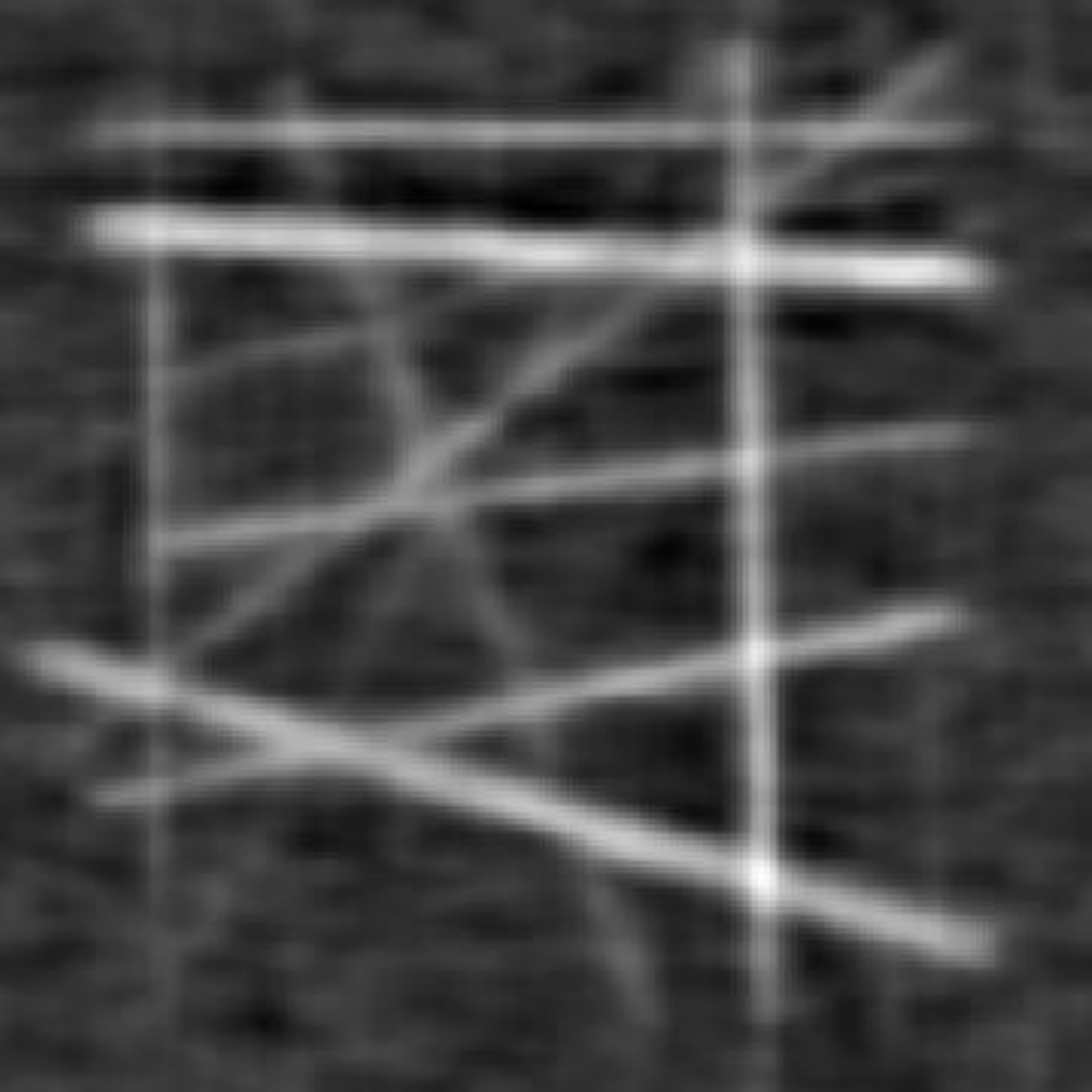}
                \caption{$SIM(2)$ Linear Diffusion}
        \end{subfigure}
        \caption
{Comparison of enhancement via linear diffusion using Gaussian estimates for heat kernels on $SE(2)$ and $SIM(2)$.
Left:Original image, Center: Enhancement via Linear diffusion on $SE(2)$ group using Gaussian estimates in \cite{Duits2010}, Right: Enhancement via Linear diffusion on $SIM(2)$ group using estimate in \eqref{eq:Gaussian_Est}. 
Parameters- OS: $N_\theta=20,N_s=5$; $SE(2)$ estimate: $D_{11}=0.05, \ D_{22}=1, \
t=3$ ;$SIM(2)$ estimate: $D_{11}=0.05, \ D_{22}=1, \ D_{44}=0.01, \
t=0.7$.} 
\label{fig:LinearDiffusionComp}
\end{figure*}

Figure \ref{fig:EnhancementKernels_Scales}
shows the typical structure of these enhancement kernels. 
\begin{remark}
When cascading group convolutions and transformations in the co-domain of the scores one can generalize the scattering operators by Mallat et al. on $\mathbb{R}^n$, see \cite{Mallat2012,Sifre2013}, to left-invariant scattering operators on affine Lie groups such as $SIM(2)$ which would provide us with stability under local deformations. 
\end{remark}
In practice, medical images exhibit complicated structures which require local adaptivity per group location via gauge frames. This brings us to non-linear diffusions that we will solve numerically in the next section.
\subsection{Nonlinear Left Invariant Diffusions on $SIM(2)$}
\begin{figure}[t]
\centering
\includegraphics[scale=0.4]{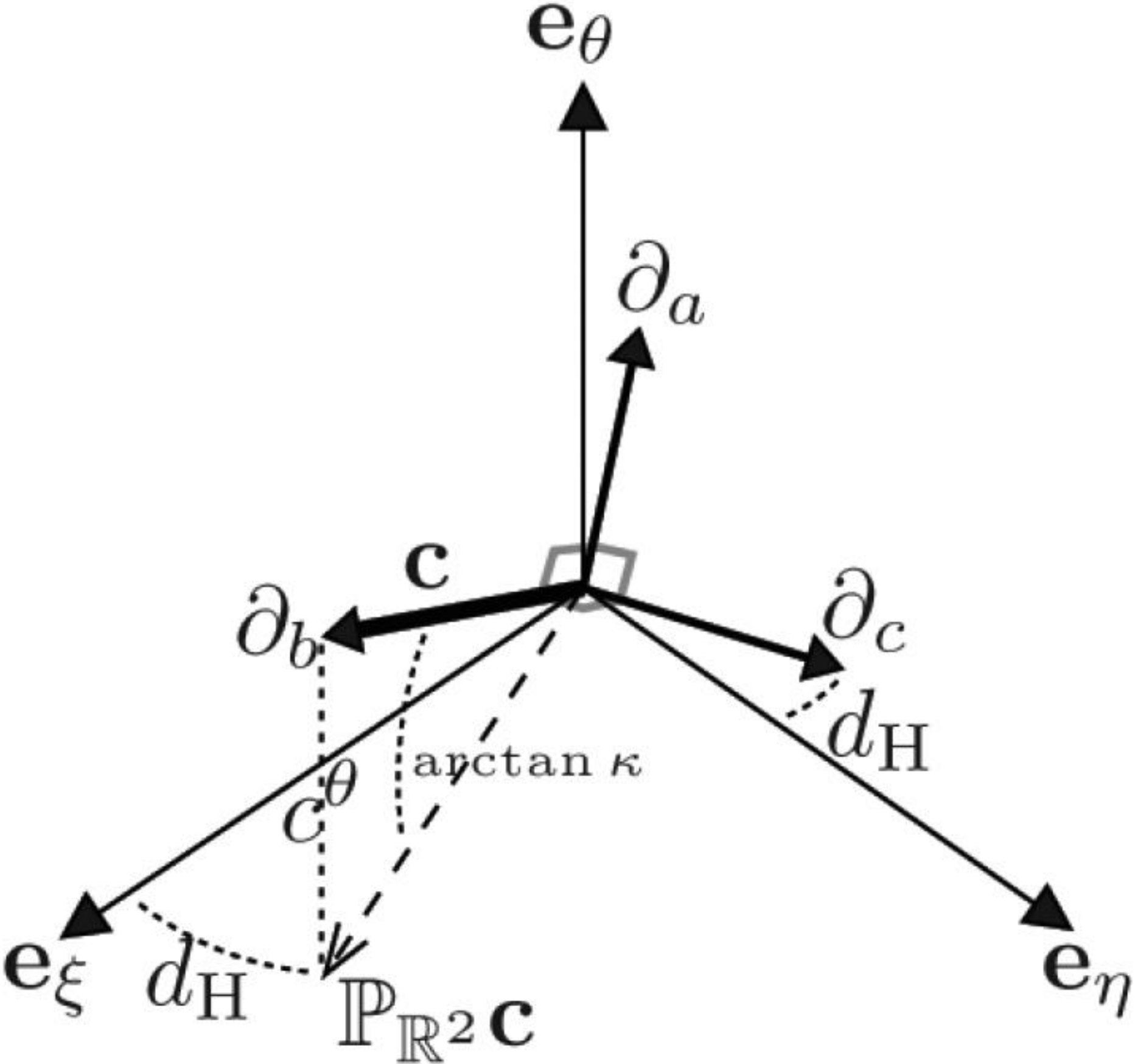}
\caption
{Illustration of curvature $\kappa$ and deviation from horizontality $d_H$ and the gauge frame \eqref{eq:GaugeFrames}. Note that $e_\theta\leftrightarrow\partial_\theta, \ e_\xi\leftrightarrow\partial_\xi, \ e_\eta\leftrightarrow\partial_\eta$.} 
\label{fig:CEDOS_QAM}
\end{figure}
Adaptive nonlinear diffusion on the $2D$ Euclidean motion group $SE(2)$ called coherence enhancing diffusion on orientation scores (CED-OS)  was introduced in \cite{Franken2009,Duits2010}. We wish apply this adaptive $SE(2)$ diffusion to each scale in our scale-OS, which is possible because at a fixed scale the scale-OS is a function on the $SE(2)$ group. Next we present a brief outline of the CED-OS algorithm and then apply it to our case of interest. 

\subsubsection{CED-OS - Brief Outline}
CED-OS involves the following two steps:
\begin{itemize}
\item{\textit{Curvature Estimation.} Curvature estimation of a spatial curve is based on the optimal exponential curve fit at each point. We find the best exponential curve fit (and the corresponding curvature) to the data $(x,y,\theta)\mapsto\left|(\mathcal{W}_\psi f)(x,y,\theta)\right|$ by using the exponential map,
\begin{align*}
t\mapsto \exp\big{(}t(c^1_*A_1+c^2_*A_2+c^3_*A_3)\big{)} \text{ with } (c^1_*)^2+\beta^2(c^2_*)^2+\beta^2(c^3_*)^2=1,
\end{align*}
via the techniques explained in \cite{Duits2010a,Franken2009} and summarized in \ref{App:CurvatureEstimation}.}
\item{\textit{Adaptive curvatures based diffusion scheme using gauge coordinates.} The best exponential curve fit mentioned above is parametrized by $\mathbf{c}_*=(c^1_*,c^2_*,c^3_*)\in\mathbb{R}^3$ which provides us the curvature $\kappa$ (and deviation from horizontality $d_H$ if we do not impose $c^3_*=0$). In fact it furnishes a whole set of gauge frames $\{\partial_a,\partial_b,\partial_c\}$ as can be seen in Figure~\ref{fig:CEDOS_QAM}. The gauge frames in spherical coordinates read
\begin{align}\label{eq:GaugeFrames}
\partial_a&=-\cos\alpha\cos d_H\partial_\xi-\cos\alpha\sin d_H\partial_\eta+\beta\sin\alpha\partial_\theta,\nonumber\\
\partial_b&=\sin\alpha\cos d_H\partial_\xi+\sin\alpha\sin d_H\partial_\eta+\beta\cos\alpha\partial_\theta,\\
\partial_a&=-\sin d_H\partial_\xi+\cos d_H\partial_\eta\nonumber.
\end{align}
The resulting nonlinear evolution equations on orientation scores is
\begin{align}
\begin{cases}
\partial_t U(g,t)=\begin{pmatrix}\beta\partial_\theta & \partial_\xi &  \partial_\eta\end{pmatrix}M^T_{\alpha,d_H}
 \begin{pmatrix}
  D_{aa} & 0 &  0 \\
  0 & D_{bb} & 0 \\
  0 & 0 &  D_{cc}
 \end{pmatrix}M_{\alpha,d_H}\begin{pmatrix}\beta
 \partial_\theta \\ \partial_\xi \\  \partial_\eta\end{pmatrix}U(g,t), \ t>0,
\\
U(g,t=0)=\mathcal{W}_\psi[f](g) \text{ for all } g\in SE(2),
\end{cases}\label{Evo_27}
\end{align}
where we use the shorthand notation $D_{ii}=(D_{ii}(U))(g,t)$, for $i=a,b,c$. The matrix 
\begin{align*}
M_{\alpha,d_H}=\begin{pmatrix}
\sin\alpha & -\cos\alpha\cos d_H &  -\cos\alpha\sin d_H \\
\cos\alpha & \cos d_H\sin\alpha &  \sin\alpha\sin d_H \\
  0 & -\sin d_H &  \cos d_H
 \end{pmatrix}
\end{align*}
is the rotation matrix in $SO(3)$ that maps the left-invariant vector fields $\{\beta\partial_\theta,\partial_\xi,\partial_\eta\}$ onto the gauge frame $\{\partial_a,\partial_b,\partial_c\}$.
} 
\end{itemize}
\begin{figure}[t]
\centering
\includegraphics[scale=0.7]{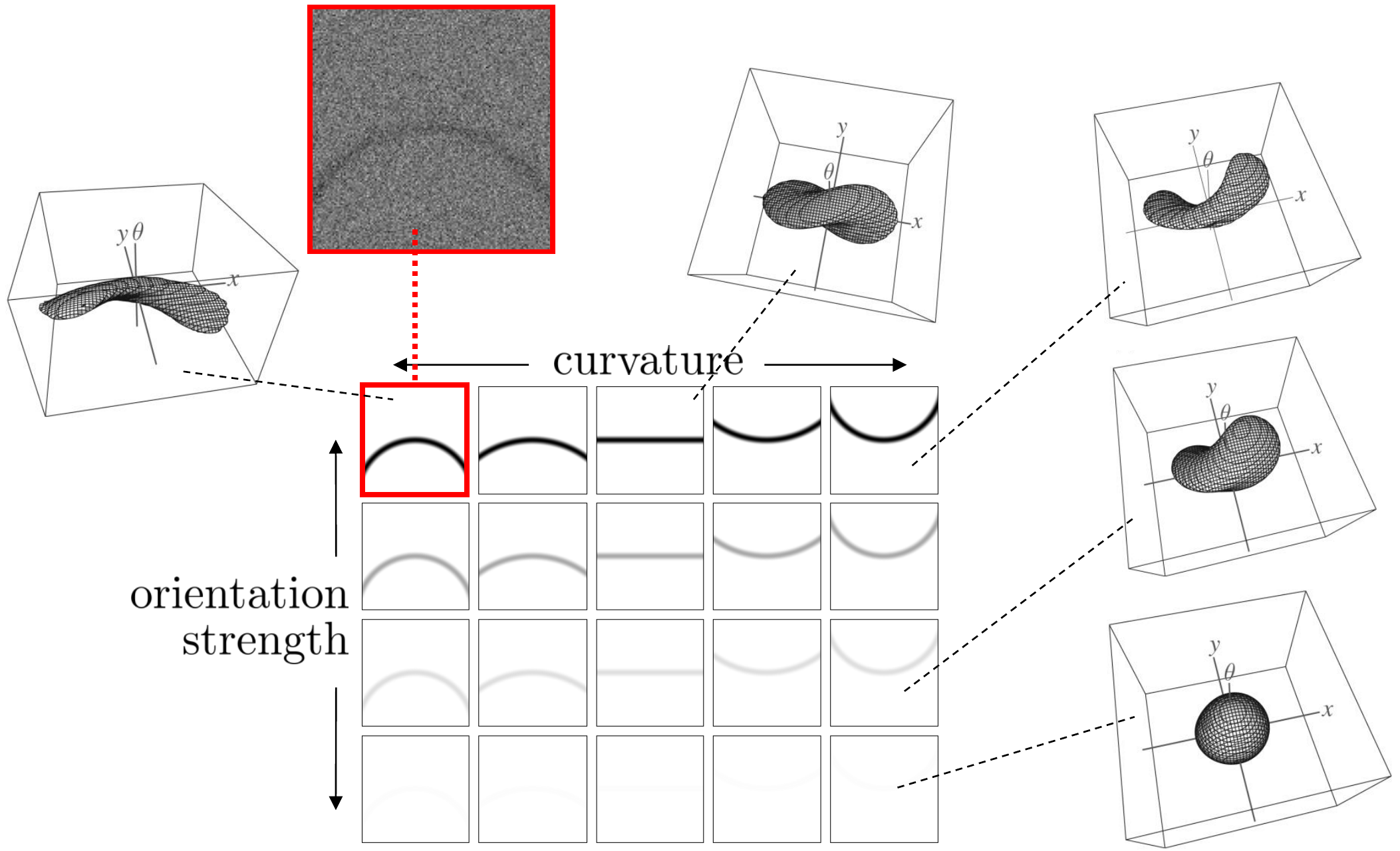}
\caption
{Illustration of heat kernels $K_t^D:SE(2)\rightarrow \mathbb{R}^+$ on $SE(2)$ corresponding to varying curvature and orientation strengths. As orientation strength decreases the kernels become more isotropic (needed for noise reduction) and as curvature increases the diffusion kernels in $SE(2)$ bend accordingly.} 
\label{fig:CurvatureEffect}
\end{figure}
Figure~\ref{fig:CurvatureEffect} depicts the dependence of the corresponding local linear diffusion kernel on curvature and orientation strengths. Eq.\eqref{Evo_27} is implemented  by a Euler  forward scheme (expressed by finite differences along the moving frame of reference while relying on B-spline approximations yielding better rotation covariance \cite{Franken2008}[Pg.142]) involving the parameters $(\widetilde{\rho},\rho_s,c,\beta)$, 
\begin{align*}
U_{n+1}&=\Delta t (Q U_n)+U_n,  
\end{align*}
where  $U_n=U(\cdot,\cdot,t_n)$  with discrete time steps 
and the generator is given by,
\begin{small}
\begin{align*}
(QU)(\cdot,\cdot,t)&=
\left(\beta\partial_\theta^D \ \partial_\xi^D \  \partial_\eta^D\right) M^T_{\alpha,d_H}
 \left(G_{\rho_s,\rho_s,\rho_s\beta^2}\ast\begin{pmatrix}
  D_{aa}(U(\cdot,\cdot,t)) & 0 &  0 \\
  0 & D_{bb}(U(\cdot,\cdot,t)) & 0 \\
  0 & 0 &  D_{cc}(U(\cdot,\cdot,t))
 \end{pmatrix}\right)M_{\alpha,d_H}\begin{pmatrix}\beta
 \partial_\theta^D \\ \partial_\xi^D \\  \partial_\eta^D\end{pmatrix},
\end{align*}
\end{small}

where $G_{\rho_s,\rho_s,\rho_s\beta^2}(x,y,\theta)$ denotes a Gaussian kernel, which is isotropic on each spatial plane with spatial scale $\rho_s=\frac{1}{2}\sigma_s^2$ and anisotropic with scale $\rho_s\beta^2$ in orientation. Here $\partial^D_j$ denotes the left-invariant finite difference operator in the $j^{th}$ direction implemented via B-spline interpolation \cite{Franken2009}. Furthermore we set $D_{bb}=0$, $(D_{aa}U)(\cdot,\cdot,t)=(D_{cc}U)(\cdot,\cdot,t)=e^{\frac{-(s(|U|)(\cdot,\cdot,t))^2}{c}}, \ c>0$, with the orientation confidence $s(U(\cdot,\cdot,t))=(\partial_a^2+\partial_b^2)U(\cdot,\cdot,t)$. The derivatives $\partial_a,\partial_b$ are Gaussian derivatives, computed with scales $(\widetilde{\rho},\widetilde{\rho},\widetilde{\rho}\beta^2)$, orthogonal to locally optimal direction with $c$ being the eigenvector of the Hessian matrix corresponding to best exponential curve fit, where we enforce horizontality \cite{Franken2009} by setting $d_H=0$ (see Figure~\ref{fig:CurvatureEffect}).
\subsubsection{CED-OS on Scale-OS (CED-SOS)}\label{sec:CED_SOS}
For an image $f$ in $\mathbb{L}^2(\mathbb{R}^2)$ the corresponding scale-OS is within $\mathcal{W}_{\psi}(f)\in \mathbb{C}_{K}^{SIM(2)}\subset \mathbb{L}^{2}(SIM(2))$ and for any fixed $a\in (0,\infty)$, $\mathcal{W}_{\psi}(f)(\cdot,\cdot,a,\cdot)\in\mathbb{L}_{2}(SE(2))$. Let $\Phi_t:\mathbb{L}_{2}(SE(2))\rightarrow \mathbb{L}_{2}(SE(2))$ 
denote nonlinear adaptive diffusion (CED-OS) on the $SE(2)$ group which is the solution operator of \eqref{Evo_27} at stopping time $t$. We propose the operator $\Lambda$ on scale-OS defined as
\begin{align}\label{NonLinDiffOp}
(\Lambda[(\mathcal{W}_\psi f])(x,y,a_i,\theta)=\sum \limits_{i=1}^m(\Phi_{t_{a_i}}[\mathcal{W}_{\psi}f (\cdot,\cdot,a_i,\cdot)])(x,y,\theta) ,  
\end{align}  
where $(x,y,\theta)\in \mathbb{R}^2\times [0,2\pi)$ and $\{a_{i}\}_{i=1}^m$ is the discretization of $[a^-,a^+]$. Recall the defintion of $0<a^-<a^+<\infty$ from Section \ref{sec:DiscreteAnalogue}. Here we make the specific choice of $t_{a_i}$ with $t_{a_i}\leq t_{a_k}$ where $1\leq i<k\leq m$. The idea here is that on lower scales we have to diffuse more as noise is often dominant at lower scales and therefore lower scales need higher diffusion time. 
\subsection{Employing the gauge frames for differential invariany features - Vesselness}\label{sec:SIM-vesselness}
Retinal vasculature images are highly useful for non-invasive investigation of the quality of the vascular system which is affected in case of diseases such as (diabetic retinopathy, glaucoma, hypertension, arteriosclerosis, Altzheimer’s disease), see \cite{Patton2006}. Here accurate and robust detection of the vascular tree is of vital importance \cite{Budai2009,Abramoff2010,Xu2011,
Bankhead2012,Bekkers2012}. The vascular tree structure in these retinal images is often hard to detect due to low-contrast, noise, the presence of tiny vessels, mirco-bleedings and other abnormalities due to diseases,  bifurcations, crossings, occurring at multiple scales, and that is precisely where our continuous wavelet framework on SIM(2) (i.e. multiple scale orientation scores) and the left-invariant evolutions acting upon them comes into play. In this section we will briefly show some benefits of our framework in terms of multi-scale vesselness filtering,  see \cite{Frangi1998,Hannink2014}, enhancement and tracking as illustrated in Figures \ref{fig:Julius-Vesselness-Results} and \ref{fig:Enh-Erik}. We would like to clarify that although $SIM(2)$ vesselness filter summarized below is based on the framework presented in this article, it was developed in \cite{Hannink2014} and has been presented here to demonstrate benefits.

Let $f:\mathbb{R}^2\rightarrow\mathbb{R}$ be a greyscale image. Let $H^sf:=G_s*(Hf)$ denote the Hessian of a Gaussian derivative at scale $s=\sigma^2/2$, with eigensystem $E_{\lambda_1}=\langle\mathbf{e_1}\rangle$, $E_{\lambda_2}=\langle\mathbf{e_2}\rangle$. Then the vesselness filter is given by 
\begin{align*}
\mathcal{V}^{R,S}_{\sigma_1,\sigma_2}(f)=\begin{cases}
e^{-\frac{R^2}{2\sigma_1^2}}\bigg{(}
1-e^{-\frac{S^2}{2\sigma_2^2}}\bigg{)} & \text{ if } \lambda_2>0
\\
0 & \text{ if } \lambda_2\leq 0
\end{cases}
\end{align*}
with anisotropy measure $R=\lambda_1/\lambda_2$ and structureness $S=\lambda_1^2+\lambda_2^2$. Typically $\sigma_1=1/2$ and $\sigma_2=0.2\|S\|_\infty$.\\

Now let us rewrite this filter in gauge derivatives $(\partial_\xi^s,\partial_\eta^s)$ where $\partial_\xi^s$ takes the Gaussian derivative along $\mathbf{e_1}$ at scale $s=\sigma^2/2$ whereas $\partial_\eta^s$ takes the Gaussian derivative along $\mathbf{e_2}$ at scale $s=\sigma^2/2$. Then the matrix representation of $H^s$ w.r.t. local frame of reference becomes\footnote{Here $f_y$ denotes partial derivative of $f$ w.r.t. $y$.}
\begin{align*}
\mathcal{V}_{\sigma_1,\sigma_2}(f)=\begin{cases}
e^{-\frac{1}{2\sigma_1^2}\big{(}\frac{f^s_{\eta\eta}}{f^s_{\xi\xi}}\big{)}}\bigg{(}
1-e^{-\frac{(f^s_{\eta\eta})^2+(f^s_{\xi\xi})^2}{2\sigma_2^2}}\bigg{)} & \text{ if } f_{\eta\eta}>0
\\
0 & \text{ if } f_{\eta\eta}\leq 0.
\end{cases}
\end{align*} 
Therefore on $SE(2)$-CS transformed image where we have left-invariant frames $\{\partial_\xi,\partial_\eta,\partial_\theta\}$ with $\xi=x\cos\theta+y\sin\theta$, $\eta=-x\sin\theta+y\cos\theta$ we propose the vesselness filter on $SE(2)$-CS transformed image $U:=\mathcal{W}_\psi(f)$ 
\begin{align*}
\mathbf{x} \mapsto \frac{|\mathcal{V}_{\sigma_{1},\sigma_{2}} (U (\cdot,\theta))|(\mathbf{x})}{ 
\max \limits_{(\mathbf{x},\theta) \in SE(2)} |\mathcal{V}_{\sigma_{1},\sigma_{2}} (U(\cdot,\theta))|(\mathbf{x})
}  \ \ \   \text{ with } 
\sigma_1=\frac{1}{2}, \  \sigma_2=0.2\max\limits_{x}\Delta_x U(x,\theta)=0.2\|S\|_\infty.
\end{align*}
However, akin to the coherence enhancing diffusion on orientation scores \cite{Duits2010a} better results are obtained by using a local gauge frame $\{\partial_a^s,\partial_b^s,\partial_c^s\}$ defined in Eq. \eqref{eq:GaugeFrames} with the best exponential curve fit (whose tangent $c$ corresponds to the eigenvalue of $3$D-Hessian $H_\beta^sU$ with smallest eigenvalue). In this case for all $g\in SE(2)$ we have
\begin{align}
\widetilde{\mathcal{V}}_{\sigma_1,\sigma_2}(U)=\begin{cases}
e^{-\frac{1}{4\sigma_1^2}\frac{(U^s_{cc})^2}{(U^s_{aa}+U^s_{bb})^2}}\bigg{(}
1-e^{-\frac{(U^s_{aa}+U^s_{bb})^2+(U^s_{cc})^2}{2\sigma_2^2}}\bigg{)} & \text{ if } U^s_{aa}+U^s_{bb}>0
\\
0 & \text{ if } U^s_{aa}+U^s_{bb}\leq 0.
\end{cases}\label{eq:Gauge-Vesselness}
\end{align} 
Note that we again have anisotropy measure $R=\frac{2\lambda_1}{\lambda_2+\lambda_3}$,  structureness $S=\lambda_1^2+\frac{(\lambda_2+\lambda_3)^2}{2}$, $\sigma_1=\frac{1}{2}$ and $\sigma_2=0.2 \| U_{aa}^{s}+U_{bb}^{s} \|_{\infty}$. Then we normalize 
\begin{align}\label{eq:SE2_Vess}
(V_{\sigma_1\sigma_2}(U))(g)=\frac{(\widetilde{\mathcal{V}}_{\sigma_1,\sigma_2}(U))(g)}
{\|\widetilde{\mathcal{V}}_{\sigma_1,\sigma_2}(U)\|_\infty}
\end{align}
to get the result of $SE(2)$ vesselness filter. Note that Gaussian derivatives in the score must be isotropic with scale $s>0$ in the spatial part (in order to preserve
the commutators of the generators during regularization)
and $s=\beta^2$ in the angular part, for details see \cite{Hannink2014},\cite[Ch.5]{Franken2008}. 
\\

Finally, when extending the framework to CW transformed image (multiple-scale OS) $(\mathcal{W}_\psi f)(x,y,\theta,a_l)$ with $a_l=\exp(ls_\rho)a_0$ we apply the $SE(2)$ vesselness filter on each scale layer $(\mathcal{W}_\psi f)(\cdot,\cdot,\cdot,a_l)$ with $a_l$ fixed where we adapt $\sigma_2(a_l)=0.1 \exp(ls_\rho)\|\mathcal{U}_f(\cdot,\cdot,\cdot,a_l)\|_\infty$.
\begin{figure*}[t]
        \centering
        \begin{subfigure}[t]{0.24\textwidth}
          \centering
           \includegraphics[width=\textwidth]{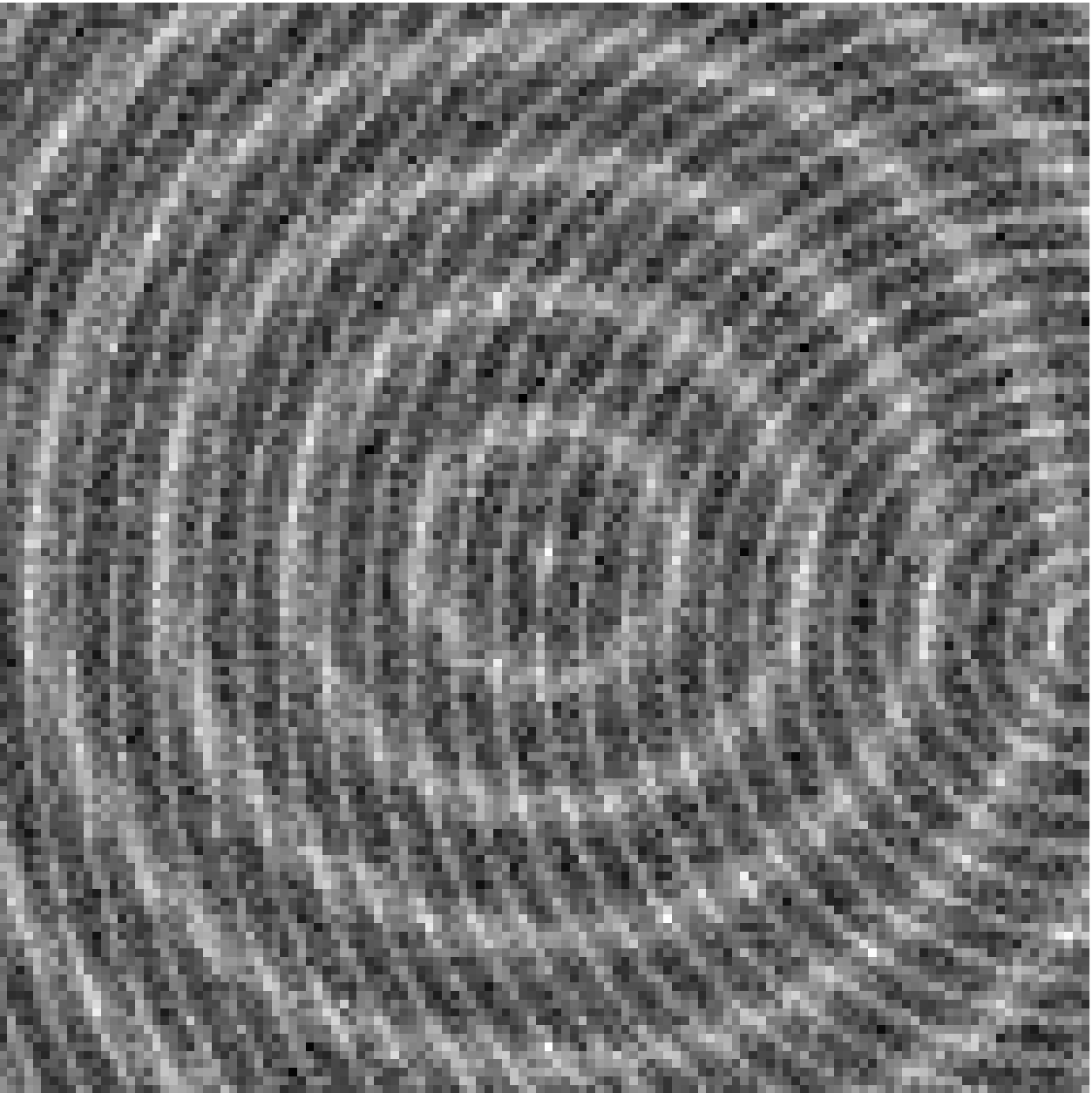}
            \caption{Noisy image}
        \end{subfigure}
        \begin{subfigure}[t]{0.24\textwidth}
                \centering
                \includegraphics[width=\textwidth]{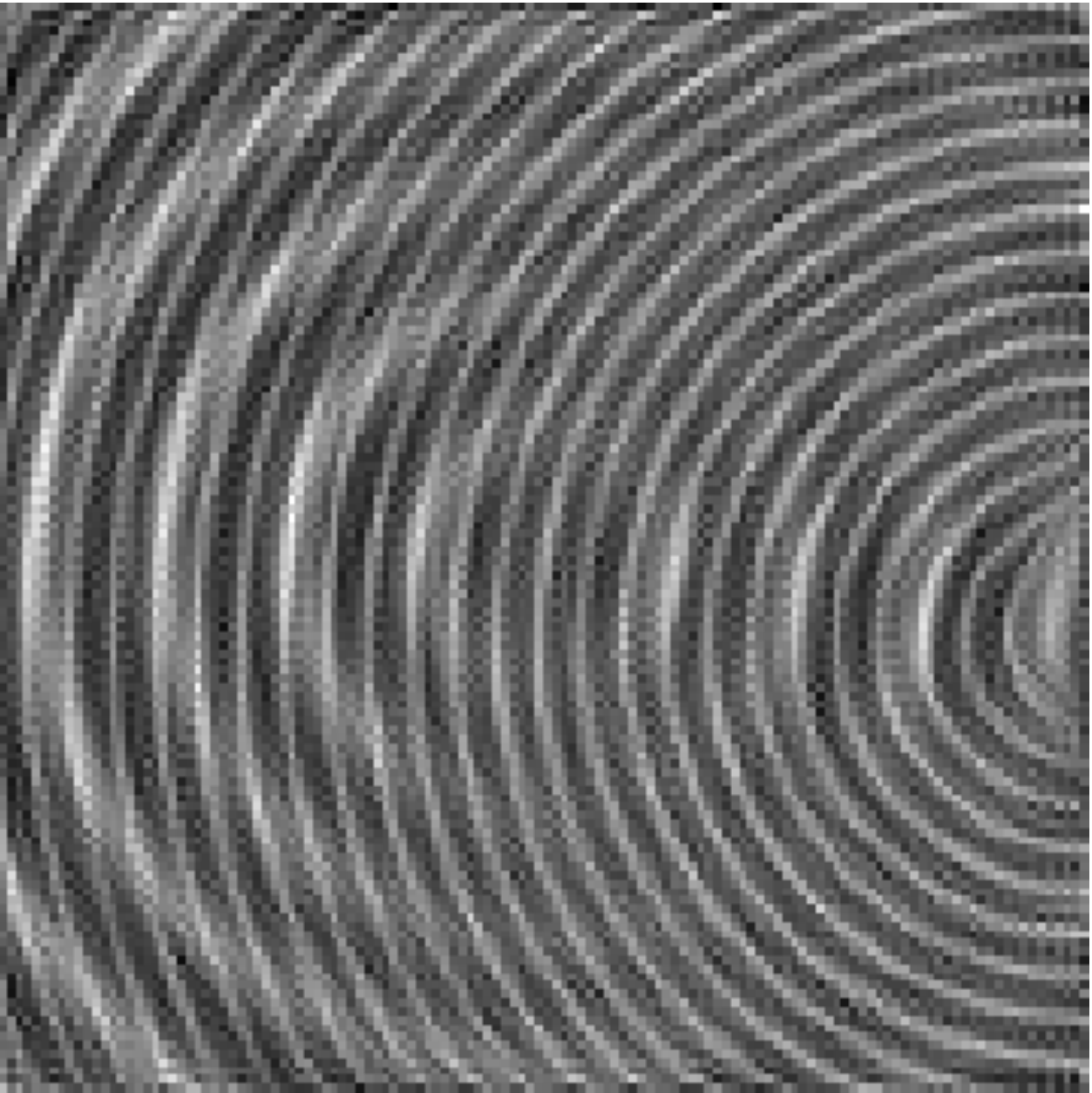}
                \caption{CED $t=15$}
        \end{subfigure}
\begin{subfigure}[t]{0.24\textwidth}
                \centering
                \includegraphics[width=\textwidth]{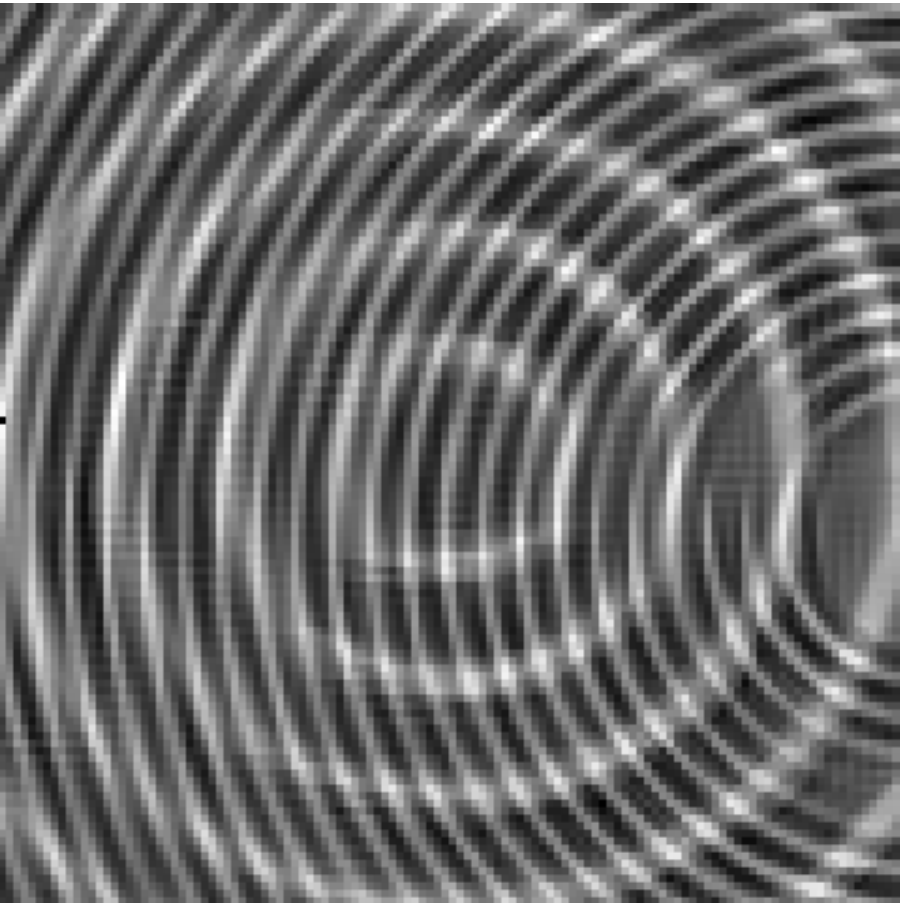}
                \caption{CED-OS $t=6$}
        \end{subfigure}
        \begin{subfigure}[t]{0.24\textwidth}
                \centering
                \includegraphics[width=\textwidth]{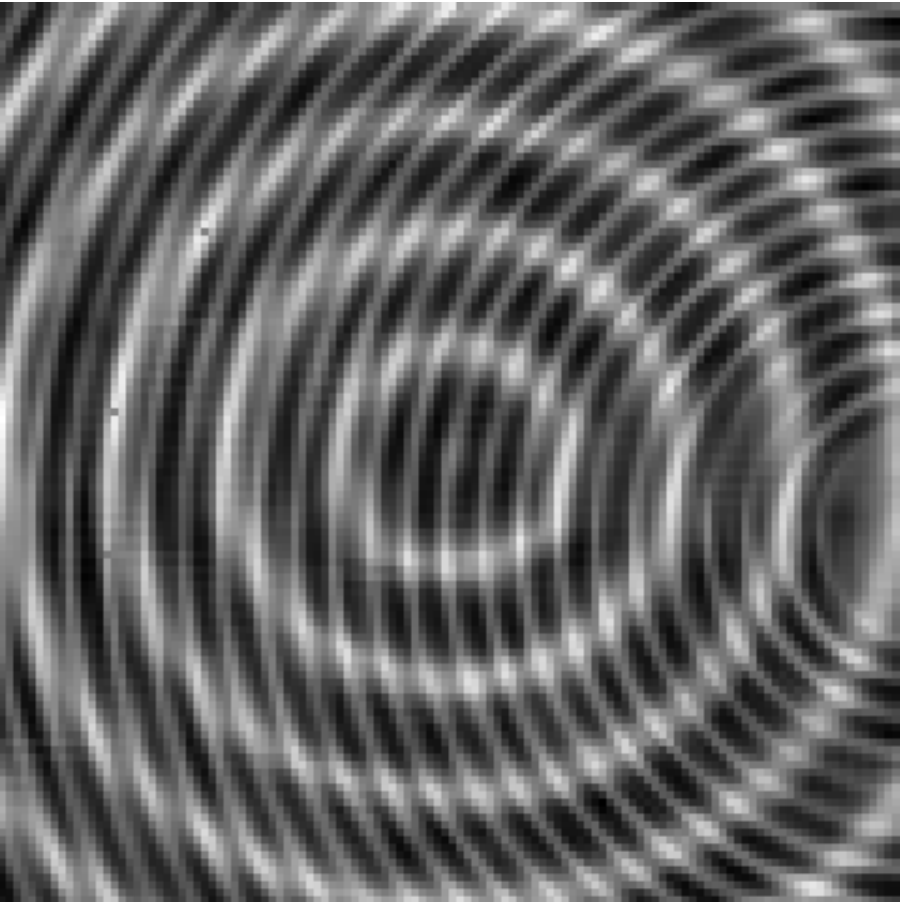}
                \caption{CED-SOS, {\footnotesize $t_i=0,2,6,12$} }
        \end{subfigure}
        \caption
{Comparison of CED \cite{Weickert1999}, CED-OS \cite{Franken2009} and our proposed method CED-SOS. In CED-SOS multi-scale crossing stuctures are better preserved. Observe that CED-OS and CED-SOS deal  appropriately with crossings, unlike CED} 
\label{fig:NonLinearDiffusionArtIm}
\end{figure*}
\section{Results}
Parameters used for creating WT on images (scale orientation score): $\text{No. of scales}=4$, $\text{No. of orientations}=20$. The parameters that we used for CED are (for explanation of parameters see\cite{Weickert1999}): $\sigma=0.5, \ \rho =4, \  C =1 \text{ and } \alpha=0.001$. The non-linear diffusion parameters for CED-OS and CED-SOS for each scale are: $\rho_s=12, \ \widetilde{\rho}= 1.5, \ \beta = 0.058 \text{ and } c = 0.08$. All the experiments in this section use these parameters (unless otherwise indicated in the caption) and the varying end times are indicated below the images. We have enforced horizontality in the experiments involving CED-OS, see \cite{Franken2009} for more details.\\

Figure~\ref{fig:NonLinearDiffusionArtIm} compares the effect of CED, CED-OS and CED-SOS on an artificial image containing additive superimposition of two images with concentric circles of varying widths. The new proposed method yields visually better results compared to both CED and CED-OS.\\

\begin{figure*}[t]
        \centering
        \begin{subfigure}{0.24\textwidth}
          \centering
           \includegraphics[width=\textwidth]{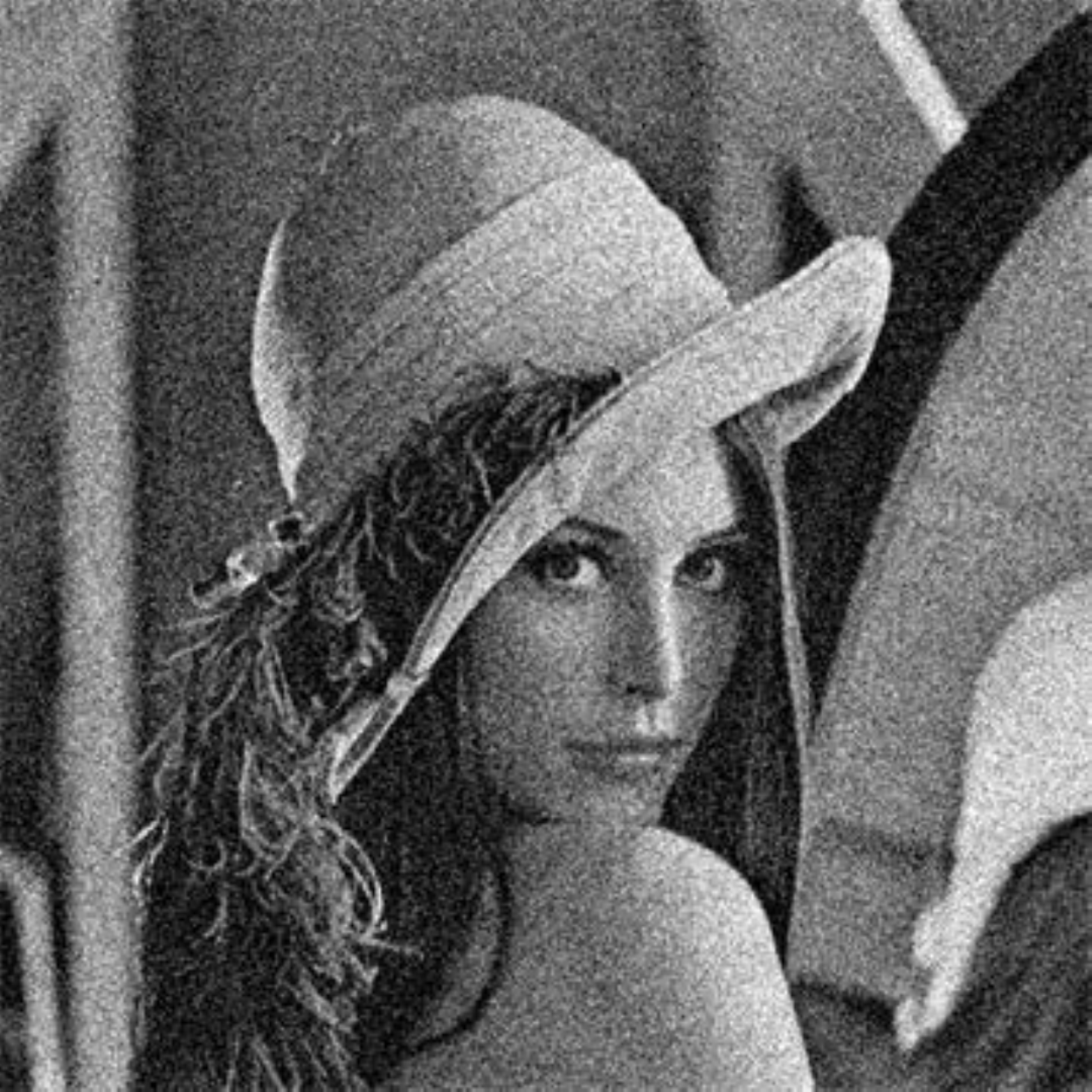}
        \end{subfigure}
        \begin{subfigure}{0.24\textwidth}
                \centering
                \includegraphics[width=\textwidth]{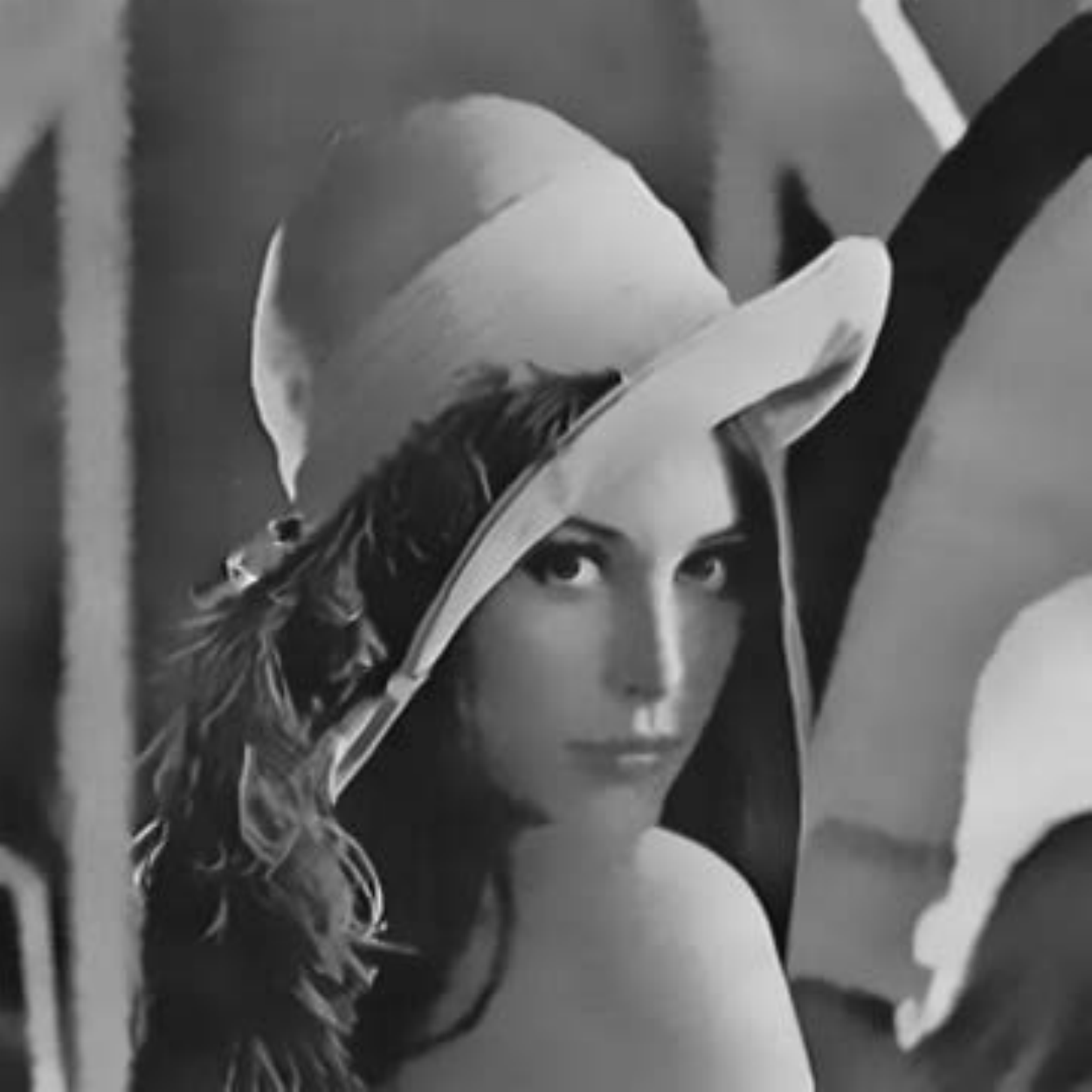}
        \end{subfigure}
\begin{subfigure}{0.24\textwidth}
                \centering
                \includegraphics[width=\textwidth]{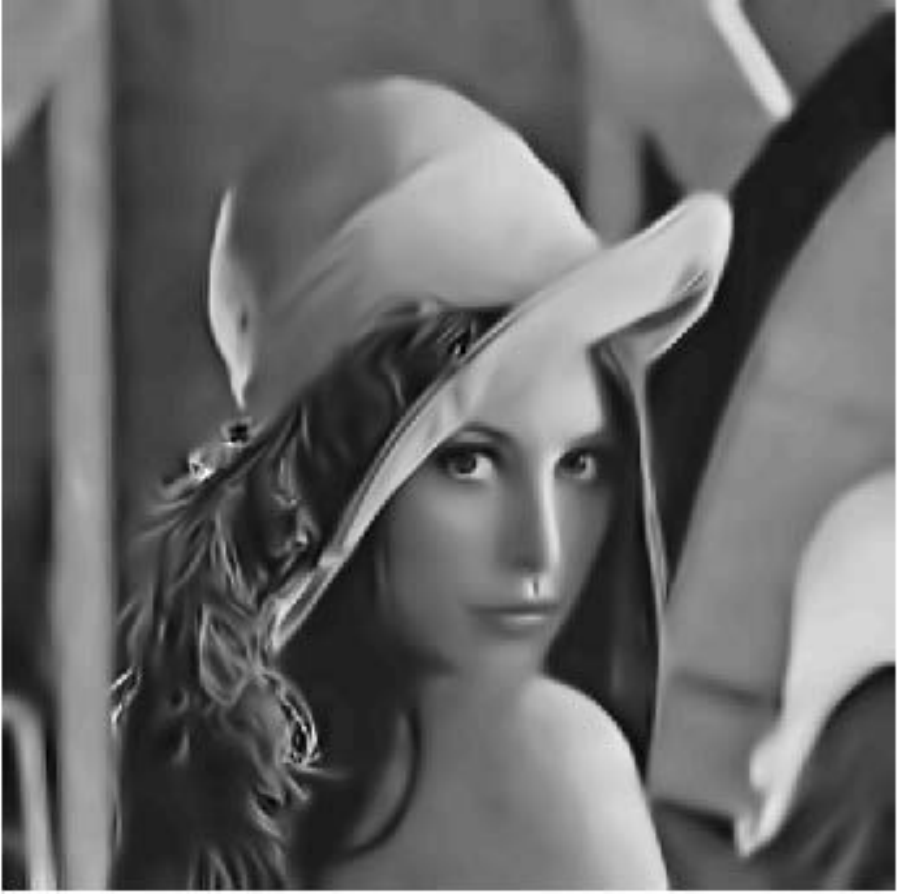}
        \end{subfigure}
        \begin{subfigure}{0.24\textwidth}
                \centering
                \includegraphics[width=\textwidth]{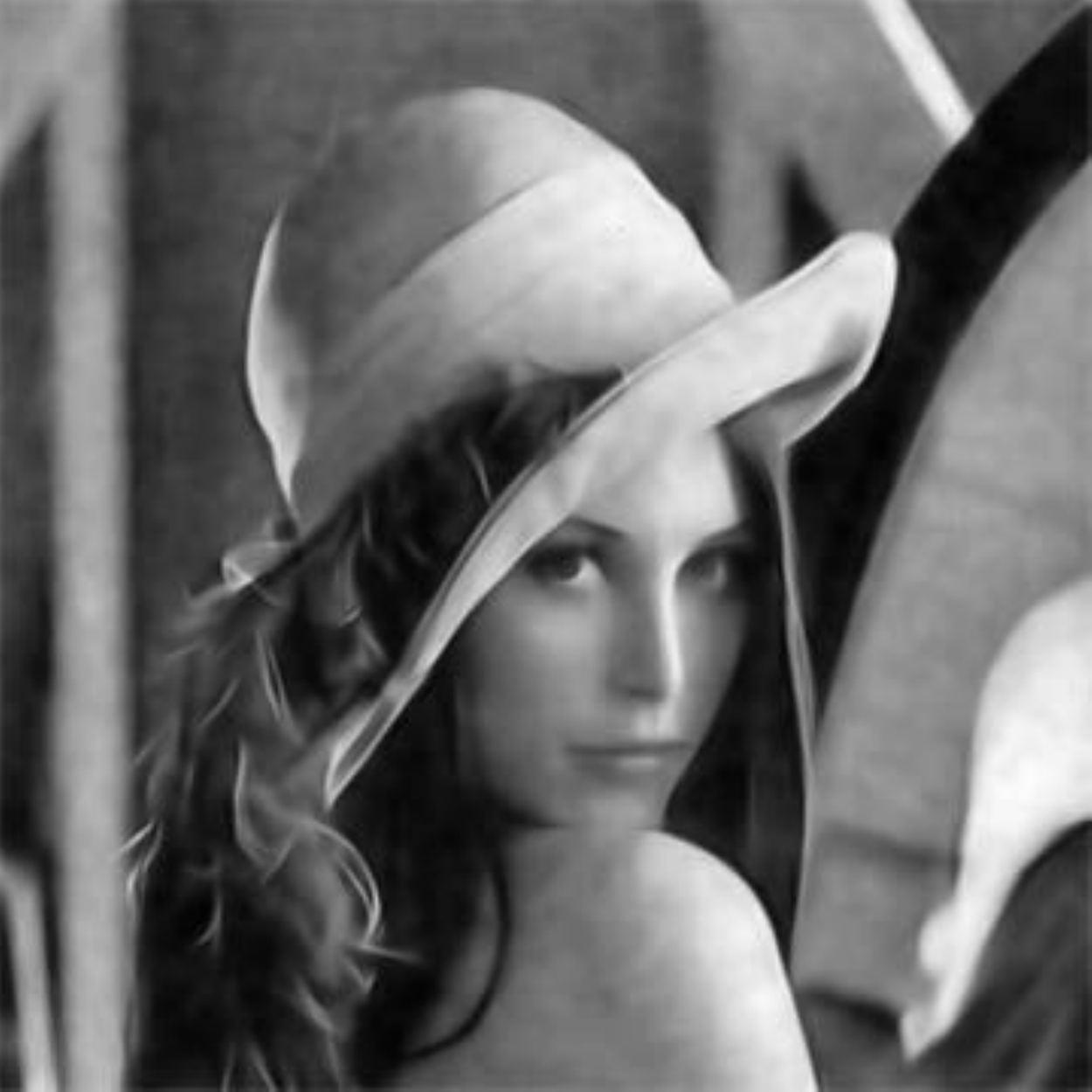}
        \end{subfigure}

        \begin{subfigure}{0.24\textwidth}
          \centering
           \includegraphics[width=\textwidth]{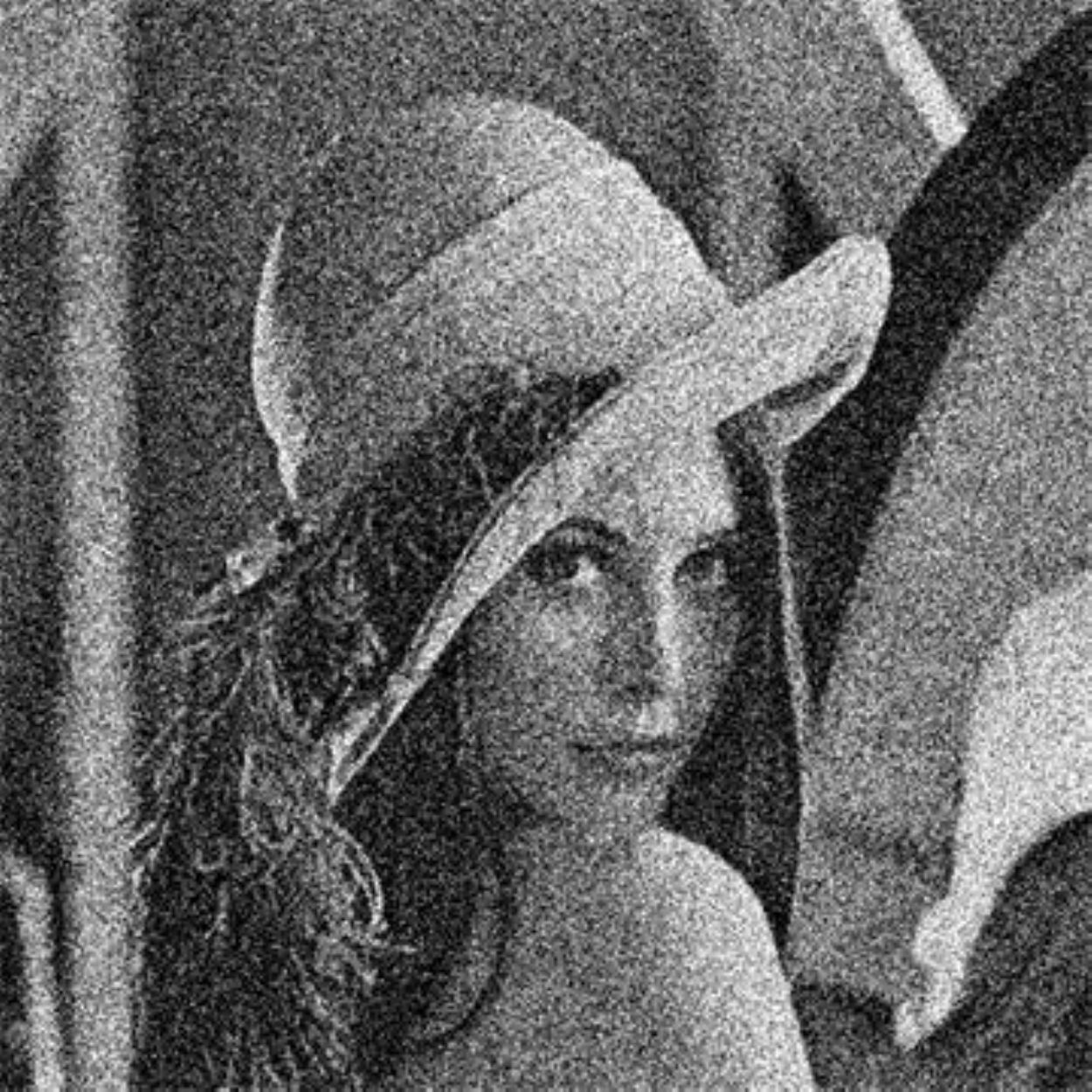}
            \caption{Noisy image}
        \end{subfigure}
        \begin{subfigure}{0.24\textwidth}
                \centering
                \includegraphics[width=\textwidth]{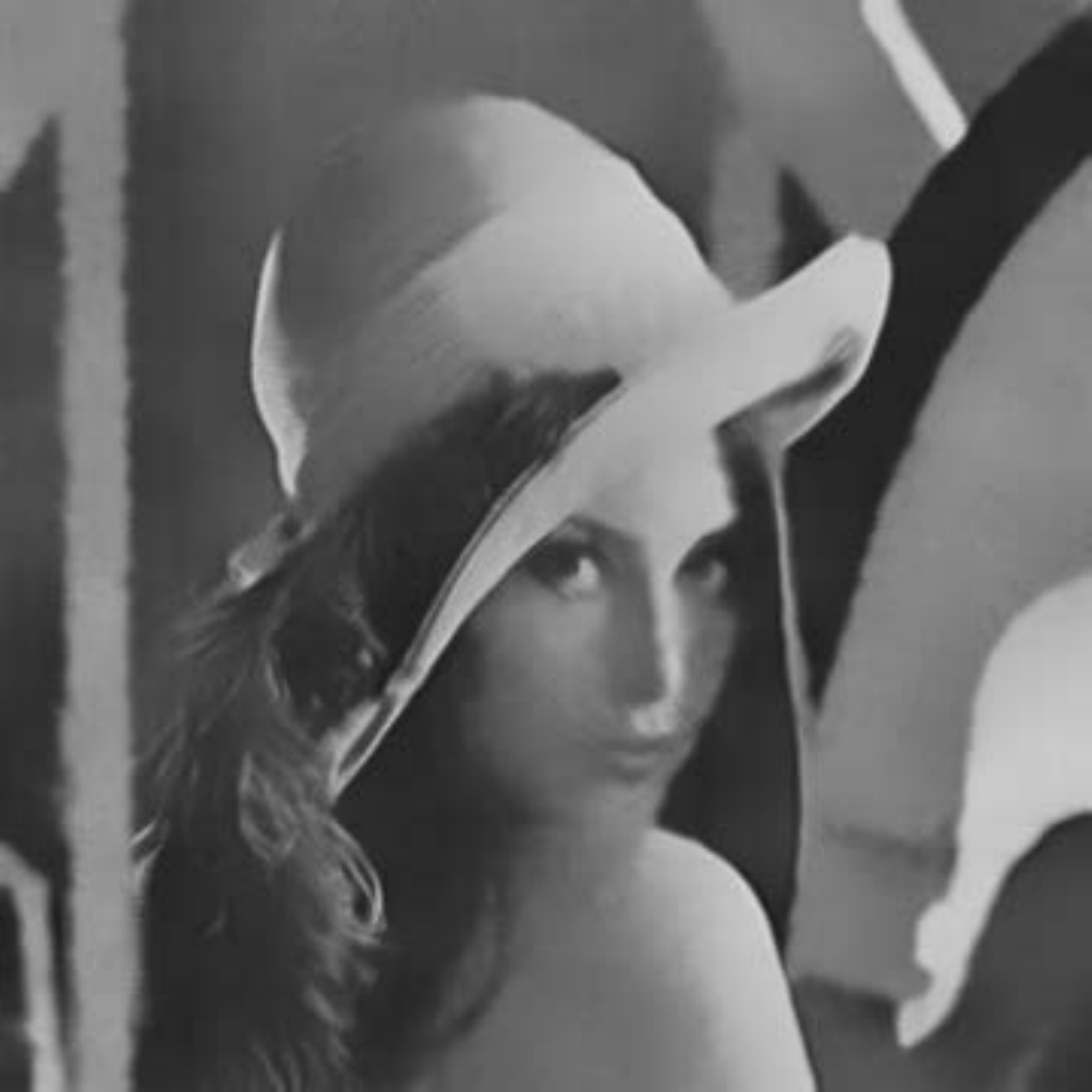}
                \caption{NLM}
        \end{subfigure}
\begin{subfigure}{0.24\textwidth}
                \centering
                \includegraphics[width=\textwidth]{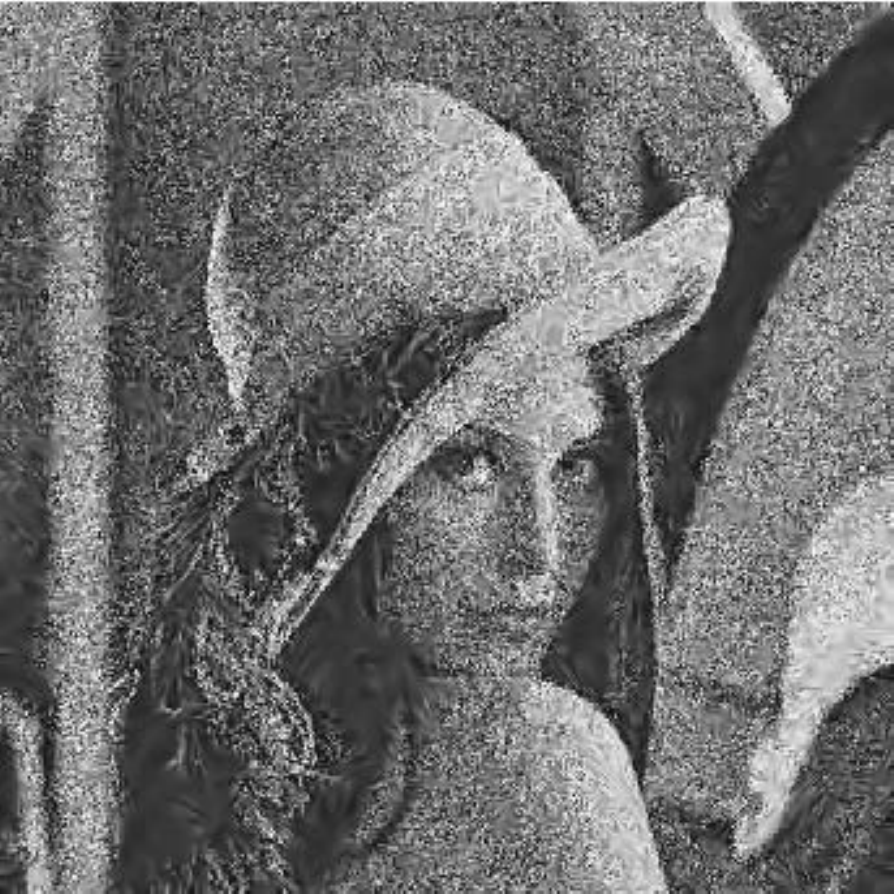}
                \caption{ISKR: Iteration 35}
        \end{subfigure}
        \begin{subfigure}{0.24\textwidth}
                \centering
                \includegraphics[width=\textwidth]{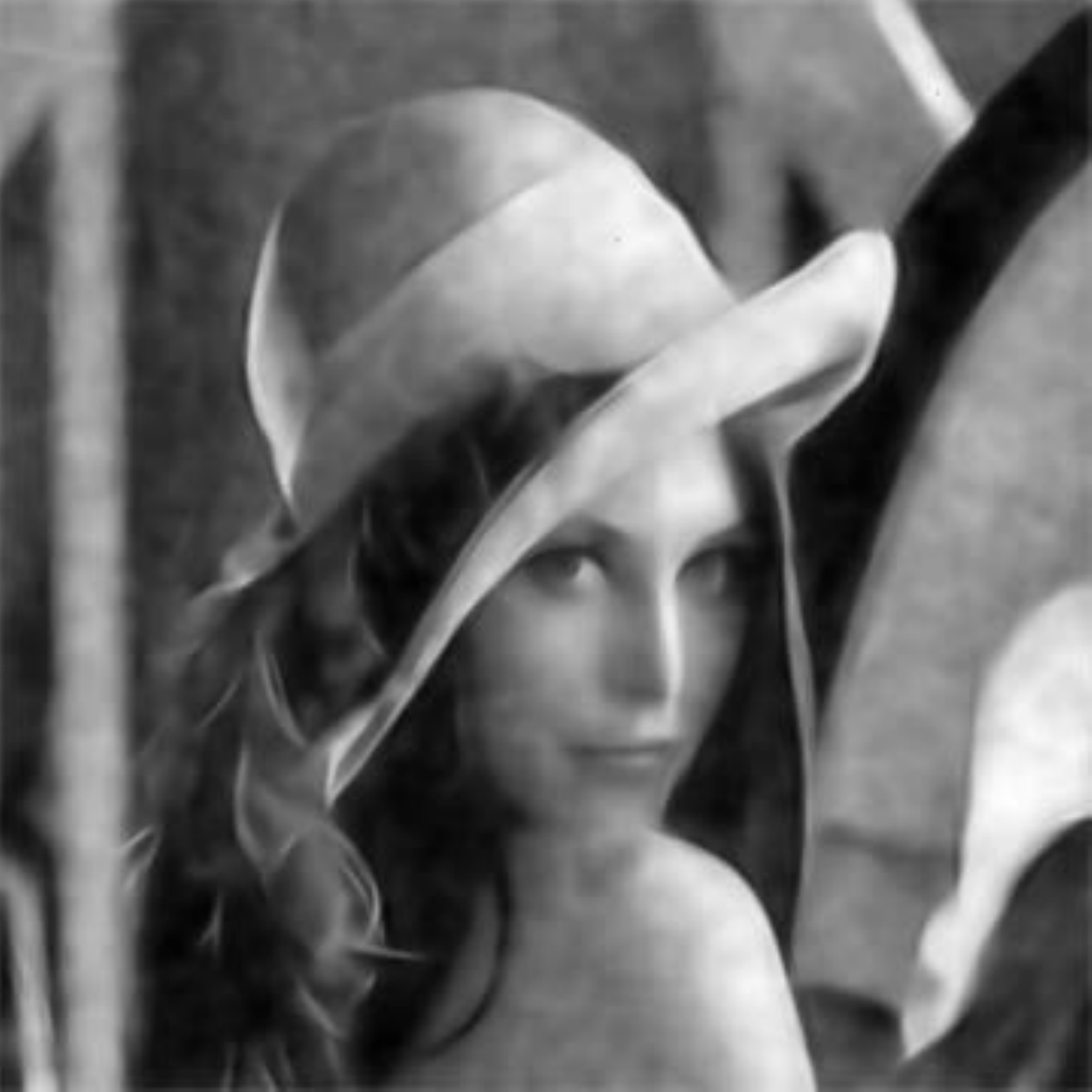}
                \caption{CED-SOS}
        \end{subfigure}
        \caption
{Comparision of CED-SOS with the state of the art denoising algorithms \cite{Seo2007} with standard parameters on Lena image with gaussian noise (top: $\sigma=35$, bottom: $\sigma=60$). (a) Noisy input image (b) Non-local Means filter \cite{Buades2005}  (c) Iterative steering kernel regression \cite{Takeda2007} (d) CED-SOS, top: $t_i=0,2,8,9$, bottom: $t_i=0,2,12,14$ .} 
\label{fig:FullLena_LN_Comp}
\end{figure*}

In Figure~\ref{fig:FullLena_LN_Comp} 
we compare the proposed algorithm with state of the art denoising algorithms \cite{Seo2007} on a noisy Lena image. We present comparisions  with Non-local Means filter \cite{Buades2005} and Iterative steering kernel regression \cite{Takeda2007} using standard parameter settings.
Although our method is designed to deal with very noisy medical images containing crossing lines and various scales 
which differ from the Lena image, our method nonetheless performs very well compared to the state-of the art. In fact, our method shows better results in terms of robustness and higher noise levels, whereas the other methods show better results at lower noise levels. At the smaller noise levels further visual improvements of our method could be obtained by including median filtering techniques rather than diffusion techniques at large flat (isotropic) areas.\\

\begin{figure*}[t]
        \centering
        \begin{subfigure}{0.24\textwidth}
          \centering
           \includegraphics[width=\textwidth]{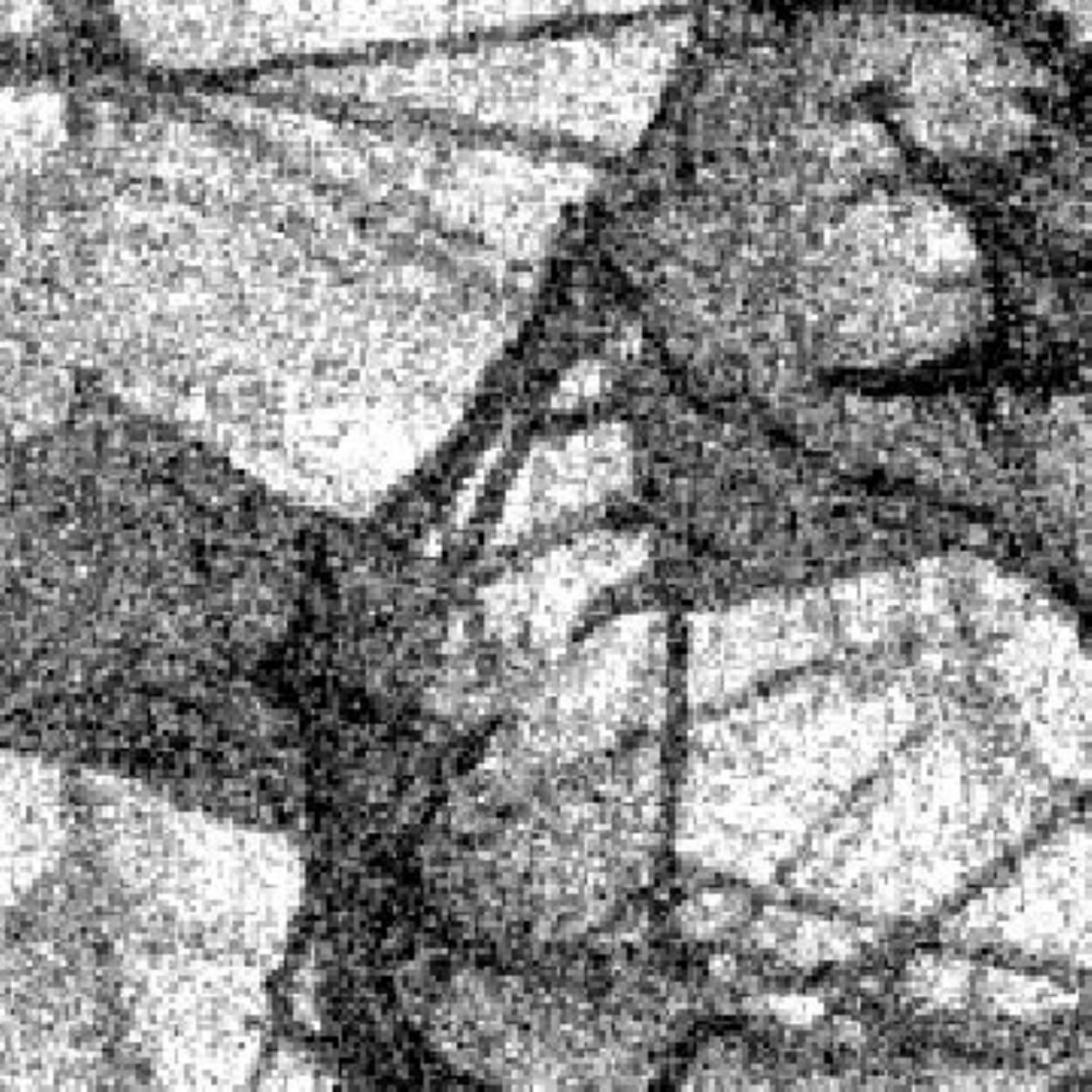}
            \caption{Noisy image}
        \end{subfigure}
        \begin{subfigure}{0.24\textwidth}
                \centering
                \includegraphics[width=\textwidth]{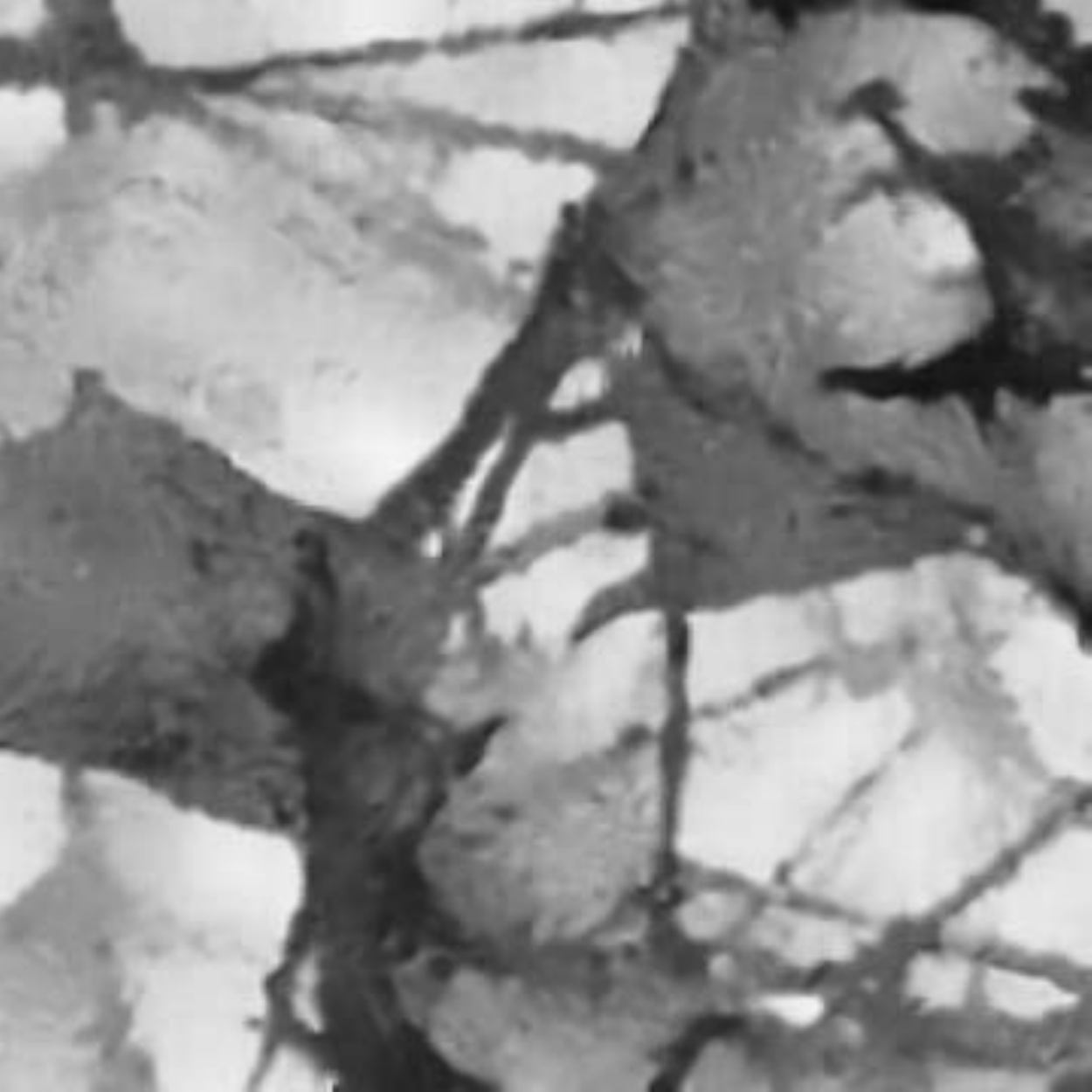}
                \caption{NLM}
        \end{subfigure}
		\begin{subfigure}{0.24\textwidth}
                \centering
                \includegraphics[width=\textwidth]{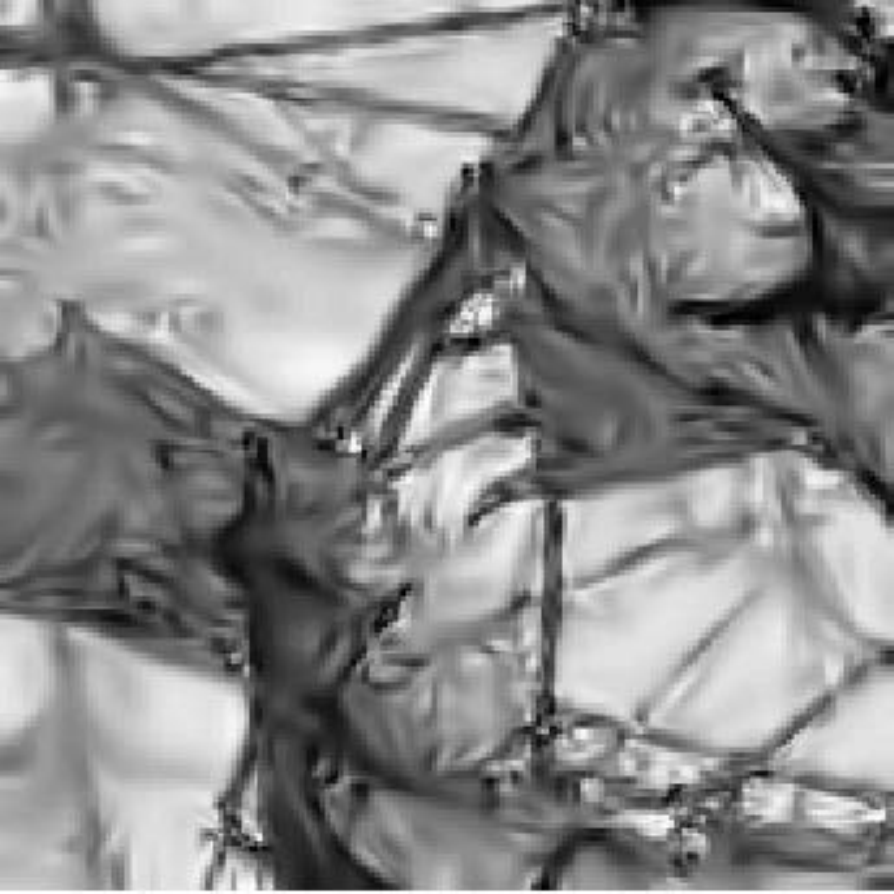}
                \caption{ISKR: Iter.35}
        \end{subfigure}
     \begin{subfigure}{0.24\textwidth}
                \centering
                \includegraphics[width=\textwidth]{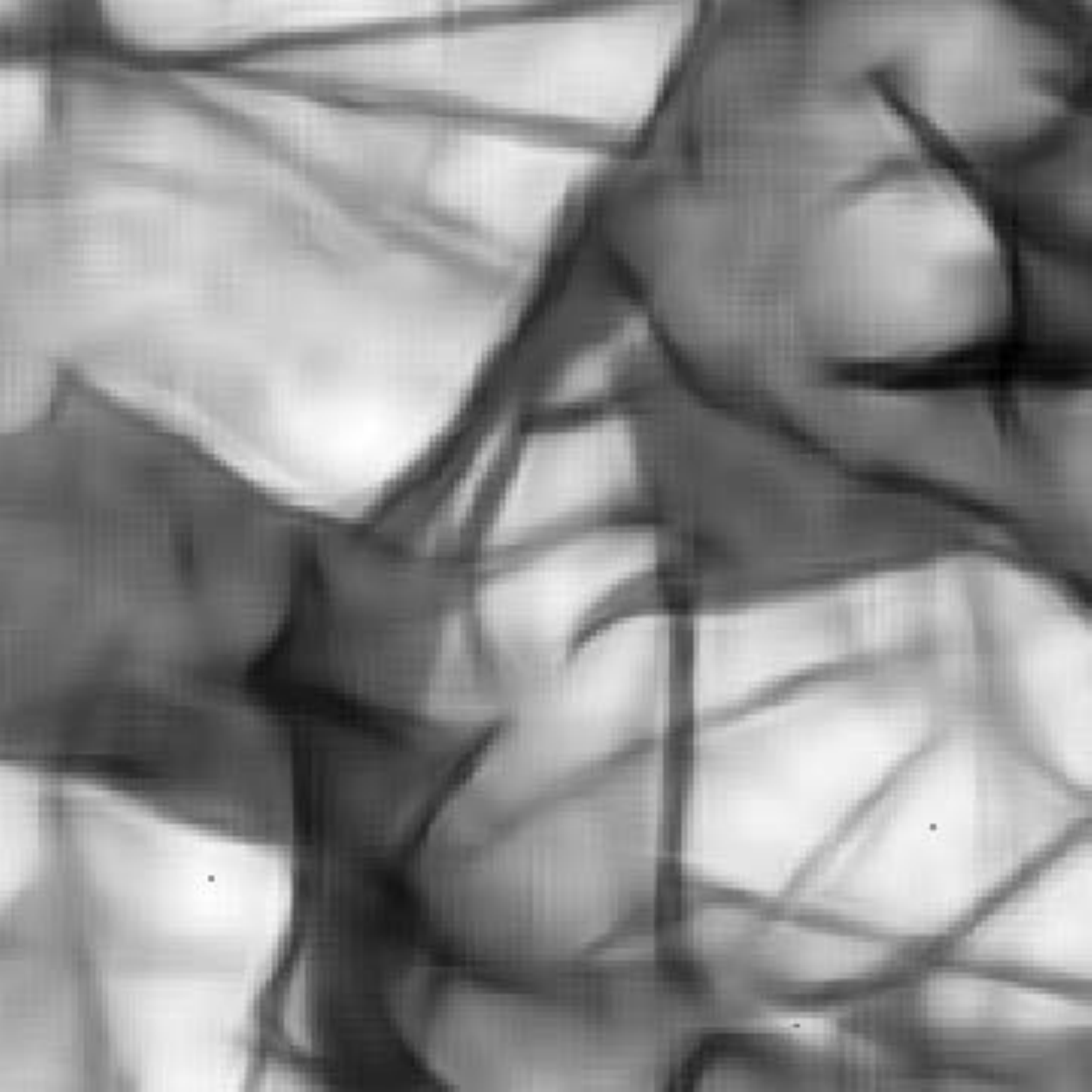}
                \caption{CED-SOS}
        \end{subfigure}
        \caption
{Comparision of CED-SOS with the state of the art denoising algorithms with standard parameters on a Collagen tissue from a cardiac wall. (a) Noisy input image (b) Non-local Means filter  (c) Iterative steering kernel regression, Iteration-35  (d) CED-SOS} 
\label{fig:Collagen_SOA_Comp}
\end{figure*}

\begin{figure*}[t]
        \centering
        \begin{subfigure}[t]{0.24\textwidth}
          \centering
           \includegraphics[width=\textwidth]{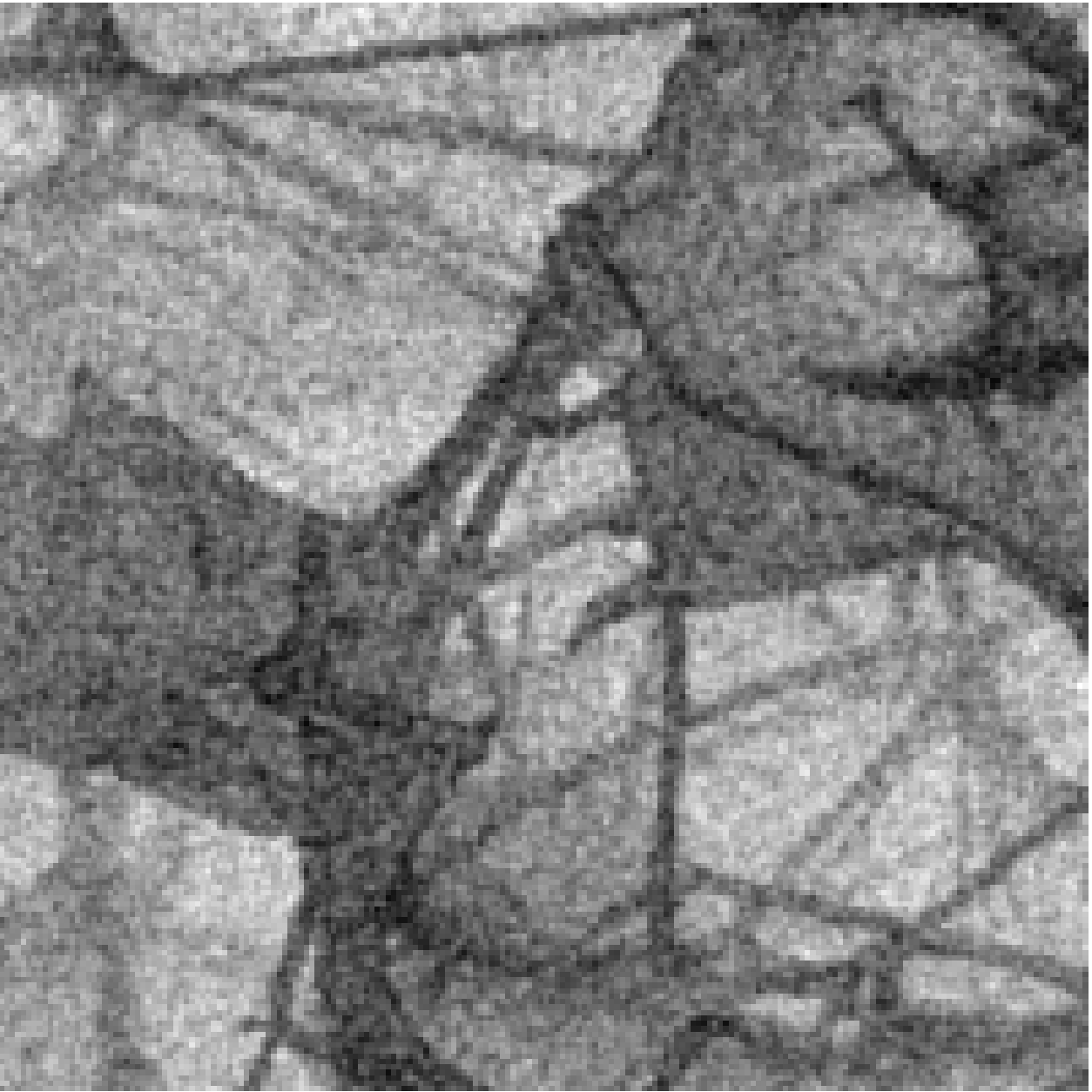}
            \caption{Noisy image}
        \end{subfigure}
        \begin{subfigure}[t]{0.223\textwidth}
                \centering
                \includegraphics[width=\textwidth]{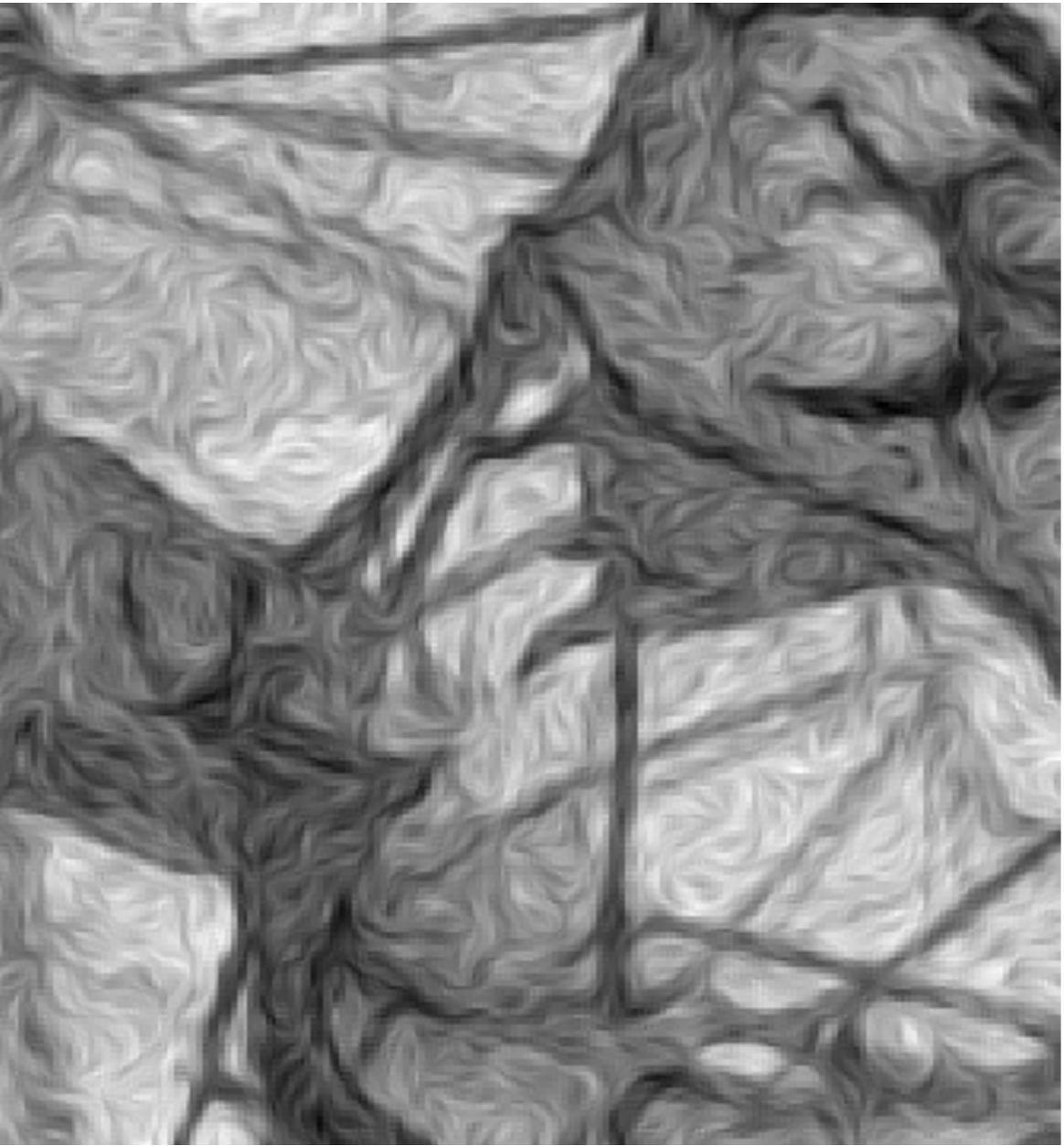}
                \caption{CED $t=4$}
        \end{subfigure}
\begin{subfigure}[t]{0.24\textwidth}
                \centering
                \includegraphics[width=\textwidth]{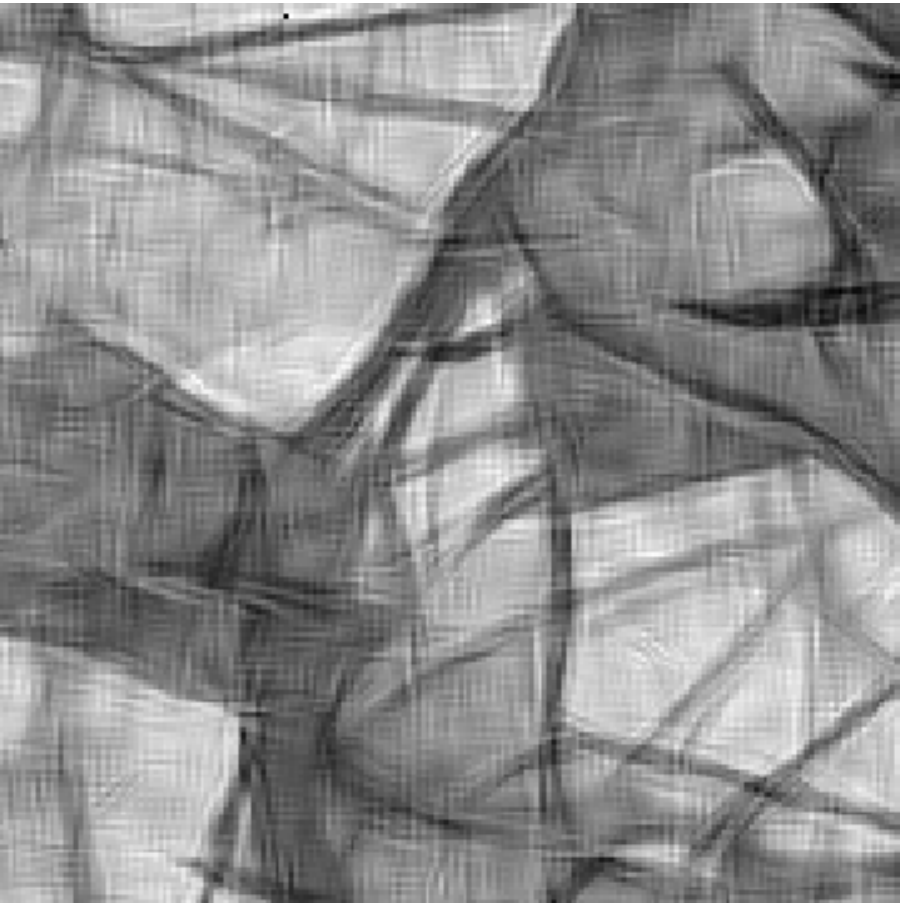}
                \caption{CED-OS $t=4$}
        \end{subfigure}
        \begin{subfigure}[t]{0.24\textwidth}
                \centering
                \includegraphics[width=\textwidth]{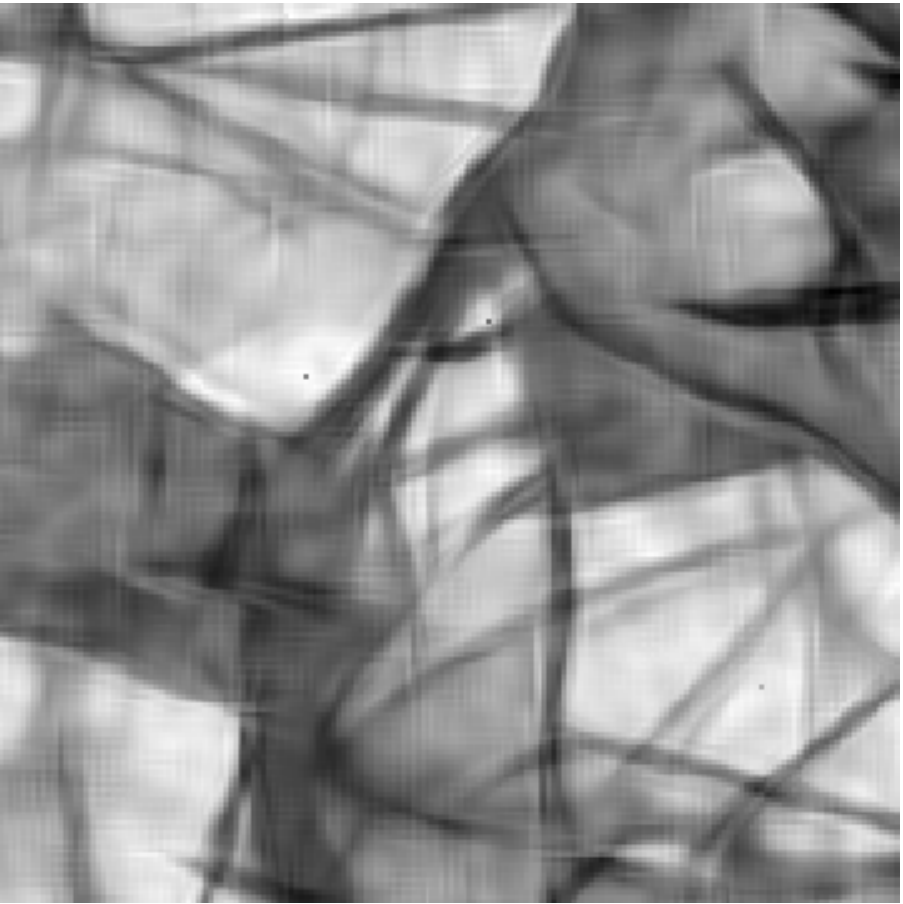}
                \caption{CED-SOS, {\footnotesize $t_i=0,1,3,16$}}
        \end{subfigure}
        
         \begin{subfigure}{0.24\textwidth}
                \centering
                \includegraphics[width=\textwidth]{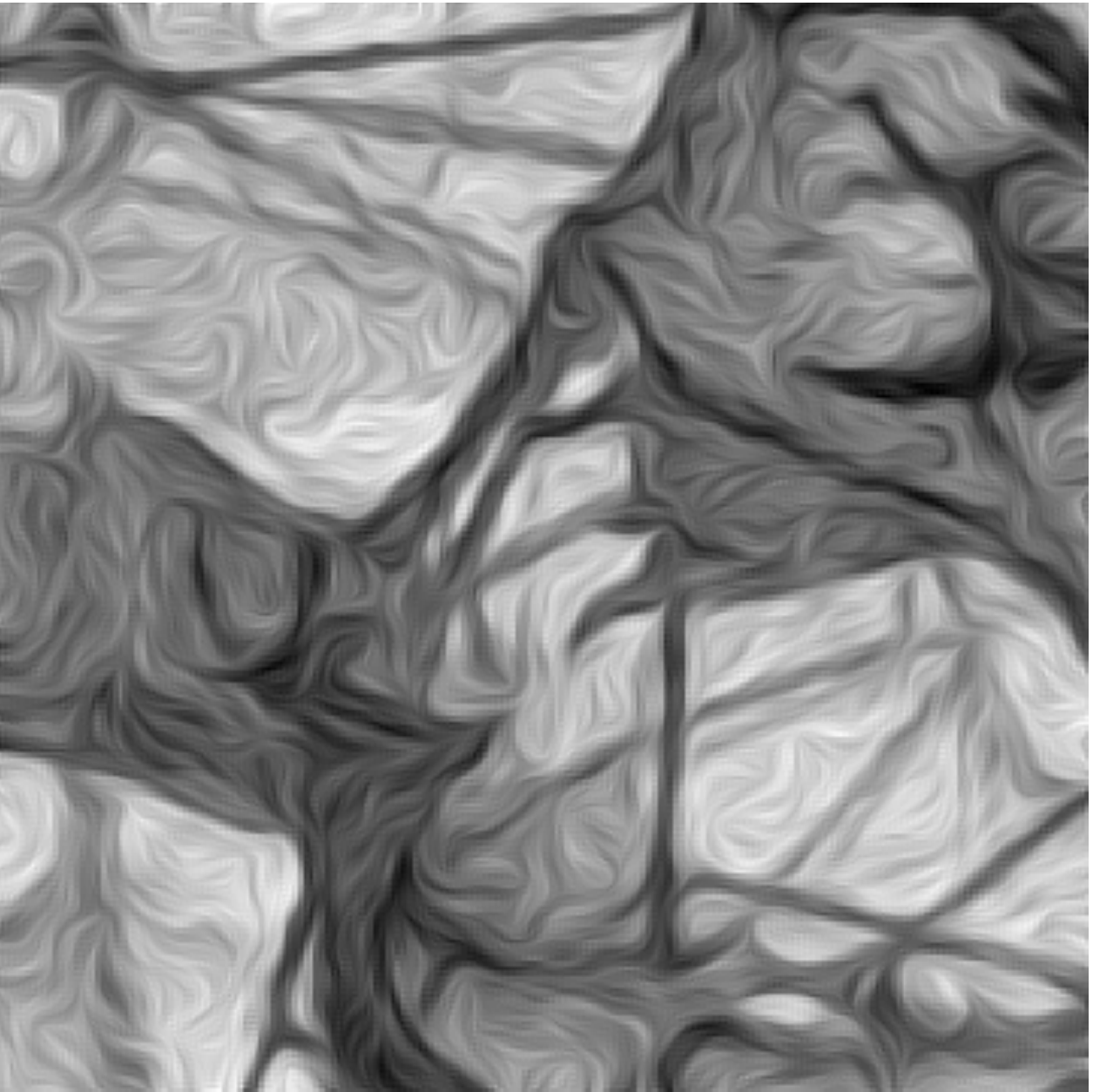}
                \caption{CED $t=10$}
        \end{subfigure}
        \begin{subfigure}{0.24\textwidth}
                \centering
                \includegraphics[width=\textwidth]{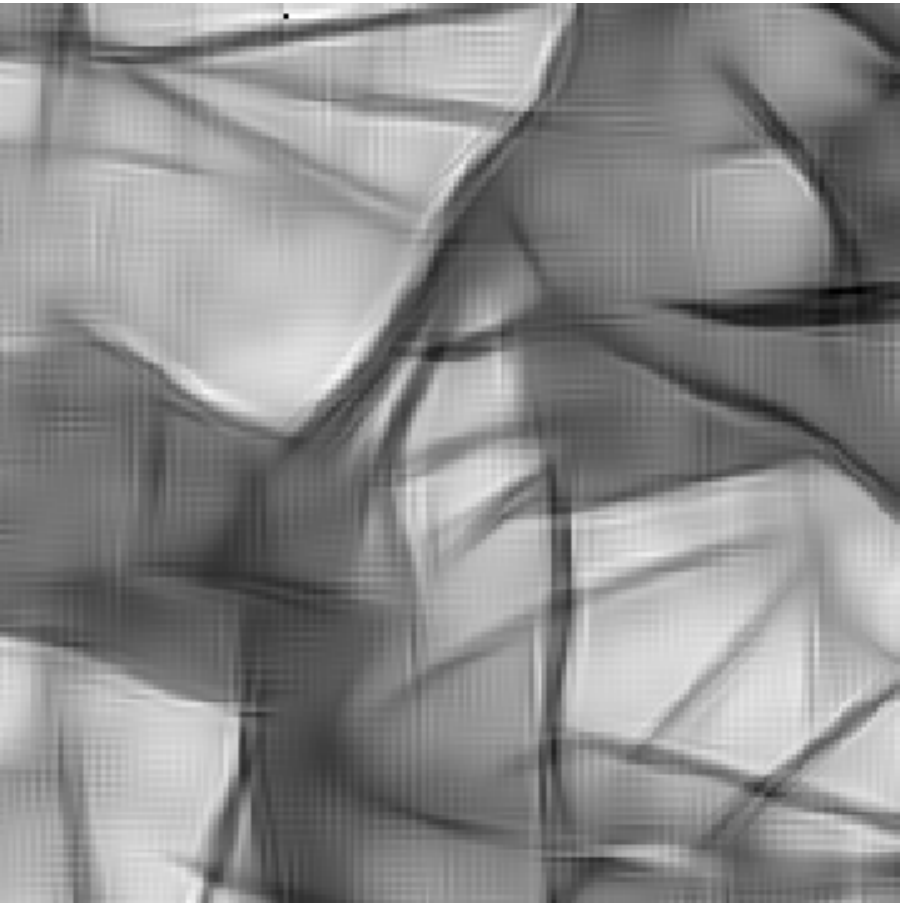}
                \caption{CED-OS $t=10$}
        \end{subfigure}
        \caption
{Results of CED, CED-OS and CED-SOS on a  
microscopy image of a Collagen tissue in the cardiac wall. CED-OS exhibits artefacts for small end time, see (c), which vanish if the algorithm is allowed to run longer at the cost of small-scale information, see (f). Therefore best visual results are obtained with CED-SOS.
} 
\label{fig:NonLinearDiffusionCollagen}
\end{figure*}

In Figure~\ref{fig:Collagen_SOA_Comp} we consider a typical example of a medical imaging application containing crossing lines at various scales in a low-contrast
medical (2-photon microscopy) image for which we actually designed our algorithm.  Note that these kind of images are acquired in tissue engineering research where the goal is to e.g. create artificial heart valves \cite{Rubbens2009}.  On these challenging type of images we obtain considerable improvements over the state of the art denoising algorithms as we appropriately deal with multiple scale crossings. Furthermore,
to illustrate the potential of the currently proposed algorithm over previous works on crossing preserving diffusion and coherence enhancing diffusion we have included Figure ~\ref{fig:NonLinearDiffusionCollagen}.\\
%

\begin{figure*}[h]
        \centering
        \begin{subfigure}{0.24\textwidth}
          \centering
           \includegraphics[width=\textwidth]{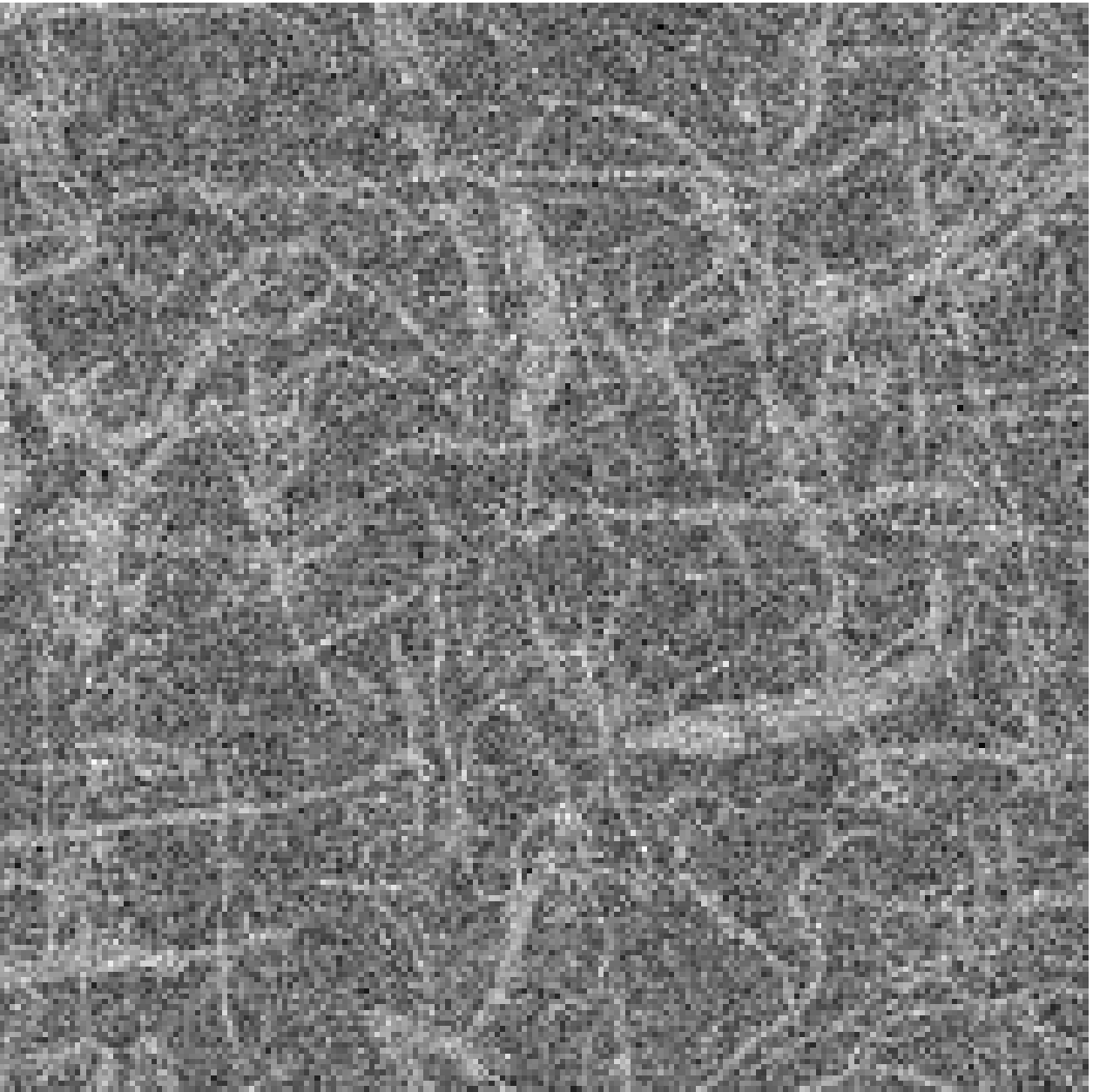}
            \caption{Noisy image}
        \end{subfigure}
        \begin{subfigure}{0.24\textwidth}
                \centering
                \includegraphics[width=\textwidth]{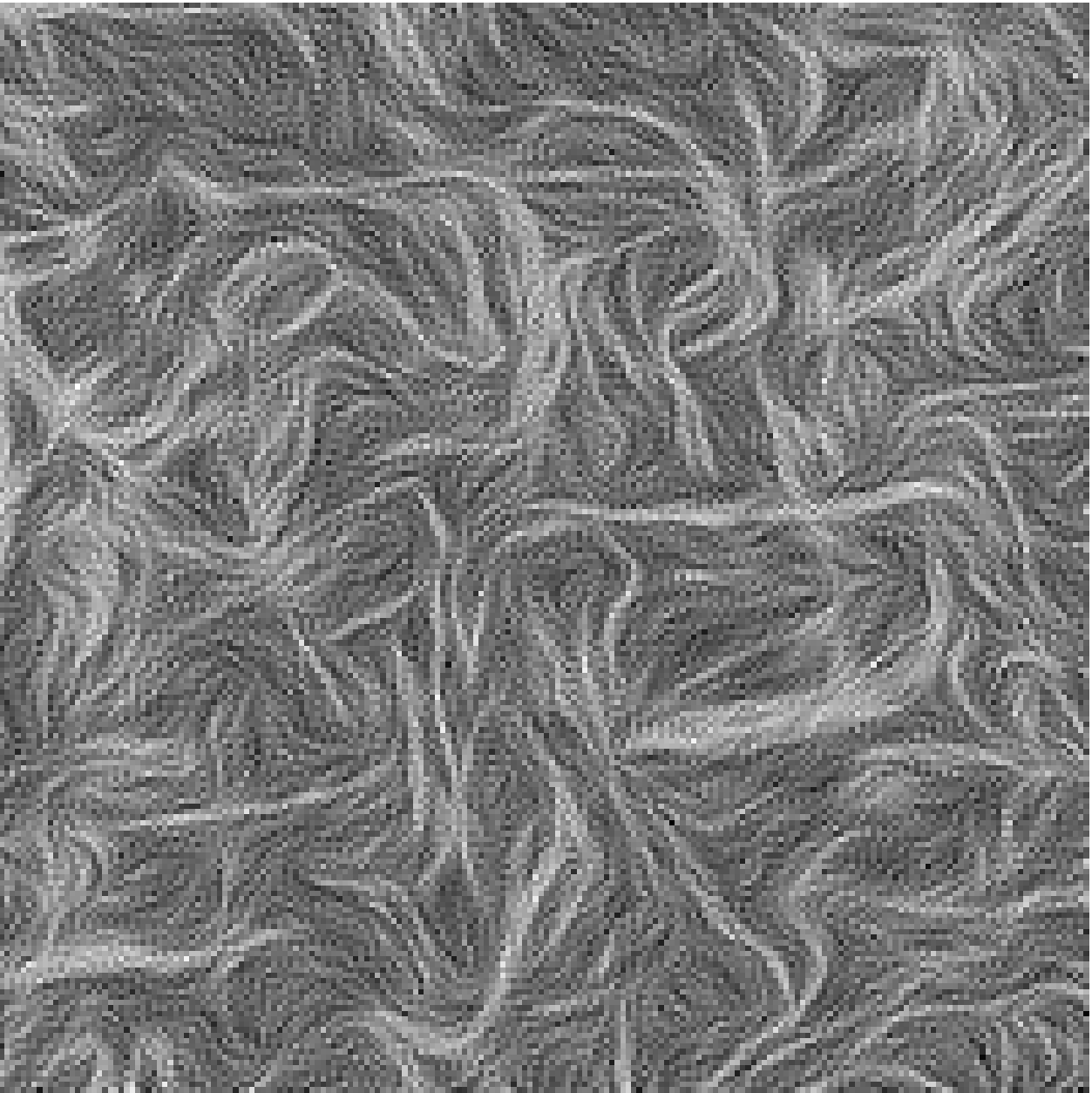}
                \caption{CED $t=6$}
        \end{subfigure}
\begin{subfigure}{0.24\textwidth}
                \centering
                \includegraphics[width=\textwidth]{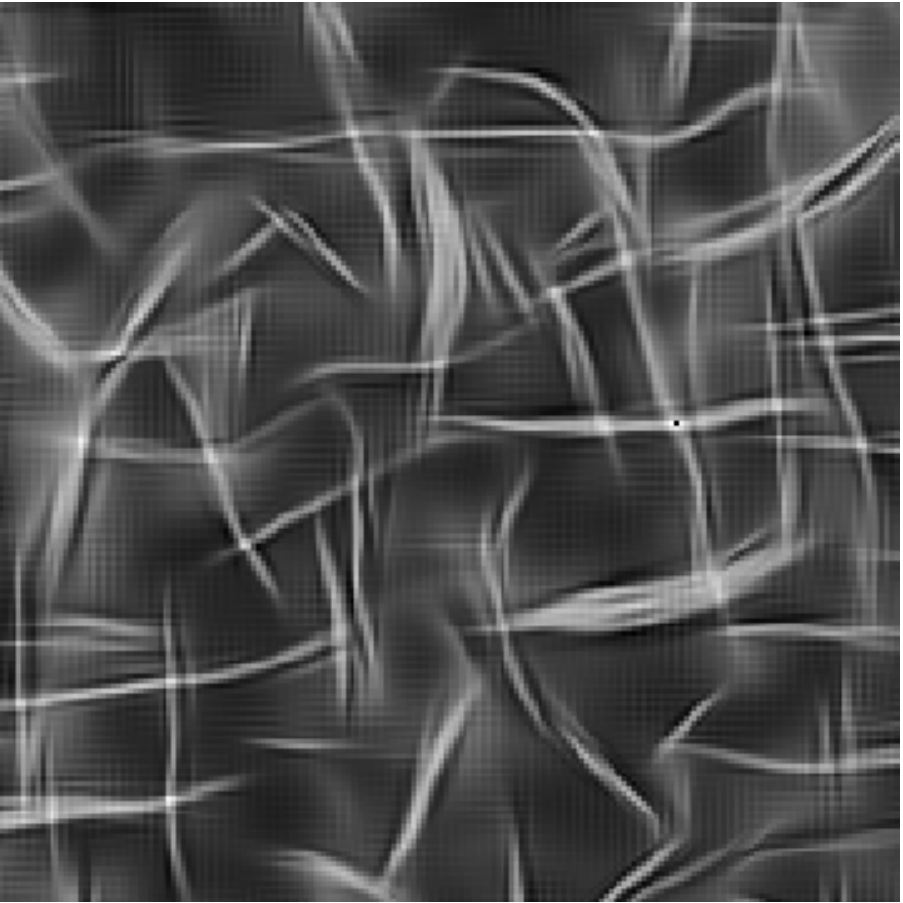}
                \caption{CED-OS $t=7$}
        \end{subfigure}
       \begin{subfigure}{0.24\textwidth}
                \centering
                \includegraphics[width=\textwidth]{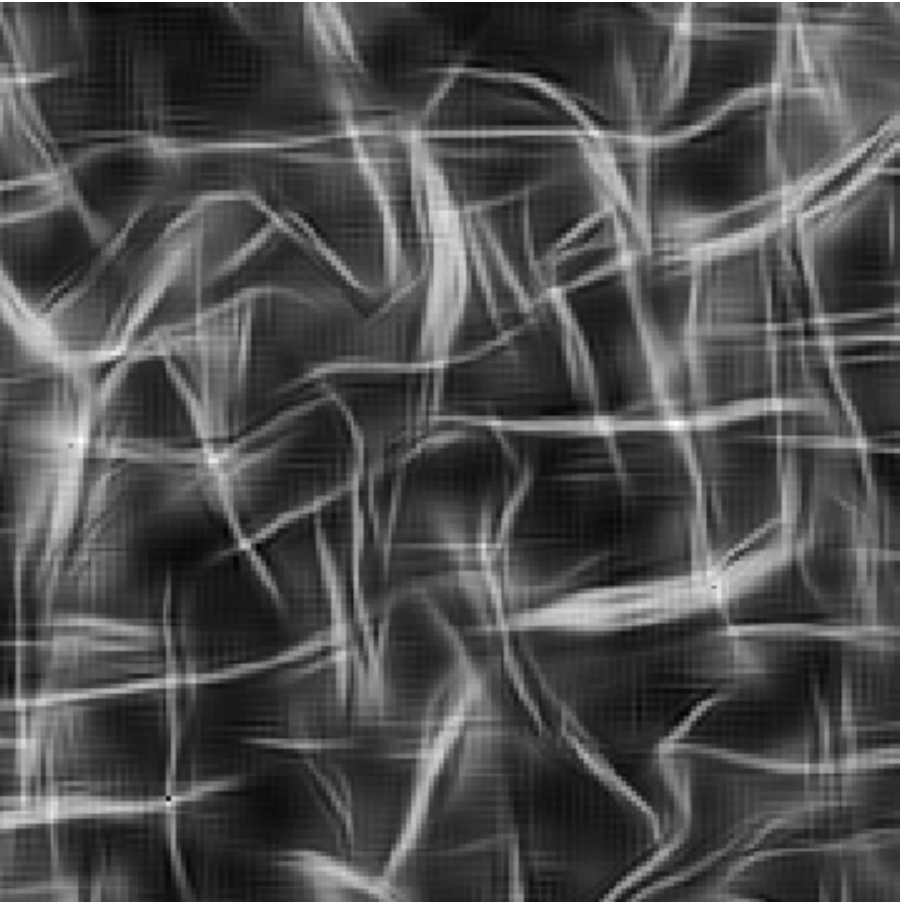}
                \caption{CED-SOS, {\footnotesize $t_i=0,1,6,10$} }
        \end{subfigure}
                \begin{subfigure}[t]{0.35\textwidth}
                \centering
                \includegraphics[width=\textwidth]{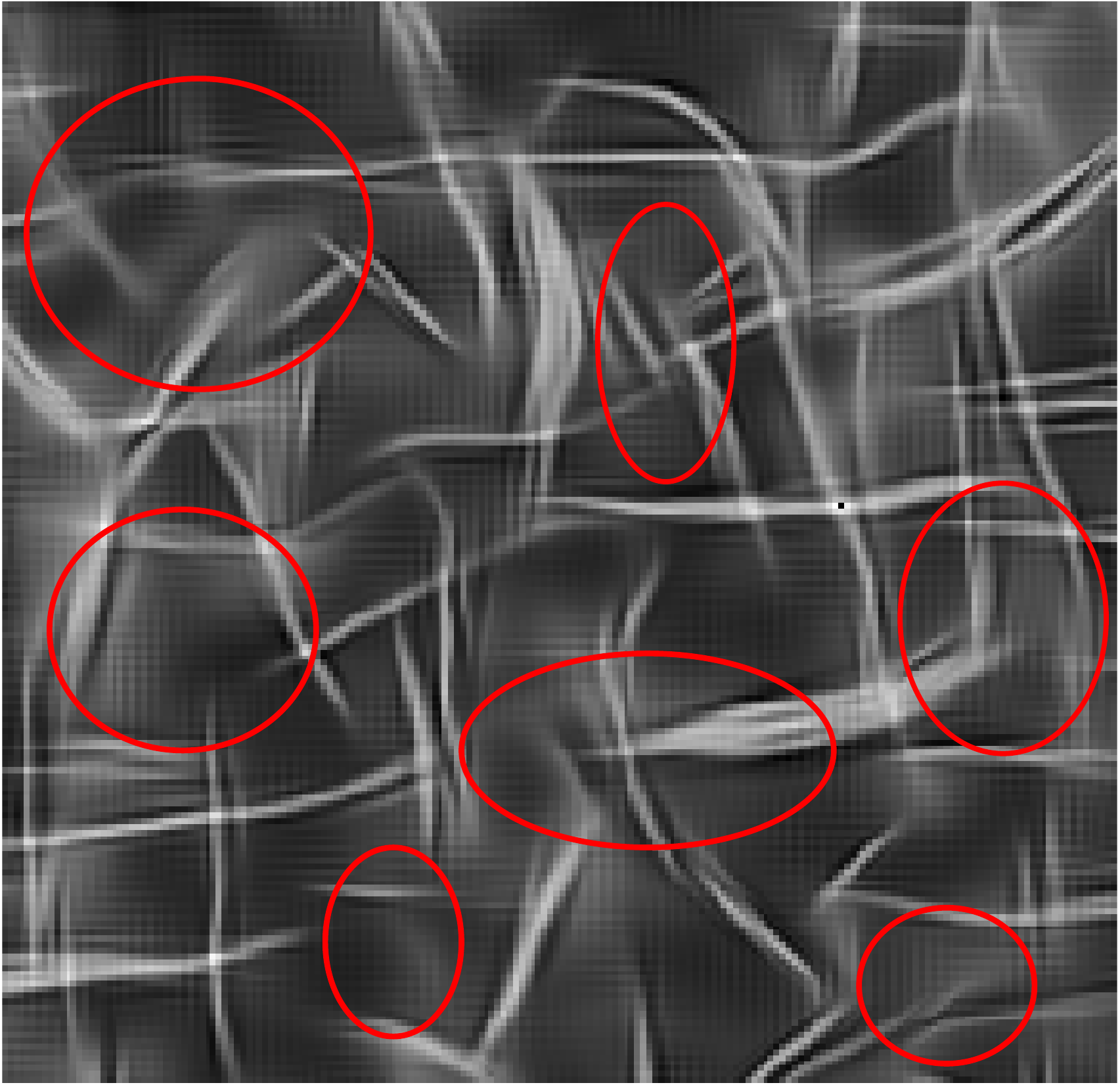}
                \caption{CED-OS}
        \end{subfigure}
        \begin{subfigure}[t]{0.35\textwidth}
                \centering
                \includegraphics[width=\textwidth]{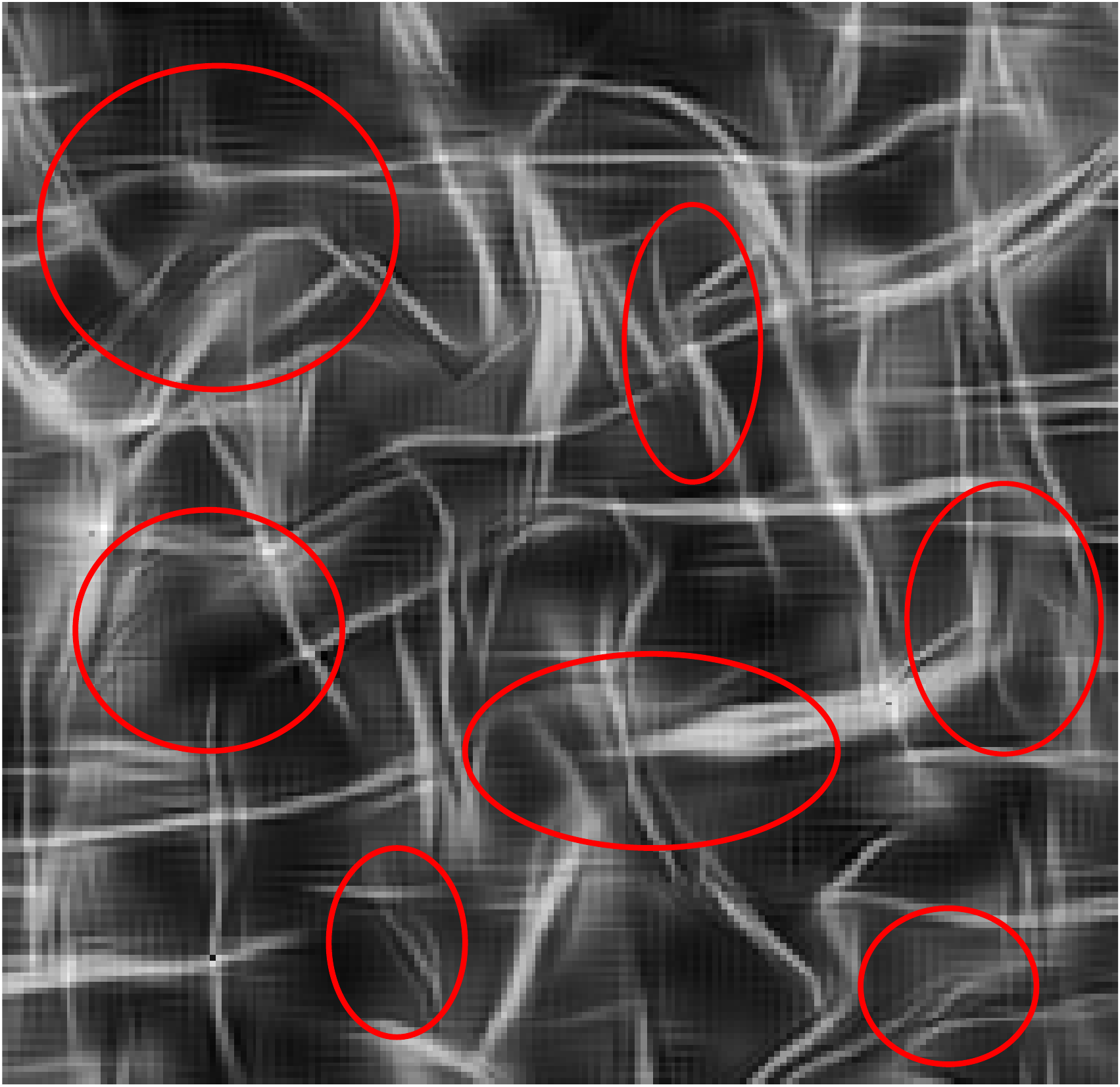}
                \caption{CED-SOS}
        \end{subfigure}

        \caption
{Top: Results of enhancement via of CED, CED-OS and CED-OS on Scale-OS on an image containing collagen fibres from a bone tissue. Bottom: Circled regions show loss of small-scale information in CED compared CED-SOS.
} 
\label{fig:Brodatz-Comp}
\end{figure*}

\begin{figure*}[h]
        \centering
%
  \begin{subfigure}[t]{0.35\textwidth}
                \centering
                \includegraphics[width=\textwidth]{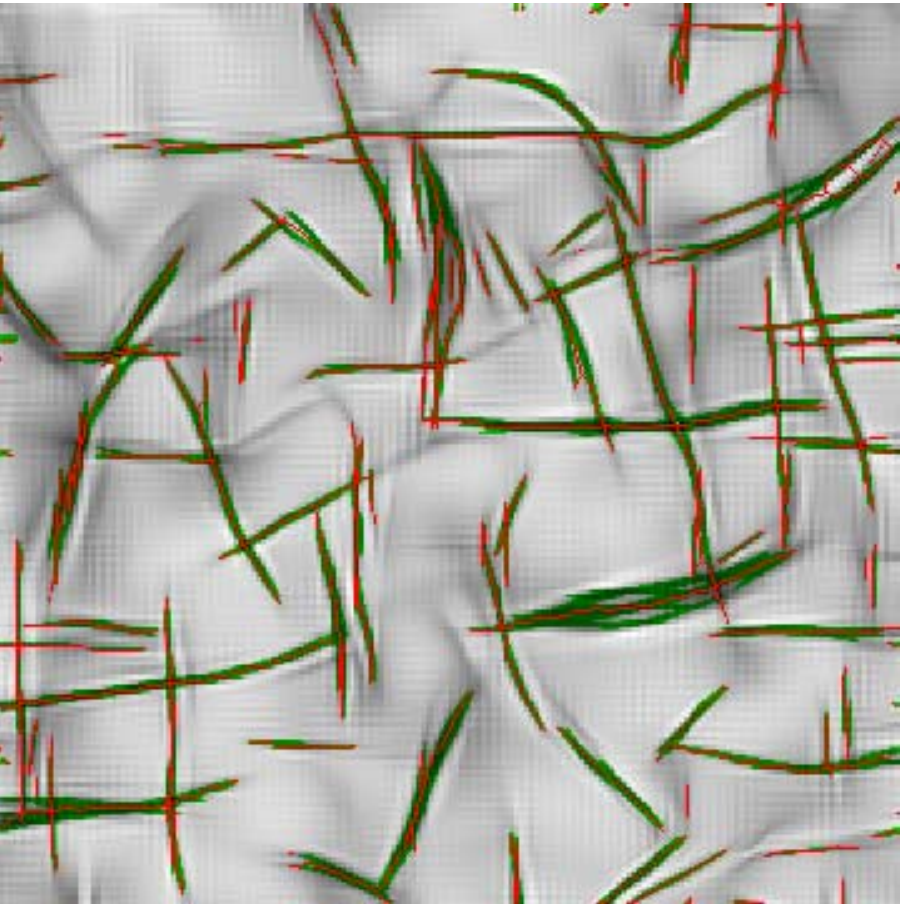}
                \caption{Vessel tracking on CED-OS}
        \end{subfigure}
        \begin{subfigure}[t]{0.35\textwidth}
                \centering
                \includegraphics[width=\textwidth]{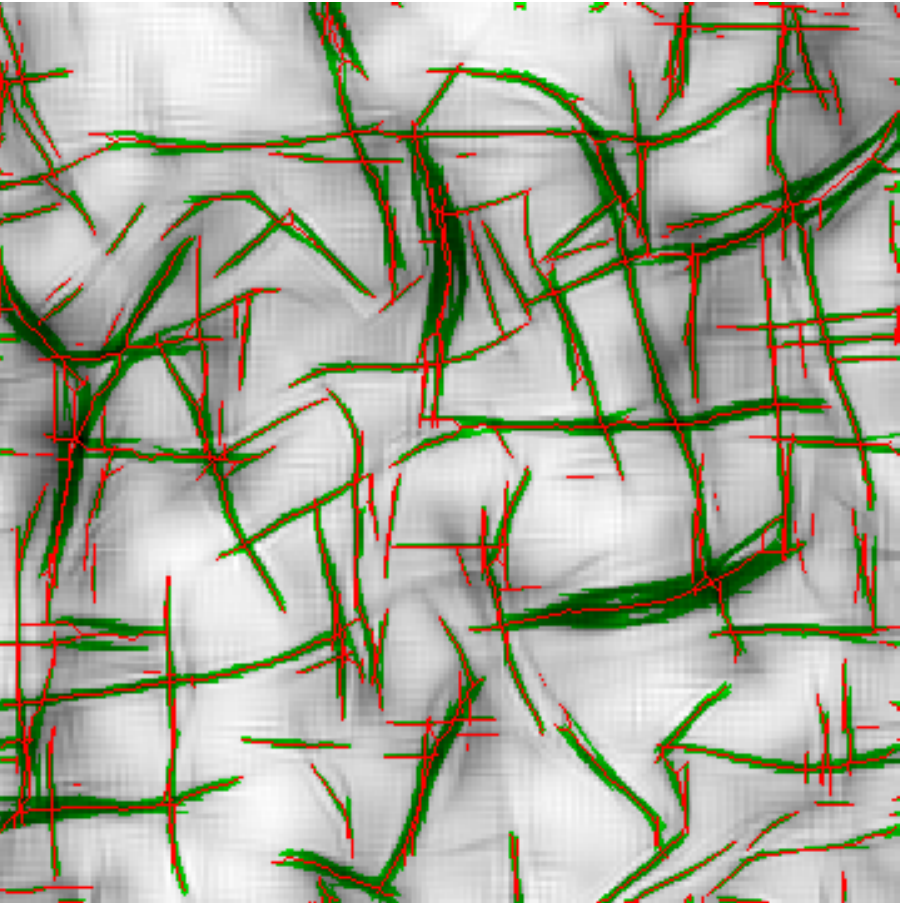}
                \caption{Vessel tracking on CED-SOS}
        \end{subfigure}
                
%
        \caption
{Results of a vessel tracking method \cite{Bekkers2012} on the images enhanced by CED-OS and CED-SOS. As expected tracking on CED-SOS enhanced images give better results.} 
\label{fig:Brodatz-Track}
\end{figure*}

For further comparison and in order to support the potential and wide applicability of our method in medical imaging, we also ran our algorithm on an image containing collagen fibers within bone tissue in Figure ~\ref{fig:Brodatz-Comp}. Here the loss of
small-scale data in CED-OS compared to CED-SOS is evident, an indeed a basic tracking algorithm (based on
[34] involving thresholding $\partial_{\eta}^{2}\mathcal{W}_{\psi}f$ and a basic morphological component method) reveals the benefit: 
``Tracking on CED-SOS outperforms tracking on CED-OS", see Figure~\ref{fig:Brodatz-Track}. 
On the other hand CED-OS has the advantage that it is faster and its implementation consumes far less memory (a factor in the order of the number of discrete scales).\\
%

\begin{figure*}[ht!]
        \centering
        \begin{subfigure}{0.48\textwidth}
          \centering
           \includegraphics[width=\textwidth]{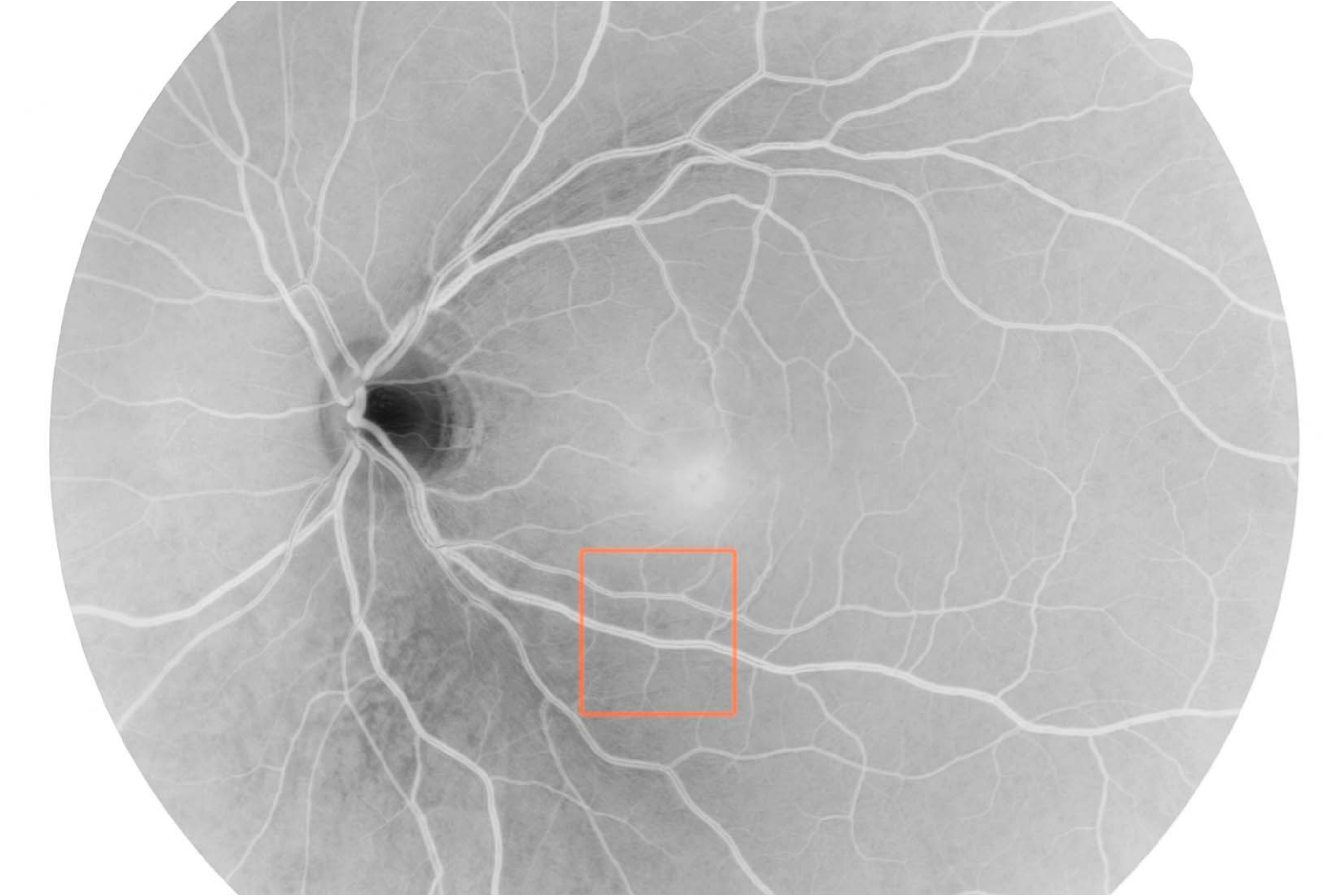}
            \caption{Optical image of the retina where bloodvessels depart from the optic disk}
        \end{subfigure}
        \begin{subfigure}{0.48\textwidth}
                \centering
                \includegraphics[width=\textwidth]{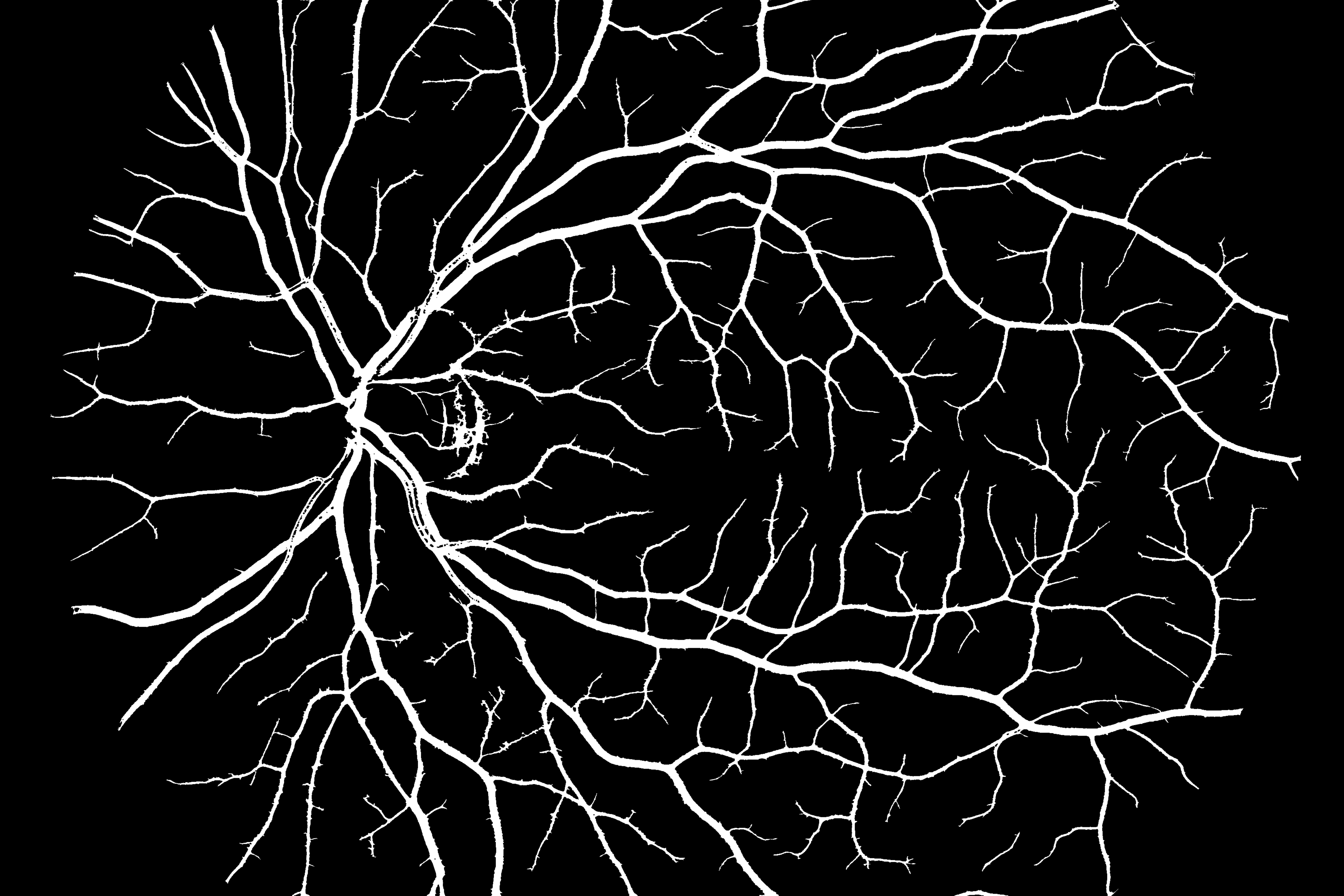}
                \caption{Segemented image via morphological components \cite{Hannink2014} based on the $SIM(2)$ vesselness filter.}
        \end{subfigure}
\begin{subfigure}{1\textwidth}
                \centering
                \includegraphics[width=\textwidth]{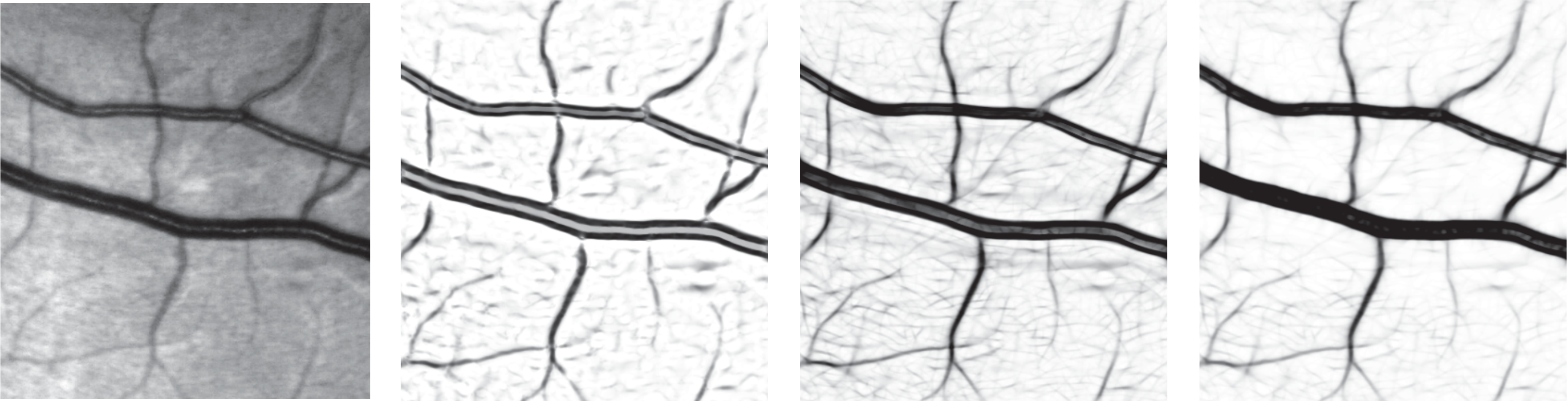}
                \caption{Comparision of vesselness filters. Left to Right: fragment of retinal image; Frangi vesselness filter \cite{Frangi1998} on image; $SIM(2)$ vesselnes filter with  frame $\partial_\theta,\partial_\xi,\partial_\eta, (\partial_\tau)$; $SIM(2)$ vesselnes filter with gauge frame $\partial_a,\partial_b,\partial_c, (\partial_\tau)$, see Eq.\eqref{eq:Gauge-Vesselness}.}
        \end{subfigure}
\begin{subfigure}{1\textwidth}
                \centering
                \includegraphics[width=\textwidth]{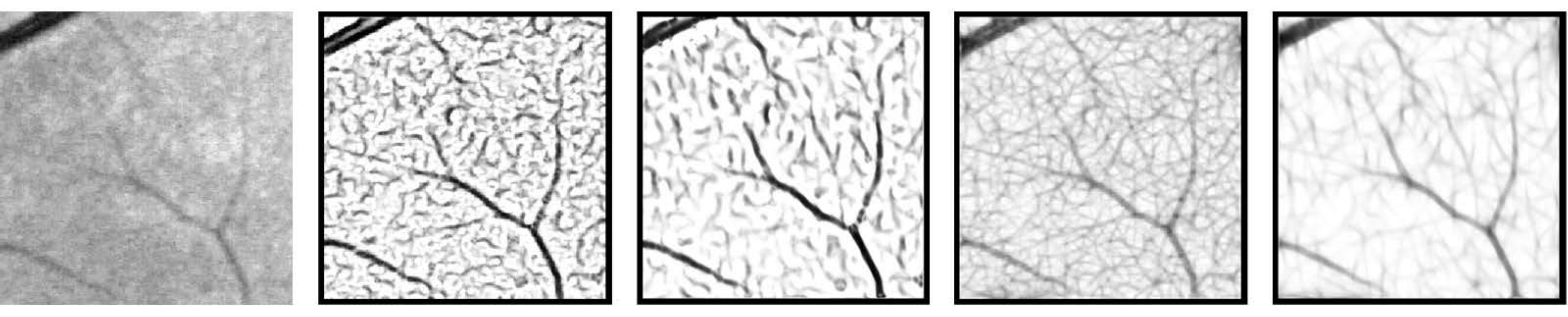}
                \caption{Comparision of vesselness filters on a tough retinal image patch with/without enhancements. Left to Right: typical fragment of retinal image where vesselness with standard parameters for the entire image fails; Frangi vesselness filter; Frangi vesselness filter after pre-enhancement via CED-SOS; $SIM(2)$ vesselnes filter; $SIM(2)$ vesselnes filter after pre-enhancement via CED-SOS.}
        \end{subfigure}
        \caption
{Figures (b) shows the result of $SIM(2)$ vesselness (Sec.\ref{sec:SIM-vesselness})  on an entire retinal fundus image (HRF database). See \cite{Hannink2014} for details. Figure (c) shows the advantage of $SIM(2)$ vesselness over commonly used Frangi approach. Figure (d) shows the advantage of enhancing the image prior to applying vesselness filters.} 
\label{fig:Julius-Vesselness-Results}
\end{figure*}

In Figure~\ref{fig:Julius-Vesselness-Results} we show some first results of $SIM(2)$ vesselness-filter described in Section~\ref{sec:SIM-vesselness}. These experiments illustrate
that the $SIM(2)$ vesselness filter ourperforms the multiple scale vesselness filtering acting directly in the image domain, as 
again no problems arise in areas with crossings/bifurcations.  For further underpinning via qualitative comparisons on the HRF-benchmark sequence see \cite{Hannink2014}. Note that Figure~\ref{fig:Julius-Vesselness-Results} in addition to presenting the benefits of our continuous wavelet transform on SIM(2) also shows the potential of including multiple scale crossing-preserving enhancements prior to the vesselness filtering (see Figure~\ref{fig:Julius-Vesselness-Results}(d) where we have selected a problematic patch). \\
%

\begin{figure}
\centering
\includegraphics[scale=0.4]{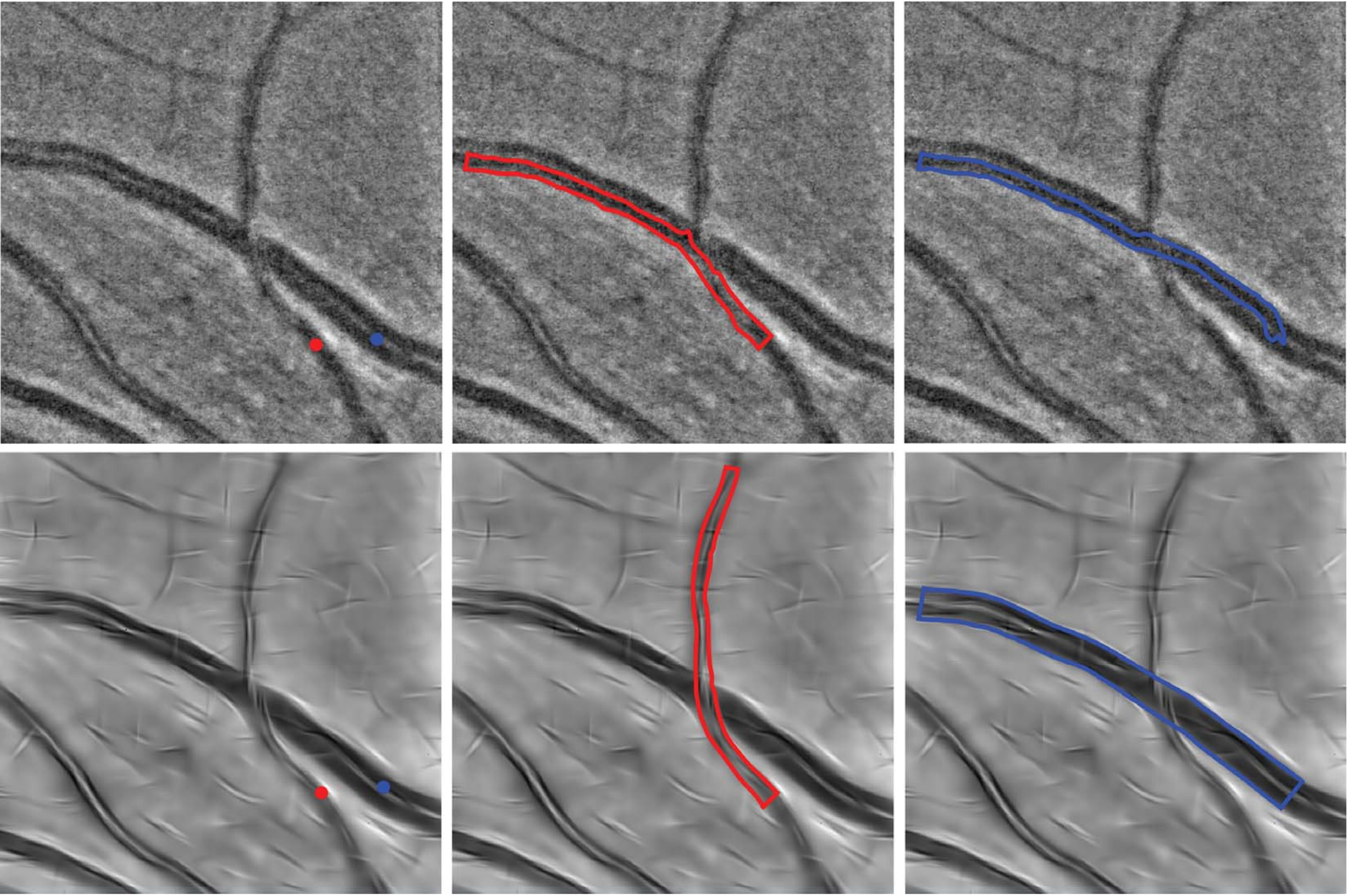}
\caption{Advantage of $SIM(2)$ crossing preserving multi-scale enhancement in the tracking of retinal vessels via state of the art tracking algorithm \cite{Bekkers2012}.
Top: Tracking applied to original image data fails due to noise. Bottom: Tracking applied after enhancement via CED-SOS is successful. }
\label{fig:Enh-Erik}
\end{figure}

In Figure~\ref{fig:Enh-Erik} we show the advantage of using $SIM(2)$ crossing preserving multi-scale enhancements CED-SOS prior to tracking retinal vessel via a  state of the art tracking algorithm \cite{Bekkers2012}.


\section{Conclusions}
There are two different tasks achieved in this article.
\begin{enumerate}
\item{Designing  an invertible score on $SIM(2)$, which also allows for accurate and efficient implementation of subsequent enhancement by contextual flows.}
\item{Construction of contextual  flows in the wavelet domain for the enhancement of elongated multi-scale crossing/bifurcating structures, via left-invariant PDEs on $SIM(2)$.}
\end{enumerate}

Regarding the first task we have presented  a generalized unitarity result for  algebraic affine Lie groups. This result is then used to design a multi-scale orientation score appropriate for subsequent left-invariant flows.

For the second task we have shown that only left-invariant flows on multi-scale orientation scores (defined on $SIM(2)$) robustly relate to Euclidean invariant operations on images. Furthermore we have  provided a differential-geometric and probabilistic interpretation of left-invariant PDEs on $SIM(2)$ which provides a strong intuitive rationale for the choice of left-invariant PDEs and involved diffusion parameters.

We have also derived analytic approximations of Green's function of linear diffusion on $SIM(2)$. Finally, we have presented crossing multiple-scale (curvature adaptive) flows via non-linear left-invariant diffusions on invertible scores.
Our preliminary results indicate that including the notion of scale in the framework of invertible orientation scores indeed has advantages over existing PDE techniques (CED, CED-OS). 

An interesting next step would be to include scale interactions in numerical implementation of non-linear diffusion on $SIM(2)$. Extending the framework of scattering operators on Lie groups (Mallat et al.  \cite{Mallat2012,Sifre2013}) to  connect scattering operators to PDEs on $SIM(2)$ is an interesting topic for future work.

\vspace{5mm}

\textbf{Acknowledgements} 
The authors would like to acknowledge Erik Bekkers (Technische Universiteit Eindhoven) for fruitful discussions and ideas regarding numerical implementation of  scale orientation score and tracking method used to create Figure~\ref{fig:Brodatz-Track}. The authors would like to thank Jiang Zhang (Technische Universiteit Eindhoven) and Julius Hannink for experiments to create Figure~\ref{fig:Enh-Erik} and Figure~\ref{fig:Julius-Vesselness-Results} respectively. The authors would also like to thank the anonymous referees
for valuable suggestions and comments.  
\appendix
\section{Proof of Theorem \ref{MPsiRecon}}\label{App:MpsiRecon}
For every $\Phi\in\mathbb{C}_K ^G$ has a $\mathcal{W}_\psi$ pre-image $f\in\mathbb{L}_{2}(\mathbb{R}^d)$, i.e. $f=(\mathcal{W}_\psi)^{-1}[\Phi]$ with 
\begin{align*}
\|\Phi\|_{\mathbb{C}_K^G}^2=\left(\Phi,\Phi\right)_{\mathbb{C}_K^G}=(f,f)_{\mathbb{L}_2(\mathbb{R}^d)}=\|f\|_{\mathbb{L}_2(\mathbb{R}^d)}^2.
\end{align*}  
We need to show that $\|\Phi\|^2_{\mathbb{C}_K^G}=\|\Phi\|^2_{M_\psi}:=(\Phi,\Phi)^{\frac{1}{2}}_{M_\psi}$, where we recall \eqref{MPsiInnProd} for the definition of the $M_\psi$-inner product. We have,
\begin{align}
\int\limits_{\mathbb{R}^d}\left|f(y)\right|^2dy=\int\limits_{\mathbb{R}^d}|\hat{f}(\omega)|^2d\omega=\int\limits_{\overline{\Omega}}|\hat{f}(\omega)|^2\frac{M_\psi(\omega)}{M_\psi(\omega)}d\omega=\int\limits_{\Omega}|\hat{f}(\omega)|^2\frac{M_\psi(\omega)}{M_\psi(\omega)}d\omega.\label{MPsiRecon_1}
\end{align}
Recall that $T$ is assumed to be a linear algebraic group (see \cite{Borel1991} for definition) and therefore it has locally closed (dual) orbits. This leads to $\Omega_c=\Omega_{cr}$ in \cite[Ch 5]{Fuhr2005} and since $\Omega_{cr}$ is open, \cite[Prop 5.7]{Fuhr2005}, we have that $\mu(\overline{\Omega}\backslash\Omega)=0$, where $\mu$ is the usual Lebesgue measure on $\mathbb{R}^d$, giving us the final equality in \eqref{MPsiRecon_1}. We can further write \eqref{MPsiRecon_1} as,
\begin{align*}
\int\limits_{\mathbb{R}^d}\left|f(y)\right|^2dy&=\int\limits_{\Omega}|\hat{f}(\omega)|^2\frac{1}{M_\psi(\omega)}\int\limits_T\frac{\left|(\mathcal{F}\mathcal{R}_t\psi)(\omega)\right|^2}{det(\tau(t))}d\mu_T(t) d\omega
=\int\limits_{\Omega}\bigg{\{}\int\limits_T\frac{
\overline{\hat{f}(\omega)\left(\mathcal{F}\mathcal{R}_t\psi(\omega)\right)}\hat{f}(\omega)\left(\mathcal{F}\mathcal{R}_t\psi(\omega)\right)}{det(\tau(t))}\bigg{\}}\frac{d\omega}{M_\psi(\omega)}\\
&=\int\limits_\Omega\int\limits_T\overline{\left(\mathcal{F}\Phi(\cdot,t)\right)(\omega)}\left(\mathcal{F}\Phi(\cdot,t)\right)(\omega)d\mu_T(t)\frac{d\omega}{M_\psi(\omega)}=\int\limits_\Omega\bigg{(}\int\limits_T\left|\mathcal{F}\Phi(\cdot,t)(\omega)\right|^2d\mu_T(t)\bigg{)}M_\psi^{-1}(\omega)d\omega.
\end{align*}
As a result $\Phi\in \mathbb{H}_{\psi}\otimes\mathbb{L}_2(T;\frac{d\mu_{T}(t)}{det(\tau(t))})$ and we have,
\begin{align*}
(f,f)_{\mathbb{L}_2(\mathbb{R}^d)}=(\mathcal{W}_\psi f,\mathcal{W}_\psi f)_{\mathbb{C}_K^G}=\int\limits_T\left(M_\psi^{-\frac{1}{2}}\mathcal{F}\Phi(\cdot,t),M_\psi^{-\frac{1}{2}}\mathcal{F}\Phi(\cdot,t)\right)d\mu_T(t),
\end{align*}
from which the result follows. $\hfill\Box$

\section{Proof of Theorem~\ref{thm:EuclideanInv}}\label{app:Covariance-Proof}
Note that since
\begin{equation*}
\widetilde{\mathcal{W}}_{\psi}[\mathcal{U}_{g}[f]](h)=(\mathcal{U}_{h}\psi,\mathcal{U}_{g}f)_{\mathbb{L}_{2}(\mathbb{R}^d)}=(\mathcal{U}_{g^{-1}h}
\psi,f)_{\mathbb{L}_{2}(\mathbb{R}^d)}=\mathcal{L}_{g}[\widetilde{\mathcal{W}}_{\psi} [f]](h)
\end{equation*}
we have 
\begin{equation}
\widetilde{\mathcal{W}}_{\psi}\mathcal{U}_{g}=
\mathcal{L}_{g}
\widetilde{\mathcal{W}}_{\psi}, \  \text{for all }g\in G.\label{ConvEq1}
\end{equation}
Moreover,
\begin{equation*}
(\widetilde{\mathcal{W}}_{\psi}
\mathcal{U}_{g}f,U)_{\mathbb{L}_{2}
(G)}=(\mathcal{L}_{g}
\widetilde{\mathcal{W}}_{\psi},U)_{\mathbb{L}_{2}
(G)}\Leftrightarrow(f,\mathcal{U}_{g^{-1}}
(\widetilde{\mathcal{W}}_{\psi})^* U)_{\mathbb{L}_{2}
(G)}=(f,(\widetilde{\mathcal{W}}_{\psi})^*
\mathcal{L}_{g^{-1}}U)_{\mathbb{L}_{2}(G)}
\end{equation*}
for all $U\in \mathbb{L}_{2}(G)$, $f\in \mathbb{L}_{2}(\mathbb{R}^d)$, $g\in G$ and therefore we have, 
\begin{equation}
\mathcal{U}_{g}(\widetilde{\mathcal{W}}_{\psi})^*=\mathcal{L}_{g}(\widetilde{\mathcal{W}}_{\psi})^*, \ \text{for all } g\in G.\label{ConvEq2}
\end{equation}
(Necessary condition) Assuming that $(\mathbb{P}_{\psi}\circ
\Phi)$ is left invariant it follows from \eqref{ConvEq3}, \eqref{ConvEq1} and \eqref{ConvEq2} that  
\begin{eqnarray}
\Upsilon[\mathcal{U}_{g}f]&=&(\widetilde{\mathcal{W}}
_{\psi})^{*}\circ\Phi\circ\widetilde{\mathcal{W}}
_{\psi}
\circ\mathcal{U}_{g} [f] \nonumber \\
&=&(\widetilde{\mathcal{W}}
_{\psi})^{*}\circ(\mathbb{P}_{\psi}\circ
\Phi)\circ\widetilde{\mathcal{W}}
_{\psi}
\circ\mathcal{U}_{g} [f]
\nonumber\\
&=& (\widetilde{\mathcal{W}}
_{\psi})^{*}\circ(\mathbb{P}_{\psi}\circ
\Phi)\circ\mathcal{L}_{g}\circ\widetilde{\mathcal{W}}_{\psi} [f]\nonumber \\
&=& (\widetilde{\mathcal{W}}
_{\psi})^{*}\circ\mathcal{L}_{g}
\circ(\mathbb{P}_{\psi}\circ
\Phi)\circ\widetilde{\mathcal{W}}_{\psi} [f]\\
&=& \mathcal{U}_{g}\circ(\widetilde{\mathcal{W}}
_{\psi})^{*}\circ(\mathbb{P}_{\psi}\circ
\Phi)\circ\widetilde{\mathcal{W}}
_{\psi}
[f]=\mathcal{U}_{g}[\Upsilon[f]]\nonumber
\end{eqnarray}
for all $f\in \mathbb{L}_{2}(\mathbb{R}^2)$ and $g\in G$. Thus we have $\Upsilon\mathcal{U}_{g}=\mathcal{U}_{g}
\Upsilon$ for all $g\in G$.

(Sufficient condition) Now suppose $\Upsilon$ is Euclidean invariant. Then again by \eqref{ConvEq1} and \eqref{ConvEq2} we have that,
\begin{equation*}
 (\widetilde{\mathcal{W}}
_{\psi})^{*}\circ\Phi\circ\mathcal{L}_{g}
\circ\widetilde{\mathcal{W}}_{\psi} [f]
= (\widetilde{\mathcal{W}}
_{\psi})^{*}\circ\mathcal{L}_{g}
\circ\Phi\circ\widetilde{\mathcal{W}}_{\psi} [f],
\end{equation*}
for all $f\in \mathbb{L}_{2}(\mathbb{R}^2)$ and $g\in G$. Since the range of $\mathcal{L}_{g}\big{|}_{\mathbb{C}_{K}
^{G}}$ and the range of $(\mathbb{P}_{\psi}\circ
\Phi)$ is contained in $\mathbb{C}_{K}^{G}$ and since $(\widetilde{\mathcal{W}}
_{\psi})^{*}\big{|}_{\mathbb{C}_{K}^{G}}
=(\mathcal{W}
_{\psi})^{-1}$, we have $(\mathbb{P}_{\psi}\circ
\Phi)\circ\mathcal{L}_{g}
\circ\widetilde{\mathcal{W}}_{\psi} 
= \mathcal{L}_{g}\circ(\mathbb{P}_{\psi}\circ
\Phi)\circ\widetilde{\mathcal{W}}
_{\psi}$. As the range of $\widetilde{\mathcal{W}}_{\psi}$ equals $\mathbb{C}_{K}^{G}$, we have, $\mathcal{L}_g\circ(\mathbb{P}_{\psi}\circ
\Phi)=(\mathbb{P}_{\psi}\circ
\Phi)
\circ\mathcal{L}_{g}$ for all $g\in G$. 

\section{Explicit formulation of exponential and logarithm curves}\label{App:ExpCurvesExplicit}
The exponential curves $g_{0}\gamma_c(t)=(x(t),y(t),\tau(t),\theta(t))$ passing through $g_{0}=(x_{0},y_{0},e^{\tau_{0}},\theta_{0})$ at $t=0$ is,
\begin{eqnarray}
x(t)&=&
\frac{1}{c_1^2+c_4^2}
[ e^{\tau _0} c_1 \left(\left(-\sin\left[\theta 
_0\right]+e^{t c_4} \sin\left[t c_1+\theta 
_0\right]\right) c_2+\left(-\cos\left[\theta 
_0\right]+e^{t c_4} \cos\left[t c_1+\theta 
_0\right]\right) c_3\right)
+c_1^2 x_0+ \nonumber\\
&&c_4 \left(e^{\tau _0} \left(-\cos\left[\theta 
_0\right]+e^{t c_4} \cos\left[t c_1+\theta 
_0\right]\right) c_2+e^{\tau _0} 
\left(\sin\left[\theta _0\right]-e^{t c_4} 
\sin\left[t c_1+\theta _0\right]\right) c_3+c_4 
x_0\right)]\nonumber\\
y(t)&=&
\frac{1}{{c_1^2+c_4^2}}[e^{\tau _0} c_1 \left(\left(\cos\left[\theta _0\right]-e^{t c_4} \cos\left[t c_1+\theta _0\right]\right) c_2+\left(-\sin\left[\theta _0\right]+e^{t c_4} \sin\left[t c_1+\theta _0\right]\right) c_3\right)+c_1^2 y_0+ \nonumber\\
&&c_4 \left(e^{\tau _0} \left(-\sin\left[\theta _0\right]+e^{t c_4} \sin\left[t c_1+\theta _0\right]\right) c_2+e^{\tau _0} \left(-\cos\left[\theta _0\right]+e^{t c_4} \cos\left[t c_1+\theta _0\right]\right) c_3+c_4 y_0\right)]
\nonumber\\
\tau(t)&=& tc_4+\tau_{0} \nonumber\\
\theta(t)&=&
  tc_1+\theta _0.\label{2.9}
\end{eqnarray}
To explicitly determine the $\log$ map, we solve for $\{c^{1},c^{2},c^{3},c^{4}\}$
from the equality, $g=\exp\left( t\sum\limits_{i=1}^{4}c^{i}
A_{i}\right)$, where $g\in SIM(2)$ and $c^{i}$'s are as defined earlier. This is achieved by substituting $g_{0}=e$, i.e.\ $x_{0}=y_{0}=\theta_{0}=\tau_{0}=0$, in  \eqref{2.9} yielding,
\begin{align}
c^1=\theta /t, \ 
c^2=\frac{(y \theta -x \tau) +(-\theta\eta+\tau\xi)}{t \left(1+e^{2 \tau }-2 e^{\tau } \cos\theta\right)}, \ 
c^3=\frac{-(x \theta +y \tau) +(\theta\xi+\tau\eta)}{t \left(1+e^{2 \tau }-2 e^{\tau } \cos\theta\right)}, c^4=\tau /t,\label{SIM2_11}
\end{align}
where we have made use of the definition of $\xi$ and $\eta$,
\begin{equation*}
\xi=e^{\tau}(x \cos\theta + y \sin\theta), \ \eta=e^{\tau}(-x \sin\theta + y \cos\theta).
\end{equation*}

\section{Proof of Theorem~\ref{thm:Exp_Char}}\label{App:Exp_Char}
We first note that, $U\in\mathcal{D}(d\mathcal{R}(A))\Rightarrow (d\mathcal{R}(A)U)\in H$, where $H=\mathbb{L}_{2}(SIM(2))$. Then, $\mathcal{R}_{e^{tA}}U\in\mathcal{D}(d\mathcal{R}(A))$ follows as $[d\mathcal{R}(A)\mathcal{R}_{e^{tA}}](U)=[\mathcal{R}_{e^{tA}}d\mathcal{R}(A)](U),  \text{ for all }U\in\mathcal{D}(d\mathcal{R}(A))$.
\begin{align*}
[d\mathcal{R}(A)\mathcal{R}_{e^{tA}}](U)(g)&=\left[\lim\limits_{h\rightarrow 0}\frac{\mathcal{R}_{e^{hA}}-I}{h}\right](\mathcal{R}_{e^{tA}}U)(g)=\lim\limits_{h\rightarrow 0}\frac{\mathcal{R}_{e^{hA}}\mathcal{R}_{e^{tA}}(U(g))-\mathcal{R}_{e^{tA}}(U(g))}{h}\\
&=\lim\limits_{h\rightarrow 0}\frac{\mathcal{R}_{e^{tA}}\mathcal{R}_{e^{hA}}(U(g))-\mathcal{R}_{e^{tA}}(U(g))}{h}=[\mathcal{R}_{e^{tA}}d\mathcal{R}(A)](U)(g),
\end{align*}
 where we have used $\mathcal{R}_{h}(U(g))=U(gh)$ and $e^{tA}e^{hA}= e^{(t+h)A}=e^{hA}e^{tA}$. This proves (1).
 We need to prove that $e^{t d\mathcal{R}(A)}U=\mathcal{R}_{e^{tA}}U$, where $U\in \mathcal{D}(d\mathcal{R}(A))$, $\text{ for all } A\in T_{e}(G),  \text{ for all } t>0$.
By definition $e^{t d\mathcal{R}(A)}U$ is the (strong) solution for $\frac{\partial W}{\partial t}(\cdot,t)=d\mathcal{R}(A)W(\cdot,t)$ with $W(g,0)=U$. Now $\frac{\partial (\mathcal{R}_{e^{tA}}U)}{\partial t}=d\mathcal{R}(A)(\mathcal{R}_{e^{tA}}U)$, follows from,
\begin{align}
\frac{\partial (\mathcal{R}_{e^{tA}}U)}{\partial t}&=\lim\limits_{h\rightarrow 0}\frac{\mathcal{R}_{e^{(t+h)A}}U-\
\mathcal{R}_{e^{tA}}U}{h}   =
\lim\limits_{h\rightarrow 
0}\frac{\mathcal{R}_{e^{hA}} \mathcal{R}_{e^{tA}}U -\
\mathcal{R}_{e^{tA}}U}{h}\nonumber\\
&=\lim\limits_{h\rightarrow 0}\left[
\left( \frac{\mathcal{R}_{e^{hA}}-I}{h}
\right)\mathcal{R}_{e^{tA}}U
\right] = d\mathcal{R}(A)
(\mathcal{R}_{e^{tA}}U).
\label{SIM2_13}
\end{align}
Note that the limit used above is defined on the space $H=\mathbb{L}_{2}(SIM(2))$.

Let $d\mathcal{R}(A)=\mathcal{A}=\sum\limits_{i=1}^{4}c^{i}
\mathcal{A}_{i}$. Then the solution for \eqref{SIM2_14} is, $W(\cdot,t)=e^{-t  d\mathcal{R}(A)}U=\mathcal{R}_{e^{-tA}}U$, where the second equality follows from (2). Thus we have,
\begin{align}
W(g,t)=U(g e^{-tA}) \ \ \forall g\in SIM(2), \ \forall t>0 \Rightarrow W(g_{0}  e^{tA},t)=U(g_{0})\ \ \forall t>0.
\end{align}
From which the result follows.$\hfill\Box$

\section{Approximation of $SIM(2)$ by a Nilpotent Group via Contraction}\label{App:MethContraction} 
Following the general framework by ter Elst and Robinson \cite{Elst1998}, which involves semigroups on Lie groups generated by subcoercive operators, we consider a particular case by setting the Hilbert space $H=\mathbb{L}_{2}(SIM(2))$, the group $G=SIM(2)$ and the right-regular representation $\mathcal{R}$. Furthermore we consider the algebraic basis $\{\mathcal{A}_{1}=\partial _{\theta},\mathcal{A}_{2}=\partial _{\xi},\mathcal{A}_{4}=\partial _{\tau}\}$ leading to the following filtration of the Lie algebra
\begin{align}
\mathfrak{g}_{1}:=\text{span}\{\mathcal{A}_{1},\mathcal{A}_{2},\mathcal{A}_{4}\}
\subset \mathfrak{g}_{2}=[\mathfrak{g}_{1},\mathfrak{g}_{1}]=\text{span}\{\mathcal{A}_{1},\mathcal{A}_{2},\mathcal{A}_{3},\mathcal{A}_{4}\}=\mathcal{L}(SIM(2)).
\end{align}
Based on this filtration we assign the following weights to the generators:
\begin{equation}
w_{1}=w_{2}=w_{4}=1 \text{ and } w_{3}=2.
\label{Evo_3}
\end{equation}
For e.g.\ $w_{2}=1$ since $\mathcal{A}_{2}$ occurs first in $\mathfrak{g}_{1}$, while $w(3)=2$ since $\mathcal{A}_{3}$ occurs in $\mathfrak{g}_{2}$ and not $\mathfrak{g}_{1}$. Based on these weights we define the following dilations on the Lie algebra $T_{e}(SIM(2))$ (recall $A_{i}=\mathcal{A}_{i}\large{|}_{e}$),
\begin{align*}
&\gamma_{q}\left(\sum\limits_{i=1}^{4}c^{i}A_{i}
\right)=\sum\limits_{i=1}^{4}q^{w_{i}}c^{i}A_{i}, \ \forall c^{i}\in\mathbb{R},\\
&\tilde{\gamma}_{q}(x,y,\tau,\theta)=\left(
\frac{x}{q^{w_{2}}},\frac{y}{q^{w_{3}}},\frac{\tau}{q^{w_{4}}},e^{i\frac{\theta}{q^{w_1}}}
\right),
\end{align*}
with the weights $w_{i}$ defined in \eqref{Evo_3}, and for $0<q\leq 1$ we define the Lie product $[A,B]_{q}=\gamma _{q}^{-1}[\gamma _{q}(A),\gamma _{q}(B)]$. Let $(SIM(2))_t$ be the simply connected Lie group generated by the Lie algebra\footnote{Note that $(SIM(2))_q=exp_{q}(T_{e}(SIM(2)))$.} $(T_{e}(SIM(2)),[\cdot,\cdot]_{q})$. The group products on these intermediate groups $(SIM(2))_{q\in(0,1]}$ are given by,
\begin{equation}
(x,y,\tau,\theta)\cdot_{q}(x',y',\tau ',\theta ')=(x+e^{\tau q}[\cos(q\theta)x'-q \sin(q\theta)y'],y+e^{\tau q}[\frac{\sin(q\theta)}{q} x'+ \cos(q\theta)y'],\tau+\tau ',\theta +\theta ').
\end{equation}
The dilation on the Lie algebra coincides with the pushforward of the dilation on the group $\gamma _{q}=(\tilde{\gamma}_q)_*$  and therefore the left-invariant vector fields on $(SIM(2))_t$  are given by
\begin{equation*}
\mathcal{A}_{i}^{q}\large{|}_{g}=
(\tilde{\gamma}_{q}^{-1}\circ L_{g}\circ \tilde{\gamma}_{q})_{*}A_{i},
\end{equation*}
for all $q\in(0,1]$ which leads to,
\begin{align*}
\mathcal{A}_{i}^{q}\large{|}_{g}\phi &=(\tilde{\gamma}_{q}^{-1}\circ L_{g})_{*}(\tilde{\gamma}_q)_{*}A_{i}\phi=
(\tilde{\gamma}_{q}^{-1}\circ L_{g})_{*}(\gamma _q)(A_{i})\phi =t^{w_{i}}(\tilde{\gamma}_{q}^{-1}\circ L_{g})_{*}(A_{i})\phi\\
&=q^{w_{i}}(\tilde{\gamma}_{q}^{-1})_{*}\mathcal{A}_{i}
\large{|}_{g}\phi =q^{w_{i}}\mathcal{A}_{i}
\large{|}_{\tilde{\gamma}_{q}^{-1}g}(\phi \circ \tilde{\gamma}_{q}),
\end{align*}
for all smooth complex valued functions $\phi$ defined on a open neighbourhood around $g\in SIM(2)$. So for all $g=(x,y,\tau,\theta)\in SIM(3)$ we have,
\begin{align}
\mathcal{A}_{1}^{t}\large{|}_{g}&=q(\frac{1}{q}\partial _{\theta})=\partial _{\theta} \nonumber\\
\mathcal{A}_{2}^{q}\large{|}_{g}&=q\left[  e^{\tau q}\left(\frac{\cos(q\theta)}{q}\partial _{x}+\frac{\sin(q\theta)}{q^2}\partial _{y}
\right)
\right]=
e^{\tau q}\left(\cos(q\theta)\partial _{x}+\frac{\sin(q\theta)}{q}\partial _{y}
\right)\\
\mathcal{A}_{3}^{q}\large{|}_{g}&=q^{2}
\left[  e^{\tau q}\left(-\frac{\sin(q\theta)}
{q}\partial _{x}+\frac{\cos(q\theta)}{q^2}\partial _{y}
\right)
\right]=
e^{\tau q}\left(-q\sin(q\theta)\partial _{x}+ \cos(q\theta)\partial _{y}
\right)\nonumber \\
\mathcal{A}_{4}^{q}\large{|}_{g}&=\partial _{\tau}.\nonumber
\end{align}
Furthermore, $[A_{i},A_{j}]_{q}=\gamma _{q}^{-1}[\gamma _{q}(A_{i}),\gamma _{q}(A_{j})]= \gamma _{q}^{-1}q^{w_{i}+w_{j}}[A_{i},A_{j}]=\sum\limits_{k=1}^{4}q^{w_{i}
+w_{j}-w_{k}}c_{ij}^{k}A_{k}$ and therefore
\begin{align}
[A_{1},A_{2}]_{q}=A_{3},  \ \ \ [A_{1},A_{3}]_{q}=-q^{2}A_{2}, \ \ \  [A_{4},A_{2}]_{q}=q A_{2},  \ \ \  [A_{4},A_{3}]_{q}=q A_{3}.\label{Evo_4}
\end{align}
Analogously to the case $q=1$, $(SIM(2))_{q=1}=SIM(2)$, there exists an isomorphism of the Lie algebra at the unity element $T_{e}((SIM(2))_{q})$ and the left-invariant vector fields on the group $\mathcal{L}((SIM(2))_{q})$:
\begin{equation}
(A_{i}\leftrightarrow \mathcal{A}_{i}^{q} \text{ and }A_{j}\leftrightarrow \mathcal{A}_{j}^{q})\Rightarrow  [A_{i},A_{j}]_{q}\leftrightarrow [\mathcal{A}_{i}^{q},\mathcal{A}_{j}^{q}].
\end{equation}
It can be verified  that the left invariant vector fields $\mathcal{A}_{i}^{q}$ satisfy the same commutation relations as \eqref{Evo_4}.
In the case $H\equiv\lim\limits_{q\downarrow 0} \  (SIM(2))_{q}$ the left-invariant vector fields are given by
\begin{align}
\mathcal{A}_{1}^{0}=\partial _{\theta}, \  \mathcal{A}_{2}^{0}=\partial _{x}+\theta\partial _{y}, \  \mathcal{A}_{3}^{0}=\partial _{y},  \ 
\mathcal{A}_{4}^{0}=\partial _{\tau}.
\label{eq:App-Contarction-Shear}
\end{align}
So, the homogeneous nilpotent contraction Lie group equals 
\begin{equation}
H_{3}=\lim\limits_{q\downarrow 0} \  
(SIM(2))_{q} \text{ and } 
SIM(2)=(SIM(2))_{q=1}/(\{0\}\times\{0\}
\times\{0\}\times 2\pi\mathbb{Z}),
\end{equation}
with the Lie algebra $\mathcal{L}(H)=\text{span}\{\partial _{\theta},\partial _{x} +\theta\partial _{y},\partial _{y},\partial _{\tau}\}$ and 
$\mathcal{L}(SIM(2))=\text{span}\{\partial _{\theta},\partial _{\xi},\partial _{\eta},\partial _{\tau}\}$.

\section{Differential-geometric interpretation of Left-invariant Evolutions}\label{sec:App_DiffGeomInter}
In order to keep track of orthogonality and parallel transport in such diffusions we need an invariant first fundamental form $\mathcal{G}$ on $SIM(2)$, rather than the trivial, bi-invariant (i.e.\ left and right invariant), first fundamental form on $(\mathbb{R}^4,T(\mathbb{R}^4))$, where each tangent space $T_{\mathbf{x}(\mathbb{R}^4)}$ is identified with $T_{0}(\mathbb{R}^4)$ by standard parallel transport on $\mathbb{R}^4$, i.e.\ $\mathcal{G}_{\mathbb{R}^4}(\mathbf{x},\mathbf{y})=\mathbf{x}\cdot\mathbf{y}=
x^1y^1+x^2y^2+x^3y^3+x^4y^4$. 

Recall Theorem~\ref{thm:EuclideanInv} which essentially states that operators $\Phi$ on scale-OS should be left-invariant (i.e.$\ \mathcal{L}_g\circ\Phi=\Phi\circ\mathcal
{L}_g$) and not right-invariant in order to ensure that the effective operator $\Upsilon_{\psi}$ on the image is a Euclidean invariant operator. This suggests that the first fundamental form required for our diffusions on $SIM(2)$ should be left-invariant. The following theorem characterizes the formulation of a left-invariant first fundamental form w.r.t the $SIM(2)$ group.
\begin{theorem}
\label{thm:MetricSIM}
The only real valued left-invariant (symmetric, positive, semidefinite) first fundamental form $\mathcal{G}:SIM(2)\times T(SIM(2))\times T(SIM(2)\rightarrow\mathbb{C}$ on $SIM(2)$ are given by,
\begin{align}
\mathcal{G}=\sum\limits_{i=1}^4 \sum\limits_{j=1}^4 g_{ij} \  \omega^i   \otimes \omega^j   , \  g_{ij}\in\mathbb{R}. 
\end{align}  
where the dual basis $\{\omega^{1},\omega^{2},\omega^{3},\omega^{4}\} \subset(\mathcal{L}(SIM(2)))^*$ of the dual space $(\mathcal{L}(SIM(2)))^*$ of the vector space $\mathcal{L}(SIM(2))$ of left-invariant vector fields spanned by 
\begin{align}
\{\mathcal{A}_{1},\mathcal{A}_{2},
\mathcal{A}_{3},\mathcal{A}_{4}\}=\{\partial_{\theta}
,e^{\tau}(\cos\theta\partial_{x}+
\sin\theta\partial_{y}),e^{\tau}(-\sin\theta\partial_{x}+
\cos\theta\partial_{y}),\partial_{\tau}\},
\label{Evo_10}
\end{align}
obtained by applying the operator $d\mathcal{R}:T_{e}(G)\rightarrow\mathcal{L}(G)$, defined as
\begin{align}
(d\mathcal{R}(A)\phi)(g)=\lim\limits_{t\downarrow 0}\frac{(\mathcal{R}_{exp(tA)}\phi)(g)-\phi(g)}{t}, \ A\in T_{e}(G),\phi\in\mathbb{L}_{2}(G),g\in G, 
\end{align}
to the standard basis in the Lie algebra
\begin{align}
\{A_{1},A_{2},A_{3},A_{4}\}=\{\partial
_{\theta},\partial_{x},\partial_{y},\partial_{a}\}\subset T_{e}(SIM(2)),
\end{align}
is given by
\begin{align}
\{\omega^{1},\omega^{2},
\omega^{3},\omega^{4}\}=
\{d\theta ,\frac{1}{e^{\tau}} (\cos\theta dx+\sin\theta dy),\frac{1}{e^{\tau}} (-\sin\theta dx+\cos\theta dy),d\tau\}.
\label{Evo_7}
\end{align}
\end{theorem}
The proof, see \cite{Sharma2013} for details, is similar to the $SE(2)$ case \cite{Duits2010a}. 

We consider the Maurer-Cartan form on $SIM(2)$ (discussed in the remainder of the subsection) and impose the following left-invariant, first fundamental form $\mathcal{G}_{\beta}:SIM(2)\times T(SIM(2)) \times T(SIM(2))\rightarrow\mathcal{C}$ on $SIM(2)$,
\begin{align}
\mathcal{G}=\sum\limits_{i=1}^4 \sum\limits_{j=1}^4 g_{ij} \  \omega^i   \otimes \omega^j=
\alpha^2\omega^1\otimes \omega^1 +
\beta^2 \omega^2 \otimes \omega^2 +
\beta^2 \omega^3\otimes \omega^3+
\sigma^2 \omega^4\otimes \omega^4.
\label{Evo_18}
\end{align}
Here $\sigma$ tunes the cost of changing scale, $\beta$ the cost of moving spatially and $\alpha$ the cost of moving angularly. By homogeneity we set $\alpha=1$, as it is only $\frac{\beta}{\alpha}$, $\frac{\sigma}{\alpha}$ that matter.
In order to understand the meaning of the next two theorems we need some definitions from differential geometry. See \cite{Spivak1970} for more details on these concepts.
\begin{definition}
Let $M$ be a smooth manifold and $G$ be a Lie group. A principle fiber bundle $P_{G}:=(P,M,\pi,R)$ above a manifold $M$ with structure group $G$ is a tuple $(P,M,\pi,R)$ such that $P$ is a smooth manifold (called the total space of the principle bundle), $\pi:P\rightarrow M$ is a smooth projection map with $\pi(P)=M$ and $\pi(u\cdot a)=\pi(a), \ \forall u\in P, \ a\in G$, $R$ a smooth right action $R_{g}p=p\cdot g, \ p\in P, \ g,h\in G$. Finally it should satisfy the "local triviality" condition:
For each $p\in M$ there is a neighbourhood $U$ of $p$ and a diffeomorphism $t:\pi^{-1}(U)\rightarrow U\times G$ of the form
$t(u)=(\pi(u),\phi(u))$ where $\phi$ satisfies $\phi(u\cdot a)=\phi(u)a$ where the latter product is in $G$.
\end{definition}
\begin{definition}
\label{def:CarEhrConn}
A principle fiber bundle $P_{G}:=(P,M,\pi,R)$ is commonly equipped with a Cartan-Ehresmann connection form $\omega$. This by definition is a Lie algebra $T_{e}(G)$-valued $1$-form $\omega:P\times T(P)\rightarrow T_{e}(G)$ on $P$ such that
\begin{align}
&\omega(d\mathcal{R}(A))=A \text{ for all } A\in T_{e}(G)\nonumber\\
&\omega((R_h)_{*}\mathcal{A})=Ad(h^{-1})
\omega(\mathcal{A}) \text{ for all vector fields } \mathcal{A} \text{ in } G.
\end{align}
It is also common practice to relate principle fiber bundles to vector bundles. Here one uses an external representation $\rho:G\rightarrow F$ into a finite-dimensional vector space $F$ of the structure group to put an appropriate vector space structure on the fibers $\{\pi^{-1}(m)| \ m\in M\}$ in the principle fiber bundles.  
\end{definition}
\begin{definition}
\label{def:AssVecBundle}
Let $P$ be a principle fiber bundle with finite dimensional structure group $G$. Let $\rho:G\rightarrow F$ be a representation in a finite-dimensional vector space $F$. Then the associated vector bundle is denoted by $P\times_{\rho} F$ and equals the orbit space under the right action 
\begin{align*}
(P\times F)\times G\rightarrow P\times F \text{ given by }((u,X),g)\mapsto (ug,\rho(g)X),
\end{align*}
for all $g\in G, \ X\in F$ and $u\in P$.
\end{definition}
\begin{remark}
$GL(T_e(SIM(2)))$ denotes the collection of linear operators on the Lie algebra $T_{e}(SIM(2))$. Note that each linear operator $\overline{Q}\in GL(T_e(SIM(2)))$ on $T_e(SIM(2))$ is $1$-to-$1$ related to bilinear form $Q$ on $(T_e(SIM(2)))^*\times T_e(SIM(2))$ by means of
\begin{align*}
\langle B,\overline{Q}A\rangle =Q(B,A), \text{ for all } B\in (T_e(SIM(2)))^*, \ A\in T_e(SIM(2)) \text{ and }\overline{Q}=\sum\limits_{i=1}^4Q(\omega^i \lvert_e,\cdot)A_i.
\end{align*}
So a basis for $GL(T_e(SIM(2)))$ is given by 
$\{\overline{\omega ^i\lvert_e\otimes A_j} \lvert
\ i,j\in\{1,2,3,4\}\}$. For the simplicity 
of notation we omit the overline and write 
$\omega^i\lvert_e\otimes A_j$ as it is clear 
from context whether we mean the bilinear 
form or the linear mapping.
\end{remark}
\begin{theorem}
\label{thm:Res1}
The Maurer-Cartan form $\omega$ on $SIM(2)$ can be formulated as 
\begin{align}
\omega_{g}(X_g)=\sum\limits_{i=1}^{4}\langle
\omega^i\lvert_{g},X_{g}\rangle A_{i}, \ X_{g}\in T_{g}(SIM(2)),
\label{Evo_13}
\end{align}
where $\{\omega^{i}\}_{i=1}^{4}$ is given by 
\eqref{Evo_7} and $A_i=\mathcal{A}_{i}\lvert_{e}$; recall \eqref{Evo_10}. It is a Cartan-Ehresmann connection form on the principle fiber bundle $P=(SIM(2),e,SIM(2),\mathcal{L}(SIM(2)))$, where $\pi(g)=e$, $R_{g}u=ug, \ u,g\in SIM(2)$. Let $Ad$ denote the adjoint action of $SIM(2)$ on its own Lie algebra $T_{e}(SIM(2))$, i.e.\ $Ad(g)=(R_{g^{-1}}L_g)_*$, i.e.\ the push-forward of conjugation. Then the adjoint representation of $SIM(2)$ on the vector space $\mathcal{L}(SIM(2))$ of left-invariant vector-fields is given by
\begin{align}
\widetilde{Ad}(g)=d\mathcal{R}\circ Ad(g)\circ\omega.
\label{Evo_14}
\end{align}
The adjoint representation gives rise to the associated vector bundle $SIM(2)\times_{\widetilde{Ad}}\mathcal{L}(SIM(2))$. The corresponding connecting form on this vector bundle is given by
\begin{align}
\tilde{\omega}=\mathcal{A}_{2}\otimes \omega^{3}\wedge\omega^1 + \mathcal{A}_{3}\otimes \omega^{1}\wedge\omega^2 +
\mathcal{A}_{2}\otimes \omega^{4}\wedge\omega^2 +
\mathcal{A}_{3}\otimes \omega^{4}\wedge\omega^3.
\label{Evo_17}
\end{align}
Then $\tilde{\omega}$ yields the following $4\times 4$ matrix-valued $1$-form:
\begin{align}
\tilde{\omega}_{j}^{k}(\cdot):= -\tilde{\omega}(\omega^k,\cdot,\mathcal{A}_{j}), \ k,j\in\{1,2,3,4\}
\label{Evo_11}
\end{align}
on the frame bundle, where the sections are moving frames\footnote{See \cite[Ch.8]{Spivak1970} for more details on frame bundles and moving frames.}. Let $\{\mu_k\}_{k=1}^4$ denote the sections in the tangent bundle $E:=(SIM(2),T(SIM(2)))$ which coincides with the left-invariant vector fields $\{\mathcal{A}_k\}_{k=1}^4$. Then the matrix-valued $1$-form \eqref{Evo_10} yields the Cartan connection \footnote{Following the definitions in \cite{Spivak1970}, formally this is not a Cartan connection but a Koszul connection, see \cite[Ch.6]{Spivak1970} , corresponding to a Cartan connection, i.e.\ a associated differential operator corresponding to a Cartan connection. We avoid these technicalities and use ``Cartan connection" (a Koszul connection in \cite{Spivak1970}) and ``Cartan-Ehresmann connection form" (Ehresmann connection in \cite{Spivak1970}).
} $D$ on the tangent bundle $(SIM(2),T(SIM(2)))$ given by the covariant derivatives
\begin{align}
D_{X\lvert_{\gamma(t)}}(\mu(\gamma(t)))&:=
D(\mu(\gamma(t)))(X\lvert_{\gamma(t)})\nonumber\\
&=\sum\limits_{k=1}^{4}\dot{a}^{k}
\mu_{k}(\gamma(t))+
\sum\limits_{k=1}^{4}a^{k}(\gamma(t))
\sum\limits_{k=1}^{4}\tilde{\omega}^{j}_{k}
(X\lvert_{\gamma(t)})\mu_{j}(\gamma(t))
\label{Evo_12}\\
&=\sum\limits_{k=1}^{4}\dot{a}^{k}
\mu_{k}(\gamma(t))+
\sum\limits_{k=1}^{4}\dot{\gamma}^i(t)a^k
(\gamma(t))\Gamma_{ik}^{j}\mu_{j}(\gamma(t)),\nonumber
\end{align}
with $\dot{a}^k=\dot{\gamma}^i(t)(\mathcal{A}_{i}\lvert_{\gamma(t)}a^k)$, for all tangent vectors $X\lvert_{\gamma(t)}=
\dot{\gamma}^i(t)\mathcal{A}_{i}
\lvert_{\gamma(t)}$ along a curve $t\mapsto\gamma(t)\in SIM(2)$ and all sections $\mu(\gamma(t))=
\sum\limits_{k=1}^{4}a^k(\gamma(t))\mu_{k}(\gamma(t))$. The Christoffel symbols in 
\eqref{Evo_12} are constant $\Gamma_{ik}^{j}=-c_{ik}^{j}$, with $c_{ik}^{j}$ the structure constants of the Lie algebra $T_{e}(SIM(2))$.
\end{theorem}
The proof follows on the line of \cite{Duits2010a} and is given in \cite{Sharma2013}.
As seen in the theorem above, we define the notion of covariant derivatives independent of the metric $\mathcal{G}$ on $SIM(2)$, which is the underlying principle behind the Cartan connection. Although in principle these two entities need not be related, in Lemma~\ref{lem:MetComConn} we show that the connection induced above is metric compatible.
\begin{definition}
Let $(M,\mathcal{G})$ be a Riemannian 
manifold (or pseudo-Riemannian manifold) 
where $M$ and $\mathcal{G}$ denote the 
manifold and the metric defined on it 
respectively. Let $\nabla$ denote a 
connection on $(M,\mathcal{G})$. Then 
$\nabla$ is called metric compatible with 
respect to $\mathcal{G}$ if 
\begin{align}
\nabla_{Z}\mathcal{G}(X,Y)=\mathcal{G}
(\nabla_{Z}X,Y)+\mathcal{G}(X,\nabla_{Z}Y),
\end{align}
for all $X, \ Y, \ Z \in T(M)$.
\end{definition}  
\begin{lemma}
\label{lem:MetComConn}
The Cartan connection $D$ on $(SIM(2),T(SIM(2)))$ is metric compatible with respect to $\mathcal{G}_{\beta}:SIM(2)\times T(SIM(2)) \times T(SIM(2))\rightarrow\mathcal{C}$ on $SIM(2)$ defined in \eqref{Evo_18}.
\end{lemma}
\textbf{Proof}
We first note that $D_{\mathcal{A}_i} \mathcal{G}(\mathcal{A}_j,\mathcal{A}_k)=0, \ i\in\{1,2,3,4\}$ as the covariant derivative of a scalar field is the same as partial derivative. Here $\{\mathcal{A}_{i}\}_{i=1}^4$ denote the basis of the left invariant vector fields $\mathcal{L}(SIM(2))$. The following brief computation
\begin{align}
\mathcal{G}
(D_{\mathcal{A}_i}\mathcal{A}_j,\mathcal{A}_k) 
+\mathcal{G}(\mathcal{A}_j,D_{\mathcal{A}_i}
\mathcal{A}_k)=
\mathcal{G}
(c^{l}_{ij}\mathcal{A}_l,\mathcal{A}_k) 
+\mathcal{G}(\mathcal{A}_j,c^{m}_{ik}\mathcal{A}_m) 
=-c^k_{ij}-c^j_{ik},
\end{align}
where we have used the fact that $\Gamma_{ij}^k=c^k_{ji}$ along with non-zero structure constants $c_{12}^{3}=-c_{21}^{3}= -c_{13}^{2}=c_{31}^{2}=-c_{34}^{3}= c_{43}^{3}=c_{42}^{2}=-c_{24}^{2}=1$, leads to the result.
$\hfill\Box$\\

The next theorem relates the previous results on Cartan connections and covariant derivatives to the nonlinear diffusion schemes on $SIM(2)$.
\begin{theorem}\label{thm:DiffGeomInt}
covariant derivative of a co-vector 
field $\mathbf{a}$ on the manifold $(SIM(2),
\mathcal{G})$ is a $(0,2)$-tensor field with 
components 
$\nabla_{j}a_i=\mathcal{A}_ja_i-\Gamma^k
_{ij}a_k$, whereas the covariant derivative 
of a vector field $\mathbf{v}$ on $SIM(2)$ 
is a $(1,1)$-tensor field with components\footnote{We have made use of the notation $\nabla_j:= D_{\mathcal{A}_j}$, when imposing the Cartan connection.} 
$\nabla_{j'}v^i=\mathcal{A}_{j'}v^i+
\Gamma^i_{j'k'}v^{k'}$. The Christoffel symbols equal minus the structure constants of the Lie algebra $\mathcal{L}(SIM(2))$, i.e.\ $\Gamma_{ij}^k=-c_{ij}^k$. The Christoffel symbols are anti-symmetric as the underlying Cartan connection $D$ has  constant torsion. The left-invariant equations \eqref{Evo_0} with a diagonal diffusion tensor \eqref{Evo_1} can be rewritten in covariant derivatives as
\begin{align}
\begin{cases}
&\partial_sW(g,s)=\sum\limits_{i,j=1}^4
\mathcal{A}_{i}((D_{ij}(W))(g,s)\mathcal{A}_j W)(g,s)=
\sum\limits_{i,j=1}^4
\nabla_{i}((D_{ij}(W))(g,s)\nabla_j W)(g,s),\\
&W(g,0)=\mathcal{W}_{\psi}f(g), \text{ for all } g\in SIM(2), \ s>0.
\end{cases}
\label{Evo_20}
\end{align}
Both convection and diffusion in the left-invariant evolution equations \eqref{Evo_0}
take place along the exponential curves in $SIM(2)$ which are covariantly constant curves with respect to the Cartan connection. \end{theorem}
For proof see \cite{Sharma2013}.
\begin{remark}
 Though the connection $D$ is $\mathcal{G}$ compatible (Lemma \ref{lem:MetComConn}) , we have shown in the previous theorem that our connection is not torsion free, i.e.\  the torsion tensor $T(X,Y)\neq 0$ for all $X,Y\in T(SIM(2))$. As a result minimum distance curves in the group $SIM(2)$ are curves minimizing the induced Riemannian metric and do not coincide with ``straight curves" (auto-parallel curves) in the group $SIM(2)$. The auto-parallel curves are the exponential curves. A full analysis of the shortest distance curves in $SIM(2)$ is beyond the scope of this article.
\end{remark}

\section{Curvature estimation via best Exponential Curve fit}\label{App:CurvatureEstimation}
Curvature estimation of a spatial curve using $SE(2)$-OS  is based on the optimal exponential curve fit at each point. In \cite{Duits2010a,Franken2009} the authors suggest two methods for such best exponential curve fit. Below is a brief summary.
\begin{itemize}
\item{Compute the curvature of the projection $\mathbf{x}(s(t))=\mathbb{P} _{\mathbb{R}^2}(g_0\exp(t\sum_{i=1}^3c^i_*A_i))$ of the optimal exponential curve in $SE(2)$ on the ground plane from an eigenvector $\mathbf{c}_*=(c_*^\theta,c_*^\xi,c_*^\eta)$. This eigenvector of $(\tilde{H}_{\beta}|U|)^T(\tilde{H}_{\beta}|U|)$, with a $3\times 3$ Hessian
\begin{align}
\tilde{H}_{\beta}|U|=
\begin{pmatrix}
\beta^2\partial_\theta\partial_\theta |U| & \beta\partial_\xi\partial_\theta| U| & \beta\partial_\eta\partial_\theta |U| \\
 \partial_\theta\partial_\xi |U| & \partial_\xi\partial_\xi |U| & \partial_\eta\partial_\xi |U| \\
\partial_\theta\partial_\eta |U| & \partial_\xi\partial_\eta |U| & \partial_\eta\partial_\eta |U| 
\end{pmatrix},
\label{Evo_26}
\end{align}
corresponds to the smallest eigenvalue. The curvature estimation is given by,
\begin{align*}
\kappa_{est}=\|\ddot{x}(s)\|
sgn(\ddot{x}(s)\cdot \mathbf{e}_{\eta})=\frac{c_*^\theta sign(c_*^\xi)}{\sqrt{(c_*^\xi)^2+(c_*^\xi)^2}}.
\end{align*}
Note that unlike the $SIM(2)$ case curvature is constant in this case.
}
\item{
In this method the choice of optimal exponential curve is restricted to horizontal exponential curves, which are curves in the $(SE(2),\omega^3)$ sub-Riemannian manifold. The idea is to diffuse along horizontal curves because typically the mass of a $SE(2)$-OS is concentrated around a horizontal curve \cite{Franken2009} and therefore this is a fast curvature estimation method. 

Compute the curvature of the projection $\mathbf{x}(s(t))=\mathbb{P} _{\mathbb{R}^2}(g_0\exp(t\sum_{i=1}^3c^i_*A_i))$ of the optimal exponential curve in $SE(2)$ on the ground plane from the eigenvector $\mathbf{c}_*=(c_*^\theta,c_*^\eta)$. This eigenvector of $(\tilde{H}^{hor}_{\beta}|U|)^T(\tilde{H}^{hor}_{\beta}|U|)$, with a $3\times 2$ horizontal Hessian
\begin{align}
\tilde{H}^{hor}_{\beta}|U|=
\begin{pmatrix}
\beta^2\partial_\theta\partial_\theta |U| & \beta\partial_\xi\partial_\theta| U|  \\
 \partial_\theta\partial_\xi |U| & \partial_\xi\partial_\xi |U|  \\
\partial_\theta\partial_\eta |U| & \partial_\xi\partial_\eta |U|  
\end{pmatrix},
\end{align}
corresponds to the smallest eigenvalue. The curvature estimation given by,
\begin{align*}
\kappa^{hor}_{est}=\|\ddot{x}(s)\|
sgn(\ddot{x}(s)\cdot \mathbf{e}_{\eta})=\frac{c_*^\theta}{c_*^\xi}.
\end{align*}
For numerical experiments on these curvature estimates on orientation scores of noisy images, see \cite{Franken2009}.
}
\end{itemize}

\textbf{References}


{\small
\bibliographystyle{elsarticle-num}
\bibliography{lit_duits_sharma}

\begin{thebibliography}{10}
\expandafter\ifx\csname url\endcsname\relax
  \def\url#1{\texttt{#1}}\fi
\expandafter\ifx\csname urlprefix\endcsname\relax\def\urlprefix{URL }\fi
\expandafter\ifx\csname href\endcsname\relax
  \def\href#1#2{#2} \def\path#1{#1}\fi

\bibitem{Lindeberg1993}
T.~Lindeberg, Scale-Space Theory in Computer Vision, Kluwer international
  series in engineering and computer science: Robotics: Vision, manipulation
  and sensors, Springer, 1993.

\bibitem{Alvarez1993}
L.~Alvarez, F.~Guichard, P.-L. Lions, J.-M. Morel, Axioms and fundamental
  equations of image processing, Archive for rational mechanics and analysis
  123~(3) (1993) 199--257.

\bibitem{Romeny1997}
B.~ter Haar~Romeny, L.~Florack, J.~Koenderink, M.~Viergever, Scale-Space Theory
  in Computer Vision: First International Conference, Scale-Space'97, Utrecht,
  The Netherlands, July 2-4, 1997, Proceedings, Vol. 1252, Springer, 1997.

\bibitem{Lindeberg2013}
T.~Lindeberg, Generalized axiomatic scale-space theory, in: Advances in Imaging
  and Electron Physics, Vol 178, Advances in Imaging and Electron Physics,
  Elsevier, 2013, pp. 1--96.

\bibitem{Perona1990}
P.~Perona, J.~Malik, Scale-space and edge detection using anisotropic
  diffusion, Pattern Analysis and Machine Intelligence, IEEE Transactions on
  12~(7) (1990) 629 --639.
\newblock \href {http://dx.doi.org/10.1109/34.56205}
  {\path{doi:10.1109/34.56205}}.

\bibitem{Weickert1999}
J.~Weickert, Coherence-enhancing diffusion filtering, International Journal of
  Computer Vision 31~(2-3) (1999) 111--127.

\bibitem{Franken2009}
E.~Franken, R.~Duits, Crossing-preserving coherence-enhancing diffusion
  on invertible orientation scores, International Journal of Computer Vision 85
  (2009) 253--278.

\bibitem{Scharr2012}
H.~Scharr, K.~Krajsek, A short introduction to diffusion-like methods, in:
  Mathematical methods for signal and image analysis and representation,
  Computational imaging and vision, Springer London, 2012, Ch.~1.

\bibitem{Duits2007b}
R.~Duits, B.~Burgeth, Scale spaces on lie groups, in: Proceedings of the 1st
  international conference on Scale space and variational methods in computer
  vision, SSVM'07, Springer-Verlag, Berlin, Heidelberg, 2007, pp. 300--312.

\bibitem{Citti2006}
G.~Citti, A.~Sarti, A cortical based model of perceptual completion in the
  roto-translation space, Journal of Mathematical Imaging and Vision 24~(3)
  (2006) 307--326.

\bibitem{Duits2012}
R.~Duits, H.~F{\"u}hr, B.~Janssen, M.~Bruurmijn, L.~Florack, H.~van Assen,
  Evolution equations on gabor transforms and their applications, Applied and
  Computational Harmonic Analysis 35 (2013) 483--526.

\bibitem{Duits2005}
R.~Duits,
  \href{http://bmia.bmt.tue.nl/people/RDuits/THESISRDUITS.pdf}{Perceptual
  {O}rganization in {I}mage {A}nalysis}, Ph.D. thesis, Technische Universiteit
  Eindhoven (2005).
\newline\urlprefix\url{http://bmia.bmt.tue.nl/people/RDuits/THESISRDUITS.pdf}

\bibitem{DeAngelis1995}
G.~C. DeAngelis, I.~Ohzawa, R.~D. Freeman, Receptive-field dynamics in the
  central visual pathways, Trends in Neurosciences 18~(10) (1995) 451 -- 458.

\bibitem{Young1987}
R.~Young, The gaussian derivative model for spatial vision: I. retinal
  mechanisms, Spatial Vision 2~(4) (1987) 273--293.

\bibitem{Landy1991}
M.~Landy, J.~Movshon, Computational models of visual processing, Bradford book,
  Mit Press, 1991.

\bibitem{Bosking1997}
W.~Bosking, Y.~Zhang, B.~Schofield, D.~Fitzpatrick, Orientation selectivity and
  the arrangement of horizontal connections in tree shrew striate cortex, The
  Journal of Neuroscience 17(6) (1997) 2112--2127.

\bibitem{Duits2007a}
R.~Duits, M.~Felsberg, G.~Granlund, B.~Romeny, Image analysis and
  reconstruction using a wavelet transform constructed from a reducible
  representation of the euclidean motion group, Int. J. Comput. Vision 72~(1)
  (2007) 79--102.

\bibitem{Antoine1996}
J.-P. Antoine, R.~Murenzi, Two-dimensional directional wavelets and the
  scale-angle representation, Signal Processing 52~(3) (1996) 259 -- 281.

\bibitem{Antoine1999}
J.-P. Antoine, R.~Murenzi, P.~Vandergheynsta, Directional wavelets revisited:
  Cauchy wavelets and symmetry detection in patterns, Applied and Computational
  Harmonic Analysis 6~(3) (1999) 314 -- 345.

\bibitem{Antoine2004}
J.-P. Antoine, R.~Murenzi, P.~Vandergheynst, S.~Ali, Two-dimensional wavelets
  and their relatives, Cambridge University Press, 2004.

\bibitem{Jacques2011}
L.~Jacques, L.~Duval, C.~Chaux, G.~Peyr{\'e}, A panorama on multiscale
  geometric representations, intertwining spatial, directional and frequency
  selectivity, Signal Processing 91~(12) (2011) 2699--2730.

\bibitem{Ali2014}
S.~T. Ali, J.-P. Antoine, J.-P. Gazeau, Coherent States, Wavelets, and Their
  Generalizations, Springer New York, 2014.

\bibitem{Donoho1999}
D.~Donoho, E.~Cand{\`e}s, Ridgelets: a key to higher-dimensional
  intermittency?, Royal Society of London Philosophical Transactions Series A
  357 (1999) 2495.

\bibitem{Donoho2005}
D.~Donoho, E.~Cand{\`e}s, Continuous curvelet transform: I. resolution of the
  wavefront set, Applied and Computational Harmonic Analysis 19~(2) (2005) 162
  -- 197.

\bibitem{Donoho2005a}
D.~Donoho, E.~Cand{\`e}s, Continuous curvelet transform: {II}. discretization
  and frames, Applied and Computational Harmonic Analysis 19~(2) (2005) 198 --
  222.

\bibitem{Candes2006}
E.~Cand{\`e}s, L.~Demanet, D.~Donoho, L.~Ying, Fast discrete curvelet
  transforms, Multiscale Modeling and Simulation 5 (2006) 861--899.

\bibitem{Labate2005}
D.~Labate, W.-Q. Lim, G.~Kutyniok, G.~Weiss, Sparse multidimensional
  representation using shearlets, in: Optics \& Photonics 2005, International
  Society for Optics and Photonics, 2005, pp. 59140U--59140U.

\bibitem{Guo2007}
K.~Guo, D.~Labate, Optimally sparse multidimensional representation using
  shearlets, SIAM journal on mathematical analysis 39~(1) (2007) 298--318.

\bibitem{Dahlke2011}
S.~Dahlke, G.~Steidl, G.~Teschke, Shearlet coorbit spaces: compactly supported
  analyzing shearlets, traces and embeddings, Journal of Fourier Analysis and
  Applications 17~(6) (2011) 1232--1255.

\bibitem{Bodmann2013}
B.~G. Bodmann, G.~Kutyniok, X.~Zhuang, Gabor shearlets, arXiv preprint
  (arXiv:1303.6556) (2013).

\bibitem{Ali1998}
S.~T. Ali, A general theorem on square-integrability: Vector coherent states,
  J. Math. Phys. 39(8) (1998) 3954--3964.

\bibitem{Hubel1993}
D.~H. Hubel, Evolution of ideas on the primary visual cortex, 1955-1978: A
  biased historical account, Physiology or Medicine Literature Peace Economic
  Sciences, Nobel Prize Lectures (1993) 24.

\bibitem{Grossmann1985}
A.~Grossmann, J.~Morlet, T.~Paul, Transforms associated to square integrable
  group representations. {I}. {G}eneral results, J. Math. Phys. 26~(10) (1985)
  2473--2479.

\bibitem{Fuhr2005}
H.~F{\"u}hr, Abstract Harmonic Analysis of Continuous Wavelet Transforms, no.
  1863 in Lecture Notes in Mathematics, Springer, 2005.

\bibitem{Kalitzin1999}
S.~Kalitzin, B.~Romeny, M.~Viergever, Invertible apertured orientation filters
  in image analysis, International Journal of Computer Vision 31~(2-3) (1999)
  145--158.

\bibitem{Bekkers2012}
E.~Bekkers, R.~Duits, T.~Berendschot, B.~ter Haar~Romeny, A multi-orientation
  analysis approach to retinal vessel tracking, Journal of Mathematical Imaging
  and Vision (2012) 1--28.

\bibitem{Duits2008}
R.~Duits, M.~Almsick, The explicit solutions of linear left-invariant second
  order stochastic evolution equations on the {$2$D} {E}uclidean motion group,
  Quarterly of Applied Mathematics 66 (2008) 27--67.

\bibitem{Duits2010}
R.~Duits, E.~Franken, Left-invariant parabolic evolutions on {$SE(2)$} and
  contour enhancement via invertible orientation scores part {I}: Linear
  left-invariant diffusion equations on {$SE(2)$}, Quarterly of Applied
  Mathematics 68 (2010) 255--292.

\bibitem{Duits2010a}
R.~Duits, E.~Franken, Left-invariant parabolic evolutions on {$SE(2)$} and
  contour enhancement via invertible orientation scores part {II}: Nonlinear
  left-invariant diffusions on invertible orientation scores, Quarterly of
  Applied Mathematics 68 (2010) 293--331.

\bibitem{Martens1988}
F.~Martens,
  \href{http://alexandria.tue.nl/extra3/proefschrift/PRF6A/8810117.pdf}{Spaces
  of analytic functions on inductive/projective limits of {H}ilbert {S}paces},
  Ph.D. thesis, Technische Universiteit Eindhoven (1988).
\newline\urlprefix\url{http://alexandria.tue.nl/extra3/proefschrift/PRF6A/8810117.pdf}

\bibitem{Borel1991}
A.~Borel, H.~Bass, Linear algebraic groups, Vol. 126, Springer-Verlag New York,
  1991.

\bibitem{Boscain2013}
U.~Boscain, J.-P. Gauthier, R.~Chertovskih, A.~Remizov, Hypoelliptic diffusion
  and human vision, arXiv preprint (arXiv:1304.2062) (2013).

\bibitem{Bankman2008}
I.~Bankman, Handbook of Medical Image Processing and Analysis, Academic Press
  Series in Biomedical Engineering, Elsevier/Academic Press, 2008.

\bibitem{Zhang1997}
J.~Zhang, H.~Huang, Automatic background recognition and removal (abrr) in
  computed radiography images, Medical Imaging, IEEE Transactions on 16~(6)
  (1997) 762 --771.
\newblock \href {http://dx.doi.org/10.1109/42.650873}
  {\path{doi:10.1109/42.650873}}.

\bibitem{Jacques2007}
L.~Jacques, J.-P. Antoine, Multiselective pyramidal decomposition of images:
  wavelets with adaptive angular selectivity, International Journal of
  Wavelets, Multiresolution and Information Processing 5~(05) (2007) 785--814.

\bibitem{Unser1999}
M.~Unser, Splines: A perfect fit for signal/image processing, IEEE Signal
  Processing Magazine 16 (1999) 22--38.

\bibitem{Dieudonne1977}
J.~Dieudonn{\'e}, Treatise on Analysis, no. v. 5 in Pure and applied
  mathematics, Academic Press, 1977.

\bibitem{Aubin2001}
T.~Aubin, A course in differential geometry, Graduate studies in mathematics,
  American Mathematical Society, 2001.

\bibitem{Hormander1967}
L.~H{\"o}rmander, Hypoellptic second order differential equations, Acta
  Mathematica 119 (1967) 147--171.

\bibitem{Duits2011}
R.~Duits, E.~Franken, Left-invariant diffusions on the space of positions and
  orientations and their application to crossing-preserving smoothing of hardi
  images, Int. J. Comput. Vision 92~(3) (2011) 231--264.

\bibitem{Hsu2002}
E.~Hsu, Stochastic analysis on manifolds, Contemporary Mathematics, American
  Mathematical Society, 2002.

\bibitem{Nagel1990}
A.~Nagel, F.~Ricci, E.~M. Stein, Fundamental solutions and harmonic analysis on
  nilpotent groups, Bulletin of American Mathematical Society 23 (1990)
  139--144.

\bibitem{Elst1998}
A.~ter Elst, D.~Robinson, Weighted subcoercive operators on {L}ie groups,
  Journal of Functional Analysis 157~(1) (1998) 88 -- 163.

\bibitem{Mallat2012}
S.~Mallat, Group invariant scattering, Communications on Pure and Applied
  Mathematics 65~(10) (2012) 1331--1398.

\bibitem{Sifre2013}
L.~Sifre, S.~Mallat, E.~N.~S. DI, Rotation, scaling and deformation invariant
  scattering for texture discrimination, in: Proc. CVPR, 2013.

\bibitem{Franken2008}
E.~Franken, \href{http://alexandria.tue.nl/extra2/200910002.pdf}{Enhancement of
  crossing elongated structures in images}, Ph.D. thesis, Technische
  Universiteit Eindhoven (2008).
\newline\urlprefix\url{http://alexandria.tue.nl/extra2/200910002.pdf}

\bibitem{Patton2006}
N.~Patton, T.~M. Aslam, T.~MacGillivray, I.~J. Deary, B.~Dhillon, R.~H.
  Eikelboom, K.~Yogesan, I.~J. Constable, Retinal image analysis: concepts,
  applications and potential, Progress in retinal and eye research 25~(1)
  (2006) 99--127.

\bibitem{Budai2009}
A.~Budai, J.~Hornegger, G.~Michelson, Multiscale approach for blood vessel
  segmentation on retinal fundus images, Investigative Ophtalmology and Visual
  Science 50~(5) (2009) 325.

\bibitem{Abramoff2010}
M.~D. Abr{\`a}moff, M.~K. Garvin, M.~Sonka, Retinal imaging and image analysis,
  Biomedical Engineering, IEEE Reviews in 3 (2010) 169--208.

\bibitem{Xu2011}
X.~Xu, M.~Niemeijer, Q.~Song, M.~Sonka, M.~K. Garvin, J.~M. Reinhardt, M.~D.
  Abramoff, Vessel boundary delineation on fundus images using graph-based
  approach, Medical Imaging, IEEE Transactions on 30~(6) (2011) 1184--1191.

\bibitem{Bankhead2012}
P.~Bankhead, C.~N. Scholfield, J.~G. McGeown, T.~M. Curtis, Fast retinal vessel
  detection and measurement using wavelets and edge location refinement, PloS
  one 7~(3) (2012) e32435.

\bibitem{Frangi1998}
A.~F. Frangi, W.~J. Niessen, K.~L. Vincken, M.~A. Viergever, Multiscale vessel
  enhancement filtering, in: Medical Image Computing and Computer-Assisted
  Interventation— MICCAI'98, Springer, 1998, pp. 130--137.

\bibitem{Hannink2014}
J.~Hannink, R.~Duits, E.~Bekkers,
  \href{http://arxiv.org/abs/1402.4963}{Vesselness via multiple scale
  orientation scores}, arXiv preprint (2014).
\newline\urlprefix\url{http://arxiv.org/abs/1402.4963}

\bibitem{Seo2007}
H.~J. Seo, P.~Chatterjee, H.~Takeda, P.~Milanfar, A comparison of some state of
  the art image denoising methods, in: Signals, Systems and Computers, 2007.
  ACSSC 2007. Conference Record of the Forty-First Asilomar Conference on,
  IEEE, 2007, pp. 518--522.

\bibitem{Buades2005}
A.~Buades, B.~Coll, J.-M. Morel, A review of image denoising algorithms, with a
  new one, Multiscale Modeling \& Simulation 4~(2) (2005) 490--530.

\bibitem{Takeda2007}
H.~Takeda, S.~Farsiu, P.~Milanfar, Kernel regression for image processing and
  reconstruction, Image Processing, IEEE Transactions on 16~(2) (2007)
  349--366.

\bibitem{Rubbens2009}
M.~Rubbens, A.~Mol, R.~Boerboom, R.~Bank, F.~Baaijens, C.~Bouten, Intermittent
  straining accelerates the development of tissue properties in engineered
  heart valve tissue, Tissue Engineering: Part A 15 (5).

\bibitem{Sharma2013}
U.~Sharma, R.~Duits,
  \href{http://www.win.tue.nl/analysis/reports/rana13-17.pdf}{Left-invariant
  evolutions of wavelet transforms on the similitude group}, CASA Report 13-17,
  Technische Universiteit Eidnhoven (June 2013).
\newline\urlprefix\url{http://www.win.tue.nl/analysis/reports/rana13-17.pdf}

\bibitem{Spivak1970}
M.~Spivak, A Comprehensive Introduction to Differential Geometry, Vol-{II},
  {3}rd ed., Publish or Perish Inc., 1970.

\end{thebibliography}
}






\end{document}